\documentclass{amsart}
\usepackage{amssymb}
\usepackage{mathrsfs}
\usepackage[utf8]{inputenc}
\usepackage[dvipsnames]{xcolor}
\usepackage{tikz}
\usepackage{tikz-cd}
\usetikzlibrary{arrows,decorations.pathmorphing,decorations.markings,decorations.pathreplacing,calc,matrix,intersections}
\tikzset{->-/.style={decoration={markings,mark=at position #1 with {\arrow{>}}},postaction={decorate}}}
\usepackage[normalem]{ulem}
\newcommand{\midarrow}{\draw[postaction={decorate}]}

\usepackage{fullpage}
\usepackage{enumitem}
\usepackage{braket}
\usepackage{accents}
\usepackage{textcomp}
\usepackage[scr=boondox]{mathalfa}
 
 \usepackage{mathtools}

\newcommand{\newword}[1]{\textbf{\emph{#1}}}

\newtheorem{thm}{Theorem}[section]
\newtheorem{prop}[thm]{Proposition}
\newtheorem{lem}[thm]{Lemma}

\newtheorem{lemma}[thm]{Lemma}
\newtheorem{cor}[thm]{Corollary}

\theoremstyle{definition}
\newtheorem{Def}[thm]{Definition}

\newtheorem{question}[thm]{Question}
\newtheorem{ex}[thm]{Example}
\theoremstyle{remark}

\newtheorem{rem}[thm]{Remark}
\newtheorem{caution}[thm]{Caution}
\newtheorem{obs}[thm]{Observation}

\newcommand{\trop}{\mathrm{trop}}
\DeclareMathOperator{\Cone}{\mathrm{Cone}}
\DeclareMathOperator{\NCl}{\mathrm{NCl}}
\DeclareMathOperator{\Int}{\mathrm{Int}}

\newcommand{\dgm}{\mathrm{dg}}
\newcommand{\pure}{\mathrm{pure}}
\newcommand{\MG}{\mathcal{MG}}

\newcommand{\Xg}{V_g} 
\newcommand{\Sg}{S_g} 
\newcommand{\Xgn}{V_{g,n}} 
\newcommand{\Sgn}{S_{g,n}} 
\newcommand{\X}{V} 
\newcommand{\Mgnh}{\M\mathcal{H}_{g,n}} 
\newcommand{\Mbargnh}{\overline{\M\mathcal{H}}_{g,n}}
\newcommand{\Mgnmh}{h\T(V_{g,n})} 
\newcommand{\Mbargnmh}{\overline{h\T}(V_{g,n})} 

\newcommand{\Sph}{\mathscr{S}} 

\DeclareMathOperator{\len}{len}

\DeclareMathOperator{\PSL}{PSL}
\DeclareMathOperator{\SL}{SL}
\DeclareMathOperator{\Aut}{Aut}
\DeclareMathOperator{\Isot}{Isot}
\DeclareMathOperator{\Spec}{Spec}

\DeclareMathOperator{\MCG}{MCG}
\DeclareMathOperator{\hMCG}{hMCG}
\DeclareMathOperator{\Link}{Lk}
\DeclareMathOperator{\Tw}{Tw}
\DeclareMathOperator{\tw}{tw}
\DeclareMathOperator{\Out}{Out}
\DeclareMathOperator{\id}{id}
\DeclareMathOperator{\CV}{\mathcal{CV}}

\DeclareMathOperator{\Stab}{Stab}
\newcommand{\Z}{\mathbb{Z}}
\newcommand{\C}{\mathbb{C}}

\newcommand{\Q}{\mathbb{Q}}

\newcommand{\D}{\mathcal{D}}

\newcommand{\R}{\mathbb{R}}
\newcommand{\U}{\mathcal{U}}
\renewcommand{\P}{\mathbb{P}}
\newcommand{\an}{\mathrm{an}}
\newcommand{\na}{\mathrm{na}}
\newcommand{\univ}{\mathrm{univ}}
\DeclareMathOperator{\Sch}{Sch}
\newcommand{\M}{\mathcal M}

\newcommand{\abs}[1]{\left\lvert#1\right\rvert}

\newcommand{\Mbar}{\overline{\M}}
\newcommand{\T}{\mathcal{T}}
\newcommand{\Tbar}{\overline{\mathcal{T}}}

\newcommand{\into}{\hookrightarrow}

\renewcommand{\bar}[1]{\overline{#1}}
\newcommand{\sslash}{{/\!\!/}}
\newcommand{\TV}{\mathcal T(V)}
\newcommand{\TbarV}{\overline{\mathcal T}(V)}
\renewcommand{\Im}{\mathrm{Im}}

\usepackage{hyperref}

\title{Handlebodies, Outer space, and tropical geometry}
\author[Ramadas]{Rohini Ramadas}
\author[Silversmith]{Rob Silversmith}
\author[Vogtmann]{Karen Vogtmann}
\author[Winarski]{Rebecca R. Winarski}


\begin{document}
\begin{abstract}
    The moduli space of graphs $\M_{g,n}^{\trop}$ is a polyhedral object that mimics the behavior of the moduli spaces $\M_{g,n}$, $\Mbar_{g,n}$ of (stable) Riemann surfaces; this relationship has been made precise in several different ways, which collectively identify $\M_{g,n}^{\trop}$ as the ``tropicalization'' of $\M_{g,n}$. We describe how this relationship lifts to some objects that live over $\M_{g,n}$ (like Teichm\"uller space) and that live over $\M_{g,n}^{\trop}$ (like the Culler-Vogtmann space $\CV_{g,n}^*$). We introduce the notion of a stable complex handlebody, and show that $\CV_{g,n}^*$ can be viewed as the tropicalization of a certain complex manifold $\Mgnmh$ that parametrizes complex handlebodies. An important ingredient is our construction of a partial compactification $\Mbargnmh\supset\Mgnmh$, which we prove is  a simply connected complex manifold with simple normal crossings boundary. When $n=0$, $\Mgnmh$ coincides with the moduli space of Schottky groups, $\Mbargnmh$ coincides with Gerritzen-Herrlich's extended Schottky space, and $\CV_{g,0}^*$ is the simplicial completion of the original Outer space. The resulting picture fits together many familiar objects from geometric group theory and surface topology, including Harvey’s curve complex, mapping class groups of surfaces and handlebodies, and augmented Teichm\"uller space. Many of the relationships between the objects that we see in this picture already exist in the literature, but we add some new ones, and generalize several existing relationships to include a number $n>0$ of punctures/leaves.
\end{abstract}

\maketitle

\section{Introduction}
    The purpose of this paper is to tell a story relating some well-studied objects in:
    \begin{itemize}
        \item Geometric group theory (automorphisms of free groups, Culler-Vogtmann spaces),
        \item Low-dimensional topology (surfaces and handlebodies, the structure of their loops and disks, and their associated mapping class groups),
        \item Complex geometry (Teichm\"uller space, Schottky space, and adjacent constructions),
        \item Algebraic geometry
        (Moduli spaces of stable curves), and
        \item Tropical/non-Archimedean geometry (Moduli spaces of tropical curves, and their non-Archimedean counterparts).
    \end{itemize}
    In the following sections, we outline the main characters in the story, and their relationships, which fit together to form the diagram in Figure \ref{fig:MainDiagram}. Many of the parts of this story already appear scattered in the literature around automorphisms of free groups, mapping class groups, handlebodies, Teichm\"uller space, Kleinian groups, and even phylogenetic trees. Our goal in this paper is to present what we hope is a  clean,  self-contained version of the story. Along the way, we will prove some new results, generalize some existing results, and give substantively new interpretations to some previously-studied objects. (See Section \ref{sec:IntroResults}.)
    
    \subsection{Boundary complexes and the meaning of ``tropicalization''}\label{sec:IntroTrop} Over the last decade or so, there has been a lot of research on tropical moduli spaces and the precise senses in which they are tropicalizations of algebraic moduli spaces. The classic example of this is $\M_{g,n}^{\trop}$, a generalized cone complex of dimension $3g-3+n$ parametrizing tropical curves/metric graphs. This is the tropicalization of the algebraic moduli space $\M_{g,n}$, a quasiprojective variety of (complex) dimension $3g-3+n$ parametrizing Riemann surfaces (i.e. smooth projective irreducible algebraic curves) of genus $g$ with $n$ labeled points, in the following senses:
        \begin{itemize}
        \item In the genus-zero case, $\M_{0,n}$ embeds as a closed subvariety of the algebraic torus $(\C^*)^{\binom{n}{2}-n}$, and Speyer-Sturmfels \cite{SpeyerSturmfels2004} used this fact to construct a tropicalization $\M_{0,n}^{\trop}$, a simplicial cone complex with a natural embedding in $\R^{\binom{n}{2}-n}$. As Speyer and Sturmfels noted, $\M_{0,n}^{\trop}$ is isomorphic to the space of phylogenetic trees studied in \cite{BilleraHolmesVogtmann2001}.
        \item Deligne and Mumford \cite{DeligneMumford1969} defined a normal crossings compactification $\Mbar_{g,n}$ of $\M_{g,n}$. The boundary  of 
        $\M_{g,n}$ as a subspace of $\Mbar_{g,n}$ is the complement $\Mbar_{g,n}\setminus\M_{g,n}$. The associated \newword{boundary complex} is a symmetric $\Delta$-complex\footnote{A $\Delta$-complex is a cell complex whose cells are homeomorphic to simplices.  A symmetric $\Delta$-complex is similar, but a ``cell" may be the quotient of a simplex by a finite group action.} that encodes the pattern of intersections of the strata of $\Mbar_{g,n}\setminus\M_{g,n}$. (See Section \ref{sec:BoundaryComplexes}.) The boundary complex is canonically identified with the \newword{link} $\Link(\M_{g,n}^{\trop})$, the moduli space of tropical curves with total edge length 1. (See e.g. \cite{Caporaso2012}.)
        \item Abramovich-Caporaso-Payne \cite{AbramovichCaporasoPayne2012} approached tropicalization via non-Archimedean geometry, using techniques of Thuillier on toroidal embeddings \cite{Thuillier2007} to show that the compactification $\Mbar_{g,n}$ of $\M_{g,n}$ determines a skeleton of the Berkovich analytification of $\M_{g,n}$, and this skeleton is canonically identified with $\M_{g,n}^{\trop}$. (Any complex variety of dimension $N$ has an associated Berkovich analytification, a topological space that can be obtained by taking an inverse limit of dimension-$N$ cone complexes, and a skeleton is a strong deformation retract of this topological space.) Cavalieri-Chan-Ulirsch-Wise \cite{CavalieriChanUlirschWise2020} later built the necessary tools for realizing this skeleton as a ``tropical moduli stack''.
        \item More recently, Chan-Galatius-Payne \cite{ChanGalatiusPayne2021} defined a surjective proper continuous map $\lambda_{g,n}:\M_{g,n}\to\M_{g,n}^{\trop}$, canonically defined up to homotopy, which induces an isomorphism from the compactly supported rational cohomology of $\M_{g,n}^{\trop}$ to the weight-zero compactly supported rational cohomology of $\M_{g,n}$.
    \end{itemize}
    Outer space and its variations, introduced in the following section, look a lot like tropical moduli spaces.  This was pointed out by Mikhalkin \cite[Sec. 3.1]{Mikhalkin2007}, and has been further explored and developed by several authors (\cite{ ChanMeloViviani2013,  PoineauTurchetti2020, ulirsch2021non}, summarized just below in Sections \ref{sec:IntroTeichmuller} and \ref{intro:Schottky}). In this paper, we use the above relationships between $\M_{g,n}$ and $\M_{g,n}^{\trop}$ as a ``wish list'' with which to explore the relationship between Outer space and tropical geometry. As we will see, this requires leaving the realm of algebraic geometry.  We remark that it has generally been challenging to talk about tropicalization outside of algebraic geometry, but some of the above points have compelling analogues.
    
    \subsection{Outer space and its relationship to \texorpdfstring{$\M_{g,n}^{\trop}$}{Mgntrop}}\label{sec:IntroOuterSpace} In the early 1980s, M. Culler and the third author introduced  \newword{Outer space} $\CV_g$ for $g\geq 2$; this is a $(3g-4)$-dimensional contractible space on which the group $\Out(F_g)$ of outer automorphisms of a finitely-generated free group acts with finite stabilizers \cite{CullerVogtmann1986}. Outer space $\CV_g$ parametrizes metric graphs together with a \newword{marking} which identifies the fundamental group of the graph with $F_g$, and $\Out(F_g)$ acts by changing the marking. The space $\CV_g$
    decomposes into a disjoint union of open simplices, and can  be completed to a simplicial complex $\CV_g^*$, called the \newword{simplicial completion} of Outer space. The action of $\Out(F_g)$ on $\CV_g$ extends to a simplicial action on $\CV_g^*$. For $n\geq 1$, Hatcher introduced generalizations $\CV_{g,n}$ and $\CV_{g,n}^*$ of Outer space and its simplicial completion, which both have natural actions of $A_{g,n}:=\Aut(F_g)\rtimes F_g^{n-1}$. To unify later statements we set $\CV_g=\CV_{g,0}$ and $A_{g,0}=\Out(F_g)$.

    The quotient $\CV_g/\Out(F_g)$ is the \newword{moduli space $\mathcal{MG}_g$ of metric graphs} of genus $g$ (i.e. connected with Betti number $g$) whose edge lengths sum to $1$. The quotient $\MG_{g,n}:=\CV_{g,n}/A_{g,n}$ is the moduli space of metric graphs whose edge lengths sum to 1, with $n$ labeled leaves. The quotients $\CV_g^{*}/\Out(F_g)$ and $\CV_{g,n}^*/A_{g,n}$ are naturally identified with $\Link(\M_{g}^{\trop})$ and $\Link(\M_{g,n}^{\trop})$ respectively. Under these identifications, $\mathcal{MG}_g$ and $\MG_{g,n}$ are realized as open subsets of $\Link(\M_{g}^{\trop})$ and $\Link(\M_{g,n}^{\trop})$ respectively.

    The above spaces have been used extensively to study the group $\Out(F_g)$ and its cohomology. We define them rigorously in Section \ref{sec:IntroOuterSpace}. In Appendix \ref{app:SimplicialCompletion} we give two alternate descriptions of  $\CV_{g,n}^*$.

    \begin{rem}
        $\Link(\M_{g,n}^{\trop})$ is a \emph{finite} cell complex. In contrast,  $\CV_{g,n}^*$ is an \emph{infinite} simplicial complex (unless $g=0$ or $(g,n)=(1,1)$). Therefore $\CV_{g,n}^*$ cannot be the boundary complex of a normal crossings compactification of an algebraic variety --- our first hint that we must leave the realm of algebraic geometry in order to interpret $\CV_{g,n}^*$ as a tropical moduli space.     \end{rem}

     \subsection{Teichm\"uller space \texorpdfstring{$\T(S_{g,n})$}{T(Sgn)}: its relationship to \texorpdfstring{$\M_{g,n}$}{Mgn}, and the tropical perspective}\label{sec:IntroTeichmuller}
     Fix an $n$-pointed genus-$g$ oriented surface $S_{g,n}$; its Teichm\"uller space $\T(S_{g,n})$ is a complex manifold parametrizing $n$-pointed genus-$g$ \emph{Riemann} surfaces together with a marking given by a homeomorphism from $S_{g,n}$ (See Section \ref{sec:ModuliTeichmuller}). The mapping class group  $\MCG(S_{g,n})$, i.e.  the group of isotopy classes of orientation-preserving homeomorphisms of $S_{g,n}$ that fix the $n$ labeled points,  acts on $\T(S_{g,n})$ by changing the marking, and $\M_{g,n}$ is the quotient of $\T_{g,n}$ by this action.  
     Since the action has finite stabilizers, $\M_{g,n}$ is an orbifold, $\T_{g,n}$ is its orbifold universal cover, and $\MCG(S_{g,n})$ is the orbifold fundamental group of $\M_{g,n}.$  
     
    Outer space, with its action by $\Out(F_g),$ was originally introduced by Culler and Vogtmann as an analogue of Teichm\"uller space $\T_g$ with its action of $\MCG(\Sg)$. In particular, $\M_{g}$ parametrizes Riemann surfaces while $\T(S_{g})$ parametrizes Riemann surfaces with a topological marking; $\mathcal{MG}_g$ parametrizes metric graphs while $\CV_{g}$ parametrizes metric graphs with a topological marking. 
    
    This analogy has been developed further: The period map from Teichm\"uller space to Siegel space $\mathbb{H}$   has a natural analog mapping $\CV_g$ to the space of positive-definite quadratic forms (see \cite{Baker2011}, where this map is called the {\em Jacobian map}). This map was extended (under certain conditions) by Chan-Melo-Viviani \cite{ChanMeloViviani2013}, who defined a ``tropical Siegel space'' $\mathbb{H}^{\trop}$ and  a ``tropical period map'' map  from the cone over $\CV_g^*$ to $\mathbb{H}^{\trop}$. (Note: In \cite{ChanMeloViviani2013}, the cone over $\CV_g^*$ is referred to as ``tropical Teichm\"uller space'', which is at odds with the terminology introduced just below.) 

     In \cite{Bers1974} Bers introduced a partial compactification $\Tbar(S_{g,n})$ of $\T(S_{g,n})$, called \newword{augmented Teichm\"uller space}. The action of $\MCG(S_{g,n})$ extends to $\Tbar(S_{g,n})$, and the quotient is naturally identified with $\Mbar_{g,n}$. The space $\Tbar(S_{g,n})$ is not a complex manifold, but the boundary $\Tbar(S_{g,n})\setminus\T(S_{g,n})$ has a natural stratification by locally closed subsets, each of which is a connected complex manifold. The boundary complex of $\Tbar(S_{g,n})$ (appropriately defined) is the \newword{curve complex} $\mathcal{C}(S_{g,n})$, an infinite simplicial complex  indexing multicurves on $S_{g,n}$ (see Section \ref{sec:Complexes}). The quotient of the cone $\Cone(\mathcal{C}(S_{g,n}))$ over $\mathcal{C}(S_{g,n})$ by $\MCG(S_{g,n})$ can be identified with $\M_{g,n}^{\trop}$. For these reasons, the first author \cite[Sec. 3.4]{ramadas2024thurston} proposed $\Cone(\mathcal{C}(S_{g,n}))$ as an alternate candidate for ``tropical Teichm\"uller space'', suggesting that the relationship between $\mathcal{C}(S_{g,n})$ and $\T(S_{g,n})$ is more in line with properties of tropicalizations listed in Section \ref{sec:IntroTrop} than the relationship between $\CV_{g,n}^*$ and $\T(S_{g,n})$, and noting in support of this that the Chan-Galatius-Payne map $\lambda_{g,n}:\M_{g,n}\to\M_{g,n}^{\trop}$ naturally lifts to an $\MCG(S_{g,n})$-equivariant continuous map $\widetilde{\lambda}_{g,n}:\T(S_{g,n})\to\Cone(\mathcal{C}(S_{g,n}))$. (See Section \ref{sec:CGP}.)

    \subsection{Schottky space and its relationship with \texorpdfstring{$\CV_g^*$}{CVg*}}\label{intro:Schottky}
    Despite the \emph{analogy} discussed above between $\T(S_g)$ and $\CV_g$, the \emph{relationship} between $\T(S_g)$ and $\CV_g$ does not fit into any of the frameworks of ``tropicalization'' listed in Section \ref{sec:IntroTrop}. However, there is a space, called (marked) Schottky space $\Sch_g$, which lives in between $\T(S_{g})$ and $\M_{g}$, for which $\CV_g^*$ does act in some sense as a tropicalization. $\Sch_g$ is a complex manifold that parametrizes Schottky groups (a special class of free subgroups of $\PSL_2(\C)$) together with a choice of free generators. It can be obtained as the quotient of $\T(\Sg)$ by a certain infinite-index subgroup of $\MCG(S_g)$; it admits an infinite-degree covering map from $\T(\Sg)$ and an infinite-degree covering map to $\M_g$. $\Sch_g$ is known to be related to $\CV_{g}^{*}$ in two ways:
    \begin{enumerate}
        \item \textbf{Via a partial compactification.} There is a partial compactification $\Sch_g\subseteq\bar{\Sch}_g$ whose boundary complex is $\CV_g^*$. The space $\bar{\Sch}_g$, a   complex manifold (but not an algebraic variety), was constructed by Gerritzen and Herrlich \cite{GerritzenHerrlich1988} and is called \newword{extended Schottky space}. This is the avenue that we explore further in this article.
        \item \textbf{Via non-Archimedean geometry.} Gerritzen and Herrlich \cite{Gerritzen1983, Herrlich2006} also constructed a non-Archimedean analog $\Sch_g^{\na}$ of Schottky space and noted that $\Sch_g^{\na}$ admits a natural ``tropicalization'' map to the cone over $\CV_g$. Ulirsch introduced a variant/enlargement  $\Sch_g^{\na,*}$ of $\Sch_g^{\na}$ and showed that it contains the cone over $\CV_g^*$ as a skeleton (strong deformation retract) \cite{ulirsch2021non}. $\Sch_g^{\na, *}$ parametrizes genus-$g$ \textit{non-Archimedean curves} together with a topological marking by a genus-$g$ graph, and therefore Ulirsch refers to $\Sch_g^{\na, *}$ as ``non-Archimedean Teichm\"uller space'', continuing the analogy of Chan-Melo-Viviani. However, Ulirsch also mentions the connection to Schottky space, speculating that there must be meaningful sense in which $\Sch_g^{\na}$ is related to the classical Schottky space $\Sch_g$. A precise connection between $\Sch_g^{\na}$ and $\Sch_g$ was established by Poineau-Turchetti \cite{PoineauTurchetti2020}, who constructed a ``universal" analytic Schottky space $\Sch_g^{\an,\univ}$, defined over the integers, of which both the non-Archimedean $\Sch_g^{\na}$ and the classical $\Sch_g$ can be realized as specializations\footnote{More precisely, there is a topological space $M(\Z)$ (with additional structure) that parametrizes multiplicative seminorms on the ring $\Z$ that are bounded from above by the standard Archimedean norm on $\Z$. The universal Schottky space $\Sch_g^{\an,\univ}$ admits a map $\pi$ to $M(\Z)$; this can be thought of as a ``universal family" of Schottky spaces. The inclusion $\Z\into \C$ can be realized in multiple ways as an inclusion of normed rings; two of these in particular are obtained by (1) the Archimedean norm on both rings and (2) the trivial norm on both rings. Each of these inclusions of normed rings induces a map from a one-point space into $M(\Z)$. One can pull back the universal family $\pi:\Sch_g^{\an,\univ}\to M(\Z)$ along both of these maps; along (1) we obtain the classical Schottky space $\Sch_g$, while along (2) we obtain the non-Archimedean Schottky space $\Sch_g^{\na}$. This can be recast slightly more concretely as follows: the construction in \cite{PoineauTurchetti2020} can be used to construct a bundle of Schottky spaces over an open disk, most fibers of which are isomorphic to classical Schottky space, but whose fiber over the origin is isomorphic to non-Archimedean Schottky space. In this sense we can think of classical Schottky space as ``degenerating" to non-Archimedean Schottky space. The construction in \cite{PoineauTurchetti2020} was used by Dang and Mehmeti  \cite{DangMehmeti2024} to relate the dynamics of non-Archimedean Schottky groups to those of classical Schottky groups.}.
 \end{enumerate}
It is well-known that a Schottky group determines a hyperbolic \emph{handlebody}, and that therefore $\Sch_g$ can be thought of as a space of handlebodies --- see just below.

    \subsection{Main results}\label{sec:IntroResults}
\begin{figure}
        \centering
        \begin{tikzcd}\label{eq:BigTropicalizationDiagram}
    \boxed{\T(\Sgn)\subset\Tbar(\Sgn)}\arrow[ddddd,"\MCG(\Sgn)",swap]\arrow[rrrrrrr,bend left=5,rightsquigarrow,"\trop"]&&&&&&&\mathcal{C}(\Sgn)\arrow[ddddd,"\MCG(\Sgn)"]\\
    \\
    &\boxed{\T(\Xgn)\subset\Tbar(\Xgn)}\arrow[uul,hook]\arrow[d,"\Tw(\Xgn)"]\arrow[dd,bend right=75,"\MCG(\Xgn)",swap]\arrow[rrrr,bend left=5,rightsquigarrow,"\trop"]&&&&\mathcal{D}(\Xgn)\arrow[uurr,hook]\arrow[d,"\Tw(\Xgn)",swap]\arrow[dd,bend left=40,"\MCG(\Xgn)"]&&\\
    &\boxed{\Mgnmh\subset\Mbargnmh}\arrow[d,"A_{g,n}"]\arrow[rrrr,bend left=5,rightsquigarrow,"\trop"]&&&&\CV_{g,n}^*\arrow[d,"A_{g,n}",swap]&&\\
    &\boxed{\Mgnh\subset\Mbargnh}\arrow[dl]\arrow[rrrr,bend left=5,rightsquigarrow,"\trop"]&&&&\Link\M_{g,n}^{\trop}\arrow[drr,equal]&&\\
    \boxed{\M_{g,n}\subset\Mbar_{g,n}}\arrow[rrrrrrr,bend left=5,rightsquigarrow,"\trop"]&&&&&&&\Link\M_{g,n}^{\trop}\\
\end{tikzcd}
        \caption{Summary of the main objects and results of the paper. Arrows labeled by a group are quotients. Squiggly arrows labeled ``trop'' satisfy various subsets of the notions of tropicalization listed in Section \ref{sec:IntroTrop}, see Section \ref{sec:Revisit}.}
        \label{fig:MainDiagram}
    \end{figure}
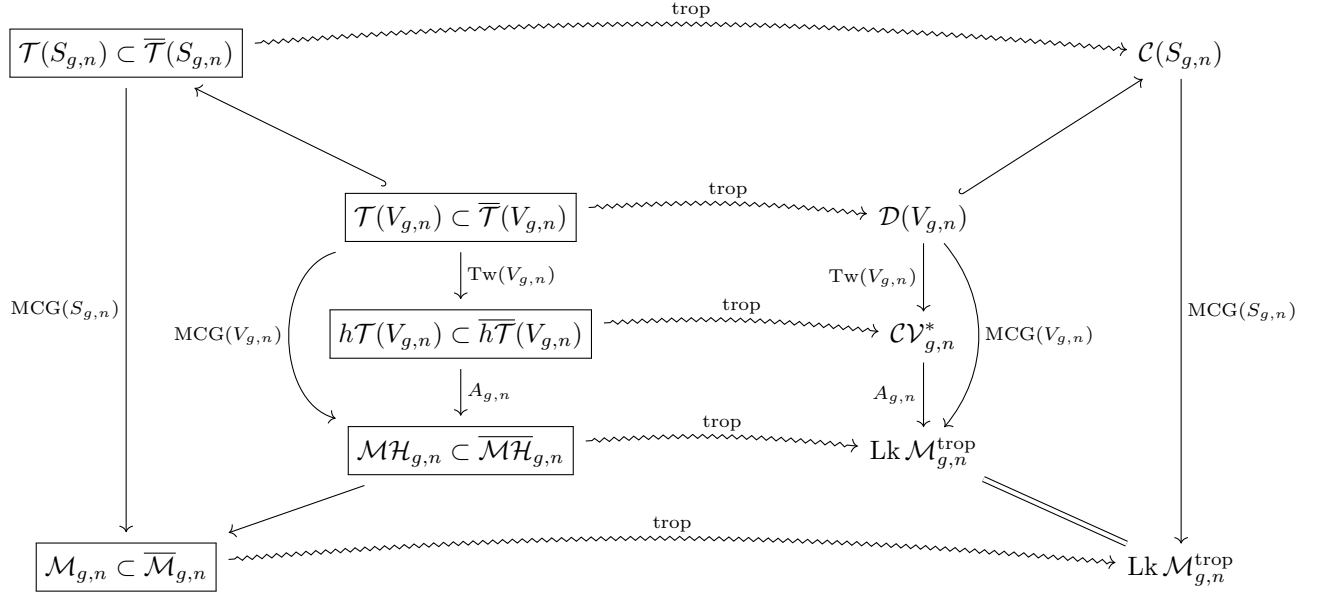

    In this paper, we aim to explore and connect the above points, to reframe them slightly, and to generalize them to allow labeled points (i.e. the case $n>0$). Schematically, from Sections \ref{sec:IntroTrop} and \ref{sec:IntroTeichmuller} we have a diagram:
    \begin{align}\label{eq:OriginalTropicalizationDiagram}
\begin{tikzcd}[ampersand replacement=\&]
    \boxed{\T(\Sgn)\subset\Tbar(\Sgn)}\arrow[d,"\MCG(\Sgn)",swap]\arrow[rr,bend left=5,rightsquigarrow,"\trop"]\&\&\mathcal{C}(\Sgn)\arrow[d,"\MCG(\Sgn)"]\\
    \boxed{\M_{g,n}\subset\Mbar_{g,n}}\arrow[rr,bend left=5,rightsquigarrow,"\trop"]\&\&\Link\M_{g,n}^{\trop}
\end{tikzcd}
\end{align}
    We fit this diagram into the larger diagram shown in Figure \ref{fig:MainDiagram}, whose pieces we now explain. A genus-$g$ \newword{handlebody} is an oriented 3-manifold with boundary that is obtained by gluing $2g$ disks on the boundary of a 3-ball to each other in pairs (via orientation-reversing homeomorphisms). Let $\Xgn$ be an $n$-pointed genus-$g$ handlebody, that is, a genus-$g$ handlebody together with $n$ distinct labeled points on its boundary. The boundary $S_{g,n}:=\partial\Xgn$ of $\Xgn$ is an oriented $n$-pointed genus-$g$ surface. 

    Associated to $\Xgn$ there is an infinite simplicial complex $\mathcal{D}(\Xgn)$ called the \emph{disk complex}, which lives naturally as a sub-complex of the curve complex $\mathcal{C}(S_{g,n})$. The mapping class group $\MCG(\Xgn)$ acts on $\mathcal{D}(\Xgn)$, and we are especially interested in the restriction of this action to the twist subgroup $\Tw(\Xgn)\subseteq\MCG(\Xgn)$, which is generated by Dehn twists around simple disks   (see Sections \ref{sec:HandlebodiesBackground} and \ref{sec:MappingClassGroups}).
    By a classical result of Luft  \cite{Luft1978}, when $n=0$ there is an exact sequence of groups, with natural maps:
    \begin{equation}
        \begin{tikzcd}
            1\arrow[r]&\Tw(\Xg)\arrow[r]&\MCG(\Xg)\arrow[r]&\Out(\pi_1(\Xg))\arrow[r]&1.
        \end{tikzcd}
    \end{equation}
    Note that $\Out(\pi_1(\Xg))\cong\Out(F_g)$. Our first result is that this sequence (appropriately modified) remains exact for arbitrary $n$:
    
    \medskip
    
    \noindent\textbf{Theorem A} (Thm. \ref{thm:LuftMarkedPoints}).  \textit{There is a natural exact sequence:
        \begin{equation}
        \begin{tikzcd}[ampersand replacement = \&]
            1\arrow[r]\&\Tw(\Xgn)\arrow[r]\&\MCG(\Xgn)\arrow[r]\&A_{g,n}\arrow[r]\&1.
        \end{tikzcd}
    \end{equation}
    (Precisely, for the third map to be canonical we should choose an identification of $A_{g,n}$ with the homotopy mapping class group of $\Xgn$, see Section \ref{sec:MappingClassGroups}).}
    
    \medskip

    We next introduce the following spaces:
    \begin{itemize}
        \item $\T(\Xgn)$ is the Teichm\"uller space of $\Xgn$ (defined in Section \ref{sec:T(Vgn)}), a complex manifold of dimension $3g-3+n$. It parametrizes ``complex handlebodies'' with a topological marking, i.e. a homeomorphism from the reference handlebody $\Xgn$. $\T(\Xgn)$ is naturally isomorphic to $\T(S_{g,n})$. 
        \item $\Tbar(\Xgn)$ is the augmented Teichm\"uller space of $\Xgn$ (defined in Section \ref{sec:Tbar(Vgn)}), a partial compactification of $\T(\Xgn)$. It parametrizes \textit{stable} complex handlebodies with a topological marking, i.e. a sufficiently well-behaved \emph{surjection} from the reference handlebody $\Xgn$. $\Tbar(\Xgn)$ is \emph{not} isomorphic to $\Tbar(\Sgn)$, but rather is naturally isomorphic to an open subset of $\Tbar(\Sgn)$. Like $\Tbar(\Sgn)$, $\Tbar(\Xgn)$ is not a complex manifold, but the boundary $\Tbar(\Xgn)\setminus\T(\Xgn)$ is stratified by locally closed subsets, each of which is a connected complex manifold. The boundary complex of $\Tbar(\Xgn)$ (appropriately defined) is the disk complex $\mathcal{D}(\Xgn)$. 
           \item   $\Mgnmh$ is the space of $n$-pointed genus-$g$ complex handlebodies marked by a \emph{homotopy equivalence} from $\Xgn$ (defined in Section \ref{sec:Mgmhbar}). It is a complex manifold of dimension $3g-3+n$, which is isomorphic to (marked) Schottky space $\Sch_g$ when $n=0$. Here ``$h\mathcal{T}$'' is short for ``homotopy Teichm\"uller space''.
        \item $\Mbargnmh$  is the space of $n$-pointed genus-$g$ stable complex handlebodies marked by a homotopy equivalence from $\Xgn$ (defined in Section \ref{sec:Mgmhbar}). This is a partial compactification of $\Mgnmh$. It is isomorphic to extended Schottky space $\bar{\Sch}_g$ when $n=0$ \cite{GerritzenHerrlich1988,Herrlich2006}.
             \item  $\Mgnh$ is the moduli space of (unmarked) $n$-pointed genus-$g$ complex handlebodies (defined in Section \ref{sec:Mghbar}). It is a complex orbifold  of dimension $3g-3+n$. In the case $n=0$, we have a natural isomorphism $\M(V_{g,0})\cong\Sch_g/\Out(F_g)$. (The latter quotient is called ``unmarked Schottky space''.) Here ``$\mathcal{MH}$'' is short for ``moduli space of handlebodies.''
        \item $\Mbargnh$ is the moduli space of (unmarked) $n$-pointed genus-$g$ stable complex handlebodies (defined in Section \ref{sec:Mghbar}), a partial compactification of $\Mgnh$. In the case $n=0$, there is a natural isomorphism $\overline{\mathcal{MH}}_{g,0}\cong\overline{\Sch}_g/\Out(F_g)$ \cite{GerritzenHerrlich1988,Herrlich2006}.

    \end{itemize}

    \begin{question}
    In the setting $n>0,$ is there a useful Schottky-group-like interpretation for $\Mgnh$, perhaps parametrizing Schottky groups together with $n$ distinct marked orbits on $\P^1\setminus(\text{limit set})$?
\end{question}
    
    \begin{rem}\label{rem:GenusZero}
    In the case that $g=0$, the three groups $\MCG(S_{0,n})$, $\MCG(\X_{0,n})$, and $\Tw(\X_{0,n})$ all coincide, as do the three spaces $\M_{0,n}$, $h\T(V_{0,n})$ and $\mathcal{MH}_{0,n}$, and the three spaces $\Mbar_{0,n}$, $\overline{h\T}(V_{0,n})$ and $\overline{\mathcal{MH}}_{0,n}$.
\end{rem}
Our results are summarized below, and some of them are depicted schematically in Figure \ref{fig:MainDiagram}:

\medskip

\noindent\textbf{Theorem B} (Thm. \ref{thm:hTBarSimplyConnected}). $\Mbargnmh$ is simply connected for all $g,n$.

\begin{rem}
    $\Mbargnmh$ is not contractible --- this is clear when $g=0$ since $\Mbargnmh\cong\Mbar_{0,n}.$ It is also not contractible when $g>0$ (see Remark \ref{rem:HTbarNotContractible}), since it contains strata isomorphic to the product of copies of $\Mbar_{0,n}$.
\end{rem}

    \medskip

    \noindent\textbf{Theorem C} (Complex Structures).
    \begin{enumerate}
        \item (Thm. \ref{thm:MarkedHandlebodiesComplexManifold}) \textit{$\Mbargnmh$ has the natural structure of a complex manifold of dimension $3g-3+n$, with respect to which the boundary $\Mbargnmh\setminus\Mgnmh$ is a countable union of complex submanifolds meeting at simple normal crossings.}
        \item (Thm. \ref{thm:UnmarkedHandlebodiesComplexOrbifold}) \textit{$\Mbargnh$ has the natural structure of a complex orbifold of dimension $3g-3+n$, with respect to which the boundary $\Mbargnh\setminus\Mgnh$ is a countable union of complex suborbifolds meeting at normal crossings.}
        \item (Thms. \ref{thm:MarkedHandlebodiesComplexManifold}(\ref{it:MgmhMghCovering},\ref{it:OutFgActionBiholomorphisms}), and \ref{thm:UnmarkedHandlebodiesComplexOrbifold}\eqref{it:MghMgEtale}) \textit{All maps in Figure \ref{fig:MainDiagram} between complex manifolds/orbifolds are holomorphic and \'etale, and all group actions on complex manifolds/orbifolds in Figure \ref{fig:MainDiagram} are actions by biholomorphisms.}
    \end{enumerate}
    \begin{rem}
        When $n=0$, Theorem C is \cite[Prop. 9]{GerritzenHerrlich1988}.
    \end{rem}
    We caution that $\Tbar(\Xgn)$ does not have a complex structure at the boundary, so its maps to $\Mbargnmh$, $\Mbargnh$ and $\Mbar_{g,n}$ are not holomorphic at the boundary. The map from $\Mbargnmh$ to $\Mbargnh$ is an orbifold covering map. The maps from $\Mbargnmh$ to $\Mbar_{g,n}$ and from $\Mbargnh$ to $\Mbar_{g,n}$ are local biholomorphisms (in the orbifold sense), but are not orbifold covering maps over the boundary of $\Mbar_{g,n}$. (The reader should have in mind the exponential map $\C\to\C$, which is a local biholomorphism but cannot be extended to a covering map.)

    \medskip

    \noindent\textbf{Theorem D.}
    \begin{enumerate}
        \item (Prop. \ref{prop:UnmarkedBijection}) \textit{The (holomorphic) map $\T(\Xgn)\to\Mgnh$ and the (continuous) map $\Tbar(\Xgn)\to\Mbargnh$ are canonically identified with the quotient maps by the mapping class group $\MCG(\Xgn)$.}
        \item (Prop. \ref{prop:MarkedBijection}) \textit{The (holomorphic) map $\T(\Xgn)\to\Mgnmh$ and the (continuous) map $\Tbar(\Xgn)\to\Mbargnmh$ are canonically identified with the quotient maps by the twist subgroup $\Tw(\Xgn)$.}
    \end{enumerate}

    \begin{rem}\label{rem:OpenMaps}
        The continuous maps $\Tbar(\Xgn)\to\Mbargnh$ and $\Tbar(\Xgn)\to\Mbargnmh$ are open, since they are quotient maps by actions of discrete groups.
    \end{rem}

    \medskip

    \noindent\textbf{Theorem E.} (Thm. \ref{thm:QuotientOfDiskComplex}) \textit{The quotient $\mathcal{D}(\Xgn)/\Tw(\Xgn)$ is canonically identified with $\CV_{g,n}^*$.}

    \medskip

    \noindent\textbf{Corollary F.}
    \begin{enumerate}
        \item (Cor. \ref{cor:BoundaryMghLinkMgTrop}) \textit{The boundary complex of $\Mgnh\subseteq\Mbargnh$ is canonically identified with $\Link(\M_{g,n}^{\trop})$.}
        \item (Cor. \ref{cor:MgmhBoundaryComplexCV}) \textit{The boundary complex of $\Mgnmh\subseteq\Mbargnmh$ is canonically identified with $\CV_{g,n}^*$.}
    \end{enumerate}

    \begin{rem}
        Surprisingly, the map $\Mbargnh\to\Mbar_{g,n}$ induces an isomorphism on boundary complexes, even though it is infinite-to-1.
    \end{rem}
    
    \medskip

    In analogy with the Chan-Galatius-Payne map $\tilde\lambda_{\Sgn}:\T(\Sgn)\to\Cone(\mathcal{C}(\Sgn))$ (see Section \ref{sec:IntroTeichmuller}), we have:

    \medskip
    
    \noindent\textbf{Theorem G.} (Prop. \ref{prop:CGPAnalogy}) \textit{There is a continuous $\MCG(\Xgn)$-equivariant surjection $$\tilde\lambda_{\Xgn}:\T(\Xgn)\to\Cone(\mathcal{D}(\Xgn))$$ which descends to an $A_{g,n}$-equivariant continuous surjection $$\lambda_{g,n}^{h\mathcal{T}}:\Mgnmh\to\Cone(\CV_{g,n}^*).$$ In turn, $\lambda_{g,n}^{h\mathcal{T}}$ descends to a continuous surjection $$\lambda_{g,n}^{\mathcal{MH}}:\Mgnh\to\M_{g,n}^{\trop}.$$}

    \medskip

    By work of Hatcher \cite{Hatcher1995}, $\CV_{g,n}^*$ also has an interpretation in terms of the sphere complex $\Sph(W_{g,n})$ of a closed 3-manifold $W_{g,n}$ called a \emph{doubled handlebody}, see Appendix \ref{app:SimplicialCompletion}. There is a natural ``doubling map'' from the disk complex $\D(\Xgn)$ to the sphere complex $\Sph(W_{g,n})$. As a consequence of Theorem D and some results from 3-manifold theory, we have:

    \medskip

    \noindent\textbf{Proposition H.} (Prop. \ref{prop:DiskComplexSphereComplex}) \textit{The doubling map  induces an isomorphism from the quotient of $\D(\Xgn)$ by $\Tw(\Xgn)$ to $\Sph(W_{g,n})$.}

    \medskip

    \subsection{Revisiting the term ``tropicalization"}\label{sec:Revisit} In Figure \ref{fig:MainDiagram}, we have five squiggly arrows labeled by ``trop". Each of these five arrows admits the interpretation ``boundary complex of (partial) compactification". However, working from the bottom to the top of the diagram, the use of the term ``tropicalization" fits into three decreasing levels of niceness:
\begin{enumerate}
    \item Nicest: $(\M_{g,n}\subset\Mbar_{g,n})\leadsto \Link\M_{g,n}^{\trop}$. 
    \begin{enumerate}
        \item $\M_{g,n}$ and $\Mbar_{g,n}$ are both \emph{algebraic stacks} as well as complex orbifolds.
        \item Since $\M_{g,n}$ is algebraic, is not compact, and is an orbifold rather than a complex manifold, there are holomorphic maps from punctured disks to $\M_{g,n}$ that do not extend across the puncture to the whole disk. This could be for one of three reasons --- either (1) the map has an essential singularity at the puncture and therefore has no limiting value, or (2) there is a limiting value at the puncture, but that limit point is in $\Mbar_{g,n}\setminus\M_{g,n}$, or (3) there is a limiting value in $\M_{g,n}$, but for the map extends across the puncture only after taking a finite cover of the punctured disc.
        \item $\Mbar_{g,n}$ is algebraic and compact, and its algebraic/complex structure extends that of the open subset $\M_{g,n}$. Therefore, every holomorphic map from a punctured disk to $\Mbar_{g,n}$ that is meromorphic at the puncture extends across the puncture to a map to $\Mbar_{g,n}$ after taking a finite cover of the punctured disk.
        \item Also, since $\Mbar_{g,n}$ is compact, there is no possibility for further enlargement to a complex orbifold in which $\Mbar_{g,n}$ is dense and open. 
        \item $\M_{g,n}$, being algebraic, has an associated non-Archimedean space, its Berkovich analytification $\M_{g,n}^{\na}$. Roughly speaking, a point of $\M_{g,n}^{\na}$ represents an ``infinitesimal one-parameter family\footnote{More precisely, a point of $\M_{g,n}^{\na}$ represents a field extension $K$ of $\C$ with a non-Archimedean valuation, together with a morphism of $\C$-schemes from $\Spec(K)$ to $\M_{g,n}$.}" of points in $\M_{g,n}$. In particular, a meromorphic map from a punctured disk to $\M_{g,n}$ determines a point of $\M_{g,n}^{\na}$. $\M_{g,n}^{\trop}$ lives inside $\M_{g,n}^{\na}$ as a skeleton, i.e. strong deformation retract. This means that a meromorphic map $\alpha$ from a punctured disk to $\M_{g,n}$ determines a point $\alpha^{\trop}$ of $\M_{g,n}^{\trop}$. The cone of $\M_{g,n}^{\trop}$ that $\alpha^{\trop}$ lies in corresponds to the boundary stratum of $\Mbar_{g,n}$ to which $\alpha$ sends the puncture.  
        \item There is the Chan-Galatius-Payne map $\lambda_{g,n}:\M_{g,n}\to\M_{g,n}^{\trop}$.
    \end{enumerate}

    \item Less nice, still pretty good: $(\Mgnh\subset\Mbargnh)\leadsto \Link\M_{g,n}^{\trop}$ and $(\Mgnmh\subset\Mbargnmh)\leadsto \CV_{g,n}^{*}$.
    \begin{enumerate}
        \item $\Mgnmh$ and $\Mbargnmh$ are complex manifolds, and $\Mgnh$ and $\Mbargnh$ are complex orbifolds (Theorem C). None of the four spaces is algebraic.
        \item $\Mgnmh$ and $\Mgnh$ are not compact. They are nonalgebraic, and it is not clear \emph{a priori} that there exist holomorphic maps from punctured disks to $\Mgnmh$ and $\Mgnh$ that do not extend across the puncture. However:
        \item The complex structures on $\Mbargnmh$ and $\Mbargnh$ extend those, respectively, on $\Mgnmh$ and $\Mgnh$. From this, we conclude that there do indeed exist holomorphic maps from punctured disks to $\Mgnmh$ and $\Mgnh$ that do not extend across the puncture but that do extend instead to maps from the disk to $\Mbargnmh$ and $\Mbargnh$. 
        \item $\Mbargnmh$ and $\Mbargnh$ are not compact. However, we \textbf{speculate} that there are no further enlargements to complex manifolds/orbifolds in which $\Mbargnmh$ and $\Mbargnh$ are dense and open (and such that the complement is a complex analytic subspace).
        \item Given that $\Mgnmh$ and $\Mgnh$ are nonalgebraic, they do not admit Berkovich analytifications the way that $\M_{g,n}$ does (as far as we are aware). However, in the case $n=0$, $h\T(V_{g,0})$ does admit a non-Archimedean \emph{analogue} $h\T(V_{g,0})^{\na}:=\Sch_{g}^{\na,*}$ \cite{Gerritzen1983, Herrlich2006, ulirsch2021non}. By \cite{ulirsch2021non}, $\Cone(\CV_g^{*})$ does live inside $h\T(V_{g,0})^{\na}$ as a skeleton. By \cite{PoineauTurchetti2020} there is an explicit relationship between  $\Sch_g^{\na}$ (which is a subset of $\Sch_{g}^{\na,*}=h\T(V_{g,0})^{\na}$) and $h\T(V_{g,0})$, but this relationship is not ``is the Berkovich analytification of"\footnote{We also wonder whether there is a universal enlargement $\Sch_g^{\an,\univ,*}$ of $\Sch_g^{\an,\univ}$ over $M(\Z)$ whose Archimedean fiber is unchanged but whose non-Archimedean fiber over the trivial norm is $\Sch_{g}^{\na,*}$, rather than $\Sch_{g}^{\na}$.
        
        Separately, Ulirsch \cite{ulirsch2021non} points out a way to add labeled points to the construction of $\Sch_g^{\na,*}$. However we are not certain whether the resulting space bears the correct relationship to $\Mgnmh$ in the above senses, e.g. it is not clear whether it contains $\CV_{g,n}^*$ as a skeleton. Nor has the construction in \cite{PoineauTurchetti2020} been carried out when labeled points are allowed. However, we would not be surprised if there is a well-behaved version of the story in this generality.}. As such, we are not aware of any interpretation of the form ``points of $h\T(V_{g,0})^{\na}$ correspond to one-parameter degenerations in $\Mgnmh$", but we would not be surprised if such an interpretation could be made meaningful/precise. Implicit in Theorem C is the following: Any holomorphic map from a punctured disk to $h\T(\Xgn)$ that extends to a map from the disk to $\bar{h\T}(\Xgn)$ determines a point of $\Cone(\CV_{g,n}^{*})$. 

        \item There are Chan-Galatius-Payne-style maps $\lambda_{g,n}^{h\mathcal{T}}:\Mgnmh\to \Cone(\CV_{g,n}^{*})$ and $\lambda_{g,n}^{\mathcal{MH}}:\Mgnh\to \M_{g,n}^{\trop}$.
    \end{enumerate}

    \item Least satisfying: $(\T(\Xgn)\subset\Tbar(\Xgn))\leadsto \mathcal{D}(\Xgn)$ and $(\T(S_{g,n})\subset\Tbar(S_{g,n}))\leadsto \mathcal{C}(S_{g,n})$. 

    \begin{enumerate}
        \item $\T(\Xgn)$ and $\T(S_{g,n})$ are (isomorphic) complex manifolds. Neither $\Tbar(\Xgn)$ nor $\Tbar(S_{g,n})$ admits a globally defined complex structure, but both are stratified by locally closed subsets, each of which is a complex manifold. None of the four spaces is algebraic.
        \item $\T(\Xgn)$ and $\T(\Sgn)$ are not compact. However, every holomorphic map from a punctured disk to $\T(\Xgn)$ or $\T(\Sgn)$ extends across the puncture. (To see this, recall that $\T(\Xgn)$ and $\T(S_{g,n})$ are both biholomorphic to the unit ball, and that a bounded holomorphic function on a punctured disk extends across the puncture.)
        \item The complex structures on $\T(\Xgn)$ and $\T(S_{g,n})$ do not extend to complex structures on $\Tbar(\Xgn)$ and $\Tbar(S_{g,n})$ respectively. Indeed, there is no complex manifold $X$ that contains $\T(S_{g,n})$ as a dense open subset and such that $X\setminus \T(S_{g,n})$ is a complex analytic subspace; if it did, then there would necessarily exist a holomorphic map from a punctured disk to $\T(S_{g,n})$ that does not extend across the puncture.
        \item $\Tbar(\Xgn)$ and $\Tbar(\Sgn)$ are also not compact.
        \item Given that $\T(S_{g,n})$ and $\T(\Xgn)$ are non-algebraic, they do not admit Berkovich analytifications. Nor are there non-Archimedean counterparts known to us -- the spaces that have been termed ``tropical Teichm\"uller spaces" in the literature have precise relationships to Schottky rather than to Teichm\"uller spaces. Given that every holomorphic map from a punctured disk to $\T(\Xgn)$ or to $\T(S_{g,n})$ extends across the puncture, from the perspective of one-parameter complex degenerations, we would expect the tropicalizations of both $\T(\Xgn)$ and of $\T(S_{g,n})$ to be one-point sets. 
        \item There are Chan-Galatius-Payne-style maps $\tilde\lambda_{g,n}:\T(S_{g,n})\to\Cone(\mathcal{C}(S_{g,n}))$ and $\tilde\lambda_{\Xgn}:\T(\Xgn)\to\Cone(\mathcal{D}(\Xgn))$. 
        
    \end{enumerate}
\end{enumerate}

\subsection{Relation to work of Hainaut-Petersen and Petersen-Wade} Since $\T(\Xgn)$ is contractible and has a properly discontinuous action of $\MCG(\Xgn),$ the quotient $\T(\Xgn)/\MCG(\Xgn)\cong\Mgnh$ is an orbifold classifying space for $\MCG(\Xgn)$. In \cite{HainautPetersen2025} (see also \cite[Sec. 3.2]{PetersenWade2024}), Hainaut and Petersen constructed a different orbifold classifying space $\mathcal{HM}_{g,n}$ for $\MCG(\Xgn)$. The space $\mathcal{HM}_{g,n}$ is an orbifold \emph{inside} $\M_{g,n}$; it is a small punctured neighborhood of a certain union of strata in $\Mbar_{g,n}$. They used $\mathcal{HM}_{g,n}$, together with results of Chan-Galatius-Payne \cite{ChanGalatiusPayne2021}, to relate the cohomology of $\MCG(\Xgn),$ to the cohomology of graph complexes.

Later, Petersen-Wade \cite{PetersenWade2024} used a ``truncated'' version of $\mathcal{HM}_{g,n}$ to prove that $\MCG(\Xgn)$ is a \emph{virtual duality group} in the sense of \cite{BieriEckmann1973}.

    \subsection{Acknowledgements} We are very grateful to Saul Schleimer for in-depth conversations around the proofs of Theorem \ref{thm:LuftMarkedPoints} and Proposition \ref{prop:DiskComplexSphereComplex}. We are grateful to Dick Hain for useful conversations and ideas around simply-connectedness and noncontractibility of $\Mbargnmh.$ We are also grateful to Nguyen-Bac Dang, Purvi Gupta, Rob Kropholler, Vlere Mehmeti, Daniele Turchetti, and Ric Wade for useful discussions.

\section{Conventions and terminology}

For the rest of the paper, we assume $(g,n)$ satisfies $2g+n\ge 3.$

\subsection{Graphs}\label{sec:GraphConventions}
All graphs in this paper will be finite and connected, possibly with loops, multiple edges and/or leaves.   A  convenient way to specify such a graph $G$ is to give a set $H$ (the half-edges of $G$),  a partition $P$ of $H$ (whose pieces correspond to the vertices of $G$) and an involution $i\colon H\to H$ (whose order 2 orbits are the edges of $G$ and fixed points are the leaves of $G$). If a graph has $n\geq 1$ leaves, they will be labeled $\{p_1,\ldots,p_n\}$.  The \newword{valence} of a partition piece $U$ (corresponding to a vertex of $G$) is the number of half-edges  in $U$. 

Two graphs $\{H_1,P_1,i_1\}$ and $\{H_2,P_2,i_2\}$ are \newword{isomorphic} if and only if there is a bijection $H_1\to H_2$ sending $P_1$ to $P_2$ and $i_1$ to $i_2$.

 \begin{ex}The $n$-thorned rose $R_{g,n}$ is the graph with one vertex, $g$ loops and $n$ labeled leaves.  Thus $H$ has $2g+n$ elements, $P$ is the trivial partition and $i$ has  $g$ orbits of order $2$ and $n$ fixed points.  
 \end{ex}

A \newword{metric graph} is a graph together with a positive real length assigned to each edge.  At times we will normalize so that the sum of the edge-lengths is equal to 1. Note that leaves do not have lengths.  An isomorphism of metric graphs is a graph isomorphism that preserves edge-lengths.

A \newword{weighted graph} is a graph together with a nonnegative integer ``weight" at each vertex.  The \newword{genus} of a weighted graph $G$ is the sum of all vertex weights, plus the first Betti number $b_1(G).$    
A weighted graph is called \newword{stable} if every weight-zero vertex has valence $\ge3$ and every weight-one vertex has valence $\ge1.$  An isomorphism of weighted graphs is a graph isomorphism that preserves weights.

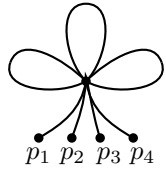
\begin{figure}
\begin{center}
\begin{tikzpicture} [thick, scale=.25]
\newcommand\vertex[1]{\fill #1 circle (.25)}
\vertex{(0,0)};
\draw  (0,0) to [out=135, in=180]  (0,4.1);
\draw  (0,0) to [out=45, in=0]  (0,4.1); 
\draw   (0,0) to [out=45, in=110]  (4,1);
\draw [thick] (0,0) to [out=-30, in=-70] (4,1);
\draw   (0,0) to [out=135, in=70]  (-4,1);
\draw   (0,0) to [out=210, in=-110] (-4,1);
\draw (0,0) to [out=260, in = 30] (-2.5,-3);
\vertex{(-2.5,-3)};
\node [below] (P1) at (-2.5,-3) {$p_1$};
\draw (0,0) to [out=265, in = 70] (-.75,-3);
\vertex{(-.75,-3)};
\node [below] (P2) at (-.7,-3) {$p_2$};
\draw (0,0) to [out=-265, in = 110] (.75,-3);
\vertex{(.75,-3)};
\node [below] (P3) at (1.25,-3) {$p_3$};
\draw (0,0) to [out=-80, in = 150] (2.5,-3);
\vertex{(2.5,-3)};
\node [below] (P4) at (3,-3) {$p_4$};

\end{tikzpicture}
\end{center}
\caption{The topological realization of $R_{3,4}$.}\label{fig:rose}
\end{figure}

The above combinatorial description of a graph is useful for efficiently defining a graph and distinguishing between its edges and its leaves, but we will sometimes want to consider a graph to be a topological space. We do this by forming its \newword{topological realization}.  This is a 1-dimensional cell complex with a 0-cell for each vertex and a 1-cell for each edge.  If the graph has leaves, then in order to make a legitimate cell complex   we must add an additional 0-cell labeled $p_i$ ``at the end of"  each leaf labeled $p_i$.    A picture of the topological realization of $R_{3,4}$ is shown in Figure \ref{fig:rose}. 

We define a homotopy equivalence of graphs to be a homotopy equivalence of their topological realizations.  If the graphs have labeled leaves, the homotopy equivalence must send the $i$-th leaf  to the $i$-th leaf.

 \subsection{Markings: a disambiguation}\label{sec:Markings} In this article, we deal with the following different types of \newword{markings} of geometric objects:
\begin{enumerate}
    \item A marking of a Riemann surface is a homeomorphism from a reference topological surface, up to isotopy. Note that this is the same data as a homotopy equivalence from the reference surface, up to homotopy, by \cite[Thm. 8.9]{FarbMargalit2011}. This is introduced in Section \ref{sec:ModuliTeichmuller}
    \item A marking of a  metric graph is a homotopy equivalence from a reference  graph, up to homotopy. For convenience, the reference graph is picked to be a thorned rose. This is introduced in Section \ref{sec:OuterSpaceBackground}.
    \item A marking of a complex handlebody is a \emph{homeomorphism} from a reference topological handlebody, up to \emph{isotopy}. This is introduced in Section \ref{sec:T(Vgn)}.
    \item A homotopy-marking of a genus-$g$ complex handlebody is a \emph{homotopy equivalence} from a reference handlebody, up to \emph{homotopy}. 
    Note that this is the same data as a homotopy equivalence from a genus-$g$ rose graph, up to homotopy. Note also that a homotopy-marking of a complex handlebody does not uniquely determine a marking of that handlebody. This is introduced in Section \ref{sec:Mgmhbar}.
\end{enumerate}
\begin{caution}
    In the literature it is common to use the term ``marked point" for a labeled point on a surface.  To avoid confusion with markings as defined above, we will use the term ``labeled point'' when we need to emphasize the labeling.
\end{caution}

\subsection{\texorpdfstring{$n$}{n}-pointed maps}\label{sec:nPointedMaps} We will often talk about various kinds of $n$-pointed topological spaces $(X,p_1,\ldots,p_n)$, where $p_1,\ldots,p_n\in X.$ In particular, we will consider:
\begin{itemize}
    \item Surfaces with $n$ labeled points,
    \item Handlebodies with $n$ labeled points on the boundary, and
    \item Topological realizations of graphs with $n$ labeled leaves. (Here $p_i$ is the 0-cell we have added to the end of the $i$-th leaf.)
\end{itemize}
An \emph{$n$-pointed map} $(X,p_1,\ldots,p_n)\to(X',p_1',\ldots,p_n')$ is a continuous map $X\to X'$ that sends $p_i\mapsto p_i'$ for $i=1,\ldots,n.$ An \emph{$n$-pointed homotopy (resp. isotopy)} between $n$-pointed maps is a homotopy (resp. isotopy) in which every intermediate map is also $n$-pointed.

\subsection{Terminology from algebraic and complex geometry}
Many of the objects in this paper can be viewed either through the lens of complex geometry or algebraic geometry (over the complex numbers). Here are some notes about our usage of terminology from both fields.
\begin{enumerate}
    \item All algebraic varieties in this paper are over the complex numbers.
    \item Every smooth algebraic variety over $\C$ has the natural structure of a complex manifold, though not every complex manifold is an algebraic variety. The algebraic analogue of an orbifold is the notion of \emph{Deligne-Mumford stack}, see \cite{Kresch2009} for an exposition. In particular, every smooth Deligne-Mumford stack over $\C$ has the natural structure of a complex orbifold, though not every complex orbifold is a Deligne-Mumford stack.
    \item A smooth projective irreducible algebraic curve (over $\C$) is the same as a compact Riemann surface without boundary. The space $\M_{g,n}$ can be viewed either as a moduli space of algebraic curves or as a moduli space of Riemann surfaces. We use these interchangeably in the rest of the paper. 
    \item See \cite{AdemLeidaRuan2007} for an exposition of orbifolds and their maps. For example, orbifolds have a notion of covering map. (The induced map of underlying topological spaces need not be a covering map.) Complex orbifolds have a notion of holomorphic map. Complex orbifolds also have a notion of ``\'etale map'', which is the orbifold analogue of ``local biholomorphism''.
    \item A \emph{normal crossings hypersurface} (see \cite{ChanGalatiusPayne2021}) in a complex manifold (say of complex dimension $n$) is a complex-codimension-1 subspace that locally looks like a union of coordinate hyperplanes in $\C^n$. (A \emph{normal crossings stratification} is a stratification by complex submanifolds that locally resembles the stratification of $\C^n$ by a union of coordinate hyperplanes.) A \emph{simple normal crossings hypersurface} is a normal crossings hypersurface whose irreducible components are smooth, i.e. have no self-crossings. There is also a notion of (simple) normal crossings hypersurface/stratification in a complex \emph{orbifold}. 
\end{enumerate}

\section{Background}

\subsection{Surfaces, multicurves, and dual graphs}\label{sec:SurfacesBackground} 

An $n$-pointed  surface of genus $g$ is a tuple $\Sgn=(\Sg,x_1,\ldots,x_n),$ where $\Sg$ is a closed oriented genus-$g$ surface, and $x_1,\ldots,x_n\in \Sg$ are distinct points.  

A \newword{simple closed curve} in $\Sgn$ is a simple closed curve in $\Sg\setminus\{x_1,\ldots,x_n\}.$ Such a curve $\gamma$ is called \newword{essential} if any disk in $\Sg$ whose boundary is $\gamma$ contains at least two of the labeled points $x_i$. A \newword{multicurve} $\Gamma$ in $\Sgn$ is a set of (distinct) isotopy classes of essential simple closed curves, that are ``disjoint'' in the sense that there exists a set of pairwise-disjoint representatives for the isotopy classes.  

\begin{figure}
\begin{center}
\begin{tikzpicture} [scale=.5]
   \draw []  plot [smooth cycle, tension=.8] coordinates 
 { (-6,0)   (-4.5 ,2)  (-2,1.25)   (0,2)   (2,1.25)   (4.5,2)  (6,0) (4.5,-2) (2,-1.25)  (0,-2) (-2,-1.25) (-4.5,-2)
    };
\draw  (.70,.1) arc (25:155:.75);
\draw   (.75,.3) arc (-10:-170:.75);
\begin{scope}[xshift=3.75cm]
\draw (.70,.1) arc (25:155:.75);
\draw (.75,.3) arc (-10:-170:.75);
\end{scope}
\begin{scope}[xshift=-3.75cm]
\draw (.70,.1) arc (25:155:.75);
\draw  (.75,.3) arc (-10:-170:.75);
\end{scope}

\draw[thick, red, dashed] (-5.2,0) ellipse (.8cm and .25cm);  
  \draw[thick, red] (-4.4,0) arc(0:180:.8cm and .25cm);   
  \node[left] (s1) at (-6,0) {$\gamma_1$}; 
  \draw[thick, red, dashed] (0,1.25) ellipse (.25cm and .725cm);
    \draw[thick, red] (0,.55) arc(-90:90:.25cm and .725cm);  
    \node[above] (s2) at (0,2) {$\gamma_2$};
     \draw[thick, red, dashed] (0,-1.15) ellipse (.25cm and .8cm);
     \draw[thick, red] (0,-2) arc(-90:90:.25cm and .8cm); 
     \node[below] (s2) at (0,-2) {$\gamma_3$}; 

\fill (4.25,1.25) circle(.05);
\node  [right] (x1) at (4.25,1.25) {$  x_1$};
\fill (4.25,-1.25) circle(.05);
\node  [right] (x2) at (4.25,-1.25) {$ x_2$};

\begin{scope}[xshift=15cm]
\draw [blue] (-2,0) .. controls (-6,3) and (-6,-3) .. (-2,0);
\draw[blue] (2,0) to (4.25,1.25); \node [right] (x11) at (4.25,1.25) {$x_1$};
\draw[blue] (2,0) to (4.25,-1.25);\node [right] (x22) at (4.25,-1.25) {$x_2$};
\draw [blue] ( -2,0) .. controls ( -1,1.5) and ( 1,1.5) .. ( 2,0);
\draw [blue] ( -2,0) .. controls ( -1,-1.5) and ( 1,-1.5) .. ( 2,0);
\draw [blue,fill=white] (-2,0) circle(.25);
\draw[blue] (-2,0) node {\tiny$0$};
\draw [blue,fill=white] (2,0) circle(.25);
\draw[blue] (2,0) node {\tiny$1$};
\end{scope}
\end{tikzpicture} 
\end{center}
\caption{ A    multicurve  $\Gamma=\{\gamma_1,\gamma_2,\gamma_3\}$ in $S_{3,2}$   and its dual graph $\tau_\Gamma$.  }
\label{fig:dual}
\end{figure}
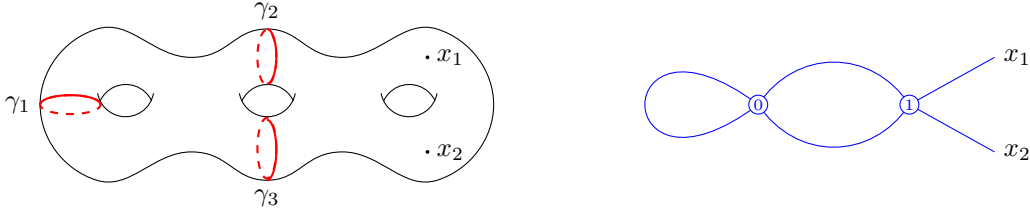

To each multicurve $\Gamma$ we will associate a {\em dual graph} $\tau_\Gamma$  whose vertices correspond to the components of $S_{g,n}\setminus \Gamma$, whose edges correspond to the curves in  $\Gamma$, and whose leaves correspond to the labeled points of $\Sgn$  (See Figure~\ref{fig:dual}). More precisely, let
$\{\gamma_i\}$ be a set of pairwise-disjoint representatives for  $\Gamma$, and let $N$ be a small tubular neighborhood of   $\{\gamma_i\}\cup\{x_j\}$; thus $N$ is homeomorphic to a disjoint union of annuli $\gamma_i\times[0,1]$ and   disks $D^2$ centered at the $x_j$.  The weighted graph $\tau_\Gamma$  is defined by 
 \begin{itemize}
 \item The half-edges $H$ are the boundary components of $N$, 
 \item Each piece $U$ of the partition consists of the elements of $H$ that are in the boundary of  a single connected component  of $S_{g,n}\setminus \Int(N)$, 
 \item  The weight on $U$ is the genus of  the corresponding component of $S\setminus \Int(N)$, and 
 \item The involution of $H$ switches the boundary components of each annulus and fixes the boundary component of each disk.  
 \end{itemize} 
 There is a continuous \newword{dual graph map} $\dgm_{S_{g,n},\Gamma}$
  from   $S_{g,n}$ to the topological realization $\overline\tau_\Gamma$ of $\tau_\Gamma$ given by sending each component of $S_{g,n}\setminus \Int(N)$ to the corresponding $0$-cell of $\tau_\Gamma$, each annulus $\gamma_i\times[0,1]$ to a $1$-cell joining the components containing the boundary of the annulus, and each disk centered at $p_i$ radially to a $1$-cell joining the the component containing the boundary of the disk to the $0$-cell labeled $p_i$.  

If $\Gamma,\Gamma'$ are multicurves and $\Gamma\subseteq\Gamma',$ then $\tau_\Gamma$ can be obtained from $\tau_{\Gamma'}$ by contracting edges (corresponding to the curves in $\Gamma'\setminus\Gamma$).

\subsection{Handlebodies, multidisks, and dual graphs} \label{sec:HandlebodiesBackground}

A \newword{genus-$g$ handlebody} $\Xg$ is an oriented topological 3-manifold obtained from a closed 3-ball by gluing $2g$ pairwise-disjoint closed disks in the boundary sphere together in pairs. The boundary $\partial \Xg$ of a genus-$g$ handlebody is an oriented genus-$g$ surface. Any two genus-$g$ handlebodies are homeomorphic. An \newword{$n$-pointed genus-$g$ handlebody} is a tuple $\Xgn=(\Xg,x_1,\ldots,x_n)$, where $\Xg$ is a genus-$g$ handlebody, and $x_1,\ldots,x_n\in \partial\Xg$ are distinct labeled points on the boundary of $\Xg$. For brevity we write $\partial\Xgn$ for the $n$-pointed surface $(\partial\Xg,x_1,\ldots,x_n)$.

A \newword{disk} in $\Xgn$ is an embedded copy of the closed disk whose boundary maps homeomorphically to a simple closed curve on $\partial\Xgn$. (In particular, a disk in $\Xgn$ does not pass through any labeled points.) Such a disk $\delta$ is called \newword{essential} if every disk contained in $\partial\Xg$ that is isotopic to $\delta$, if any such exists, contains at least two labeled points. An essential simple closed curve in $\partial\Xgn$ that bounds a disk in $\Xgn$ is called a \newword{meridian}.  A \newword{multidisk} in $\Xgn$ is a finite union of (distinct) isotopy classes of essential disks in $\Xgn.$, which have pairwise-disjoint representatives.   A multidisk is called \newword{pure} if every connected component of its complement is simply connected.  The boundary of a multidisk is a \newword{multimeridian}.

For a multidisk $\Delta$ in $ \Xgn$, the ($n$-leafed weighted) \newword{dual graph} $\tau_\Delta$ of $\Delta$ is defined in the same way as  the dual graph of a multicurve, i.e. its  vertices are the connected components of $\Xgn\setminus\Delta$  (weighted by genus), its  edges correspond to the  disks in $\Delta$, and its leaves correspond to the labeled points. (In fact it is \emph{equal} to the dual graph to the multicurve $\partial\Delta$.) The graph $\tau_\Delta$  can be formally defined using a pairwise disjoint set $\{\delta_i\}$ of essential disks representing $\Gamma$ as follows:  Let $N$ be a tubular neighborhood   of $\{\delta_i\}\cup\{x_j\}$. The elements of $H$ are then the components of   $\partial N\cap \Int(\Xgn)$ (which are all embedded open disks), the pieces of the partition are determined by  the components of $\Xgn\setminus N$ and the involution swaps the two elements of $H$ parallel to each single disc $\delta_i$ and fixes the element in the neighborhood of a labeled point.  As before, there is a continuous map $\dgm_{\Xgn,\Delta}$ from $\Xgn$ to the topological realization $\overline\tau_\Delta$ of $\tau_\Delta$.

\subsection{Mapping class groups}\label{sec:MappingClassGroups}
The \newword{ mapping class group $\MCG(\Sgn)$ of an $n$-pointed surface} is the group of orientation-preserving  homeomorphisms $\Sgn\to\Sgn$ fixing the labeled points, modulo the subgroup consisting of homeomorphisms that are $n$-pointed-isotopic to the identity\footnote{This is sometimes referred to as the \emph{pure} mapping class group, to distinguish it from the case where the homeomorphisms are allowed to permute the labeled points. We do not use this terminology as we reserve the word ``pure'' for the meaning more common in tropical geometry, see Section \ref{sec:ModuliOfGraphs} and Section \ref{sec:HandlebodiesBackground} just above.} (see Section \ref{sec:nPointedMaps}). By a classic result of Max Dehn, $\MCG(\Sgn)$ is generated by Dehn twists in simple closed curves.   

For an $n$-pointed handlebody $\Xgn,$ the \newword{mapping class group $\MCG(\Xgn)$} is the group of orientation-preserving homeomorphisms $\Xgn\to\Xgn$ fixing the labeled points, modulo the subgroup consisting of homeomorphisms that are $n$-pointed-isotopic to the identity.
\begin{prop}[{\cite[Lem. 3.1, Cor. 5.11]{Hensel2020}}]\label{prop:MCGInjection}
    The natural map $\MCG(\Xgn)\to\MCG(\partial\Xgn)$, which sends a homeomorphism $g:\Xgn\to \Xgn$ to the restriction $g|_{\partial\Xgn}$, is injective. Its image is the set of mapping classes that send meridians to meridians and non-meridians to non-meridians.
\end{prop}
 
As a result, $\MCG(\Xgn)$ is sometimes referred to as the \newword{handlebody subgroup} of $\MCG(\partial\Xgn)$. 

We define the  \newword{homotopy mapping class group} $\hMCG(\Xgn)$ to be the group of  homotopy equivalences $\Xgn\to\Xgn$ fixing the labeled points, modulo the subgroup consisting of homotopy equivalences that are $n$-pointed-homotopic to the identity.
 
Since any homeomorphism is a homotopy equivalence, there is a natural map $\MCG(\Xgn)\to\hMCG(\Xgn)$ --- this map is studied in detail in Section \ref{sec:LuftMarkedPoints}.

\begin{rem}
    The homotopy mapping class group of $\Xgn$ is an invariant of the $n$-pointed-homotopy equivalence class of $\Xgn$. Since $\Xgn$ is $n$-pointed-homotopy equivalent to the $n$-thorned rose (see Section \ref{sec:GraphConventions}), $\hMCG(\Xgn)$ is isomorphic to the group $A_{g,n}$ in Section \ref{sec:IntroOuterSpace}, see also Section \ref{sec:OuterSpaceBackground}.
\end{rem}
\begin{ex}[Running example, $(g,n)=(1,1)$]\label{ex:V110}
    For the handlebody $V_{1,1},$ we classically have $\MCG(\partial V_{1,1})\cong\SL_2(\Z),$ where the identification comes from choosing a basis of $H_1(\partial V_{1,1},\Z)\cong\Z^2$, see \cite[Sec. 2.2.4]{FarbMargalit2011}. (To match our conventions, we let $\SL_2(\Z)$ act on the right.)
    
    There are countably many isotopy classes of simple closed curves in $\partial V_{1,1}$; these correspond to the primitive elements in $\Z^2$. In $V_{1,1},$ there is a unique essential disk $\delta$ up to isotopy. The boundary $\partial\delta$ is the unique meridian in $\partial V_{1,1}$ up to isotopy. The two possible orientations on $\partial\delta$ determine two distinguished primitive elements of $H_1(\partial V_{1,1},\Z)$, which differ by a sign. Without loss of generality we may identify these elements with $$\begin{pmatrix}
        \pm1\\0
    \end{pmatrix}\in\Z^2.$$   By Proposition \ref{prop:MCGInjection}, the subgroup $\MCG(V_{1,1})\subset\MCG(\partial V_{1,1})$ is the setwise stabilizer of these two elements, i.e. under our identification $$\MCG(\X_{1,1})\cong\left\{\pm\begin{pmatrix}
        1&0\\
        a&1
    \end{pmatrix}:a\in\Z\right\}\subset\SL_2(\Z).$$
\end{ex}

\bigskip

The following proposition shows how forgetting a labeled point affects $\MCG(\Xgn)$ and $\hMCG(\Xgn)$ --- these are the analogues of the Birman exact sequence on mapping class groups, see \cite[Sec. 4.2]{FarbMargalit2011}. 
\begin{prop}[{\cite[Lem. 3.3]{Hensel2020}}]\label{prop:hMCGAddPoint}
    If $n=0,$ then $\hMCG(\Xgn)$ is naturally isomorphic to $\Out(\pi_1(\Xg))$  (which is non-canonically isomorphic to $\Out(F_g)$, where $F_g$ denotes the free group). If $n>0,$ then we have natural exact sequences that fit into a commutative diagram
    \begin{equation}
        \begin{tikzcd}
        1\arrow[r]&\pi_1(\partial\Xg\setminus\{x_1,\ldots,x_{n-1}\},x_n)\arrow[r]\arrow[d]&\MCG(\Xgn)\arrow[r]\arrow[d]&\MCG(\X_{g,n-1})\arrow[r]\arrow[d]&1\\
            1\arrow[r]&\pi_1(\Xg,x_n)\arrow[r]&\hMCG(\Xgn)\arrow[r]&\hMCG(\X_{g,n-1})\arrow[r]&1
        \end{tikzcd}
    \end{equation}
    Here the left vertical arrow is the pushforward along the inclusion map, and for a loop in $\pi_1(\partial\X\setminus\{x_1,\ldots,x_{n-1}\},x_n)$ (resp. $\pi(\Xg,x_n)$), its image in $\MCG(\Xgn)$ (resp. $\hMCG(\Xgn)$) is the ``point-pushing map'' around $\gamma$, see \cite[Lem. 3.3]{Hensel2020}. In particular, $\hMCG(\Xgn)\cong A_{g,n}.$
\end{prop}
If $\delta$ is an essential disk in $\Xgn$, the Dehn twist in its boundary curve extends to a homeomorphism of $\Xgn$, which is  called a Dehn twist  around $\delta$.  We let $\Tw(\Xgn)\subseteq\MCG(\Xgn)$ denote the \newword{twist subgroup} of $\MCG(\Xgn)$, i.e. the subgroup generated by all Dehn twists around essential disks in $\Xgn$.

\begin{prop}[{\cite[Prop. 1.7]{Otal1989}}]\label{prop:TwistSubgroupTorsionFree}
    $\Tw(\Xgn)$ is torsion-free.
\end{prop}

\subsection{Complex structures and moduli of curves}\label{sec:ModuliTeichmuller}
A \newword{complex structure} on an $n$-pointed genus-$g$ surface $\Sgn$ is a pair $(C,\phi)$, where $C=(C,p_1,\ldots,p_n)$ is an $n$-pointed Riemann surface and $\phi:\Sgn\to C$ is an $n$-pointed homeomorphism. The \newword{Teichm\"uller space} $$\T(\Sgn)=\{(C,\phi)\}/\sim$$ of $\Sgn$ is the set of complex structures on $\Sgn$, where two complex structures $(C_1,\phi_1)$ and $(C_2,\phi_2)$ are considered equivalent if there exists a commutative diagram 
\begin{align*}
    \begin{tikzcd}[ampersand replacement=\&]
        \Sg\arrow[r,"\phi_1"]\arrow[d,"\alpha"]\&C_1\arrow[d,"\beta"]\\
        \Sg\arrow[r,"\phi_2"]\&C_2
    \end{tikzcd}
\end{align*}
where $\beta:C_1\to C_2$ is an $n$-pointed isomorphism of Riemann surfaces and $\alpha:\Sgn\to \Sgn$ is a homeomorphism that is $n$-pointed-isotopic to the identity. $\T(\Sgn)$ has the natural structure of a complex manifold of dimension $3g-3+n$.

$\MCG(\Sgn)$ acts on the right on $\T(\Sgn)$ by biholomorphisms, via $(C,\phi)\cdot h=(C,\phi\circ h)$. The action is properly discontinuous, and the quotient is the \newword{moduli space $\M_{g,n}$ of $n$-pointed genus-$g$ Riemann surfaces} (actually a Deligne-Mumford stack/algebraic complex orbifold).

\subsection{Stable curves, their moduli, and augmented Teichm\"uller space}\label{sec:StableCurvesTeichmuller}
A \newword{stable $n$-pointed genus-$g$ curve} over $\C$ is a connected projective algebraic curve $(C,x_1,\ldots,x_n)$ over $\C$ of arithmetic genus $g$ with at-worst-nodal singularities (see e.g. \cite{KockVainsencher}), with $n$ distinct labeled smooth points, such that (1) every genus-zero irreducible component of $S$ contains at least three ``special points'' (labeled points or nodes), and (2) every genus-one irreducible component of $S$ contains at least one special point. Just like a multicurve on an $n$-pointed genus-$g$ surface, a stable $n$-pointed genus-$g$ curve $(C,x_1,\ldots,x_n)$ has an ($n$-leafed weighted) \newword{dual graph} $\tau_{(C,x_1,\ldots,x_n)}$, whose vertices are the  irreducible components of $C$ (weighted by geometric genus), whose  edges are the nodes $\{y_j\}$ of $C$, and whose labeled leaves are the points $x_1,\ldots,x_n$. Formally, the graph $\tau_{(C,x_1,\ldots,x_n)}$ is defined by taking a small neighborhood of $\{x_i\}\cup\{y_j\}$ and letting the half-edges be the boundary components of this neighborhood.  Thus the dual graph of a stable $n$-pointed genus-$g$ curve is a stable $n$-leafed genus-$g$ weighted graph. As before,  there is a   continuous $n$-pointed ``dual graph map'' $\dgm_C$ from 
$C$ to the topological realization of $\tau_{(C,x_1,\ldots,x_n)}$ which is canonical up to homotopy.

\begin{figure}
\begin{center}
\begin{tikzpicture} [scale=.5]
\draw (-3,0) to[out=115,in=0] (-4,.5) to[out=180,in=0] (-5,0) to[out=180,in=-90] (-5.75,1) to[out=90,in=180] (-4,2) to[out=0,in=180] (-2,1.25) to[out=0,in=180] (-.75,1.5) to[out=0,in=90] (0,1) to[out=-90,in=65] (-.75,0);
\draw (-2.9,.2) to[out=-115,in=0] (-4,-.5) to[out=180,in=0] (-5,0) to[out=180,in=90] (-5.75,-1) to[out=-90,in=180] (-4,-2) to[out=0,in=180] (-2,-1.25) to[out=0,in=180] (-.75,-1.5) to[out=0,in=-90] (0,-1) to[out=90,in=-65] (-.85,0.2);
\draw (6,0) to[out=90,in=0] (4,2) to[out=180,in=0] (2,1.25) to[out=180,in=0] (.75,1.5) to[out=180,in=90] (0,1) to[out=-90,in=115] (.75,0);
\draw (6,0) to[out=-90,in=0] (4,-2) to[out=180,in=0] (2,-1.25) to[out=180,in=0] (.75,-1.5) to[out=180,in=-90] (0,-1) to[out=90,in=-115] (.85,0.2);


\begin{scope}[xshift=3.75cm]
\draw (.70,.1) arc (25:155:.75);
\draw (.75,.3) arc (-10:-170:.75);
\end{scope}

\fill (4.25,1.25) circle(.05);
\node  [right] (x1) at (4.25,1.25) {$  x_1$};
\fill (4.25,-1.25) circle(.05);
\node  [right] (x2) at (4.25,-1.25) {$ x_2$};

\begin{scope}[xshift=15cm]

\draw [blue] (-2,0) .. controls (-6,3) and (-6,-3) .. (-2,0);
\draw[blue] (2,0) to (4.25,1.25); \node [right] (x11) at (4.25,1.25) {$x_1$};
\draw[blue] (2,0) to (4.25,-1.25);\node [right] (x22) at (4.25,-1.25) {$x_2$};
\draw [blue] ( -2,0) .. controls ( -1,1.5) and ( 1,1.5) .. ( 2,0);
\draw [blue] ( -2,0) .. controls ( -1,-1.5) and ( 1,-1.5) .. ( 2,0);
\draw [blue,fill=white] (-2,0) circle(.25);
\draw[blue] (-2,0) node {\tiny$0$};
\draw [blue,fill=white] (2,0) circle(.25);
\draw[blue] (2,0) node {\tiny$1$};
\end{scope}
\end{tikzpicture} 
\end{center}
\caption{ A stable curve $(C,x_1,x_2)$ in $\Mbar_{3,2}$   and its dual graph $\tau_{(C,x_1,x_2)}$.  }
\label{fig:dualOfCurve}
\end{figure}
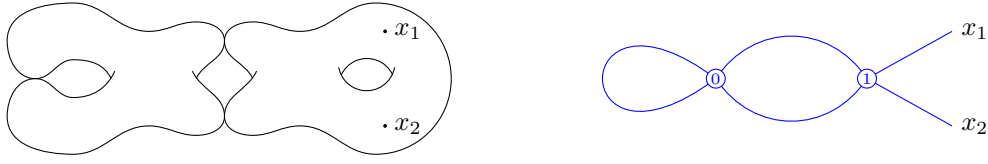

\subsubsection{The moduli space of stable curves}\label{sec:StableCurves} Deligne and Mumford showed that the \newword{moduli space $\Mbar_{g,n}$ of stable $n$-pointed genus-$g$ curves} is a smooth irreducible Deligne-Mumford stack. In particular it is a complex orbifold that contains $\M_{g,n}$ as a dense open suborbifold. The boundary $\Mbar_{g,n}\setminus\M_{g,n}$ has a stratification by closed suborbifolds indexed by stable $n$-leafed genus-$g$ weighted graphs, i.e. for any such graph $\tau,$ we have a ``locally closed \newword{boundary stratum}'' 
$$\M_\tau=\{(C,x_1,\ldots,x_n)\in\Mbar_{g,n}:(C,x_1,\ldots,x_n)\text{ has dual graph }\tau\}.$$ 
The stratum $\M_\tau$ has a recursive structure, given by a natural isomorphism 
$$\M_\tau\cong\left(\prod_{\substack{\text{$v$ internal vertex of}\\ \tau_{(C,x_1,\ldots,x_n)}}}\M_{g_v,n_v}\right)\bigg/\Aut(\tau_{(C,x_1,\ldots,x_n)}),$$ 
where $g_v$ is the geometric genus of the irreducible component $v$, and $n_v$ is the number of special points on $v$. The closure $\Mbar_\tau$ of $\M_\tau$ in $\Mbar_{g,n}$ is characterized by 
$$\Mbar_\tau=\bigsqcup_{\tau'}\M_{\tau'},$$
 where $\tau'$ runs over all stable $n$-leafed genus-$g$ weighted graphs such that $\tau$ can be obtained from $\tau'$ by internal edge contractions. We refer to $\Mbar_\tau$ as a ``closed boundary stratum''. The isomorphism above extends to give $$\Mbar_\tau\cong\left(\prod_{\substack{\text{$v$ a vertex of}\\ \tau_{(C,x_1,\ldots,x_n)}}}\Mbar_{g_v,n_v}\right)\bigg/\Aut(\tau_{(C,x_1,\ldots,x_n)}),$$ so in particular each closed boundary stratum is a closed suborbifold of $\Mbar_{g,n}.$

\subsubsection{Augmented Teichm\"uller space} \label{sec:AugmentedTeichmullerSpace} Bers \cite{Bers1974} introduced a partial compactification of $\T(\Sgn)$, as follows. The \newword{augmented Teichm\"uller space} $$\Tbar(\Sgn)=\{(C,\phi)\}/\sim$$ of $\Sgn$ is the set of pairs $(C,\phi)$ where $C=(C,p_1,\ldots,p_n)$ is an $n$-pointed genus-$g$ stable curve, and $\phi:\Sgn\to C$ is a surjective $n$-pointed map satisfying:
\begin{itemize}
    \item $\phi$ is a homeomorphism on the complement of some multicurve $\Gamma$ in $\Sgn$, and
    \item $\phi$ contracts each component of $\Gamma$ to a node of $C$,
\end{itemize}
modulo the same equivalence relation as for $\T(\Sgn)$. $\Tbar(S_{g,n})$ is contractible, see \cite{Wolpert2009} --- this follows from the fact that $\Tbar(S_{g,n})$ has a $CAT(0)$ metric. The action of $\MCG(\Sgn)$ on $\T(\Sgn)$ extends to $\Tbar(\Sgn)$. However, unlike $\T(\Sgn)$, $\Tbar(\Sgn)$ is not a complex manifold, as shown in the following example:
\begin{ex}[Running example, continued]\label{ex:V111}
 For a 1-pointed genus-1 surface $S_{1,1}$, the augmented Teichm\"uller space $\Tbar(S_{1,1})$ is homeomorphic to $\mathbb{H}\cup\Q\cup\{\infty\}$, where $\mathbb{H}$ is the (open) upper half-plane in $\C$. (Under this identification, $\T(S_{1,1})$ is identified with $\mathbb{H}$, which is indeed a complex manifold.) As in Example \ref{ex:V110}, we have $\MCG(S_{1,1})\cong\SL_2(\Z)$, which acts on $\mathbb{H}\cup\Q\cup\{\infty\}$ via the usual (right) action: $$z\cdot\begin{pmatrix}
        a&b\\
        c&d
    \end{pmatrix}=\frac{az+c}{bz+d}.$$
    Note that the action is not quite free; it has a kernel of order two. Note also that the points of $\Q\cup\{\infty\}$ have infinite stabilizers. For example, the stabilizer of $\infty$ consists of all elements of $\SL_2(\Z)$ with $b=0.$
\end{ex}
\medskip

Although $\Tbar(S_{g,n})$ is not a complex manifold (as we see in the above example), it has the structure of a countable disjoint union of complex manifolds, with a recursive stratification analogous to that of $\Mbar_{g,n}$ just above, as follows. Any element $(C,\phi)\in\Tbar(\Sgn)$ has an associated (possibly empty, contracted) multicurve $\Gamma$. Consider the locally closed \newword{boundary stratum} $$\T_\Gamma(\Sgn)=\{(C,\phi)\in\Tbar(\Sgn):(C,\phi)\text{ has associated multicurve }\Gamma\}.$$ We have an isomorphism
\begin{align}\label{eq:TSBoundaryRecursion}
    \T_\Gamma(\Sgn)\cong\prod_{\text{$v$ a vertex of $\tau_\Gamma$}}\T(\Sgn\langle v\rangle),
\end{align} where for each connected component $v$ of $\Sgn\setminus\Gamma$, $\Sgn\langle v\rangle$ is the pointed surface contracting each boundary component of (the closure in $\Sgn$ of) $v$ to a point, and labeling that point. In particular, $\T_\Gamma(\Sgn)$ has the natural structure of a complex manifold, of dimension $3g-3+n-\abs{\Gamma}$, where $\abs{\Gamma}$ is the number of simple closed curves in $\Gamma$.

A stratum $\T_{\Gamma'}(\Sgn)$ is contained in the closure of $\T_\Gamma(\Sgn)$ if and only if $\Gamma\subseteq\Gamma'.$ In particular, the closure $\Tbar_\Gamma(\Sgn)$ of $\T_\Gamma(\Sgn)$ in $\Tbar(\Sgn)$ can be characterized by:
$$\Tbar_\Gamma(\Sgn)=\bigsqcup_{\Gamma'\supseteq\Gamma}\T_{\Gamma'}(\Sgn).$$

We refer to $\Tbar_\Gamma(\Sgn)$ as a ``closed boundary stratum''. The recursive structure of $\T_\Gamma(\Sgn)$ extends: \begin{align}\label{eq:TbarSBoundaryRecursion}
    \Tbar_\Gamma(\Sgn)\cong\prod_{\text{$v$ a vertex of $\tau_\Gamma$}}\Tbar(\Sgn\langle v\rangle).
\end{align}

\begin{rem}\label{rem:HInverseMulticurve}
    The action of $\MCG(\Sgn)$ on $\T(\Sgn)$ extends naturally to $\Tbar(\Sgn)$. Note that if $(C,\phi)\in\Tbar(\Sgn)$ has contracted multicurve $\Gamma,$ then the multicurve contracted by $(C,\phi)\cdot h$ is $h^{-1}(\Gamma)$. 
\end{rem}
\begin{rem}
    Unlike the action of $\MCG(\Sgn)$ on $\T(\Sgn),$ the action of $\MCG(\Sgn)$ on $\Tbar(\Sgn)$ is not properly discontinuous. Indeed,  if $(C,\phi)\in\Tbar(\Sgn)\setminus \T(\Sgn)$ has contracted multicurve $\Gamma,$ then powers of the Dehn twist around any simple closed curve $\gamma\in\Gamma$ stabilize $(C,\phi),$ so the stabilizer of $(C,\phi)$ is infinite.
\end{rem}

We will need the following basic properties of this action.
\begin{lem}[{\cite[Prop. 2.5]{HubbardKoch2014}}]
    Let $\Gamma$ be a multicurve in $\Sgn$, and let $\Tw(\Gamma)\subseteq\MCG(\Sgn)$ be the subgroup generated by Dehn twists around simple closed curves in $\Gamma$, let $\MCG(\Sgn,\Gamma)\subseteq\MCG(\Sgn)$ denote the subgroup of mapping classes whose induced action on $\tau_\Gamma$ is trivial. (That is, $\MCG(\Sgn,\Gamma)$ consists of elements that preserve each curve in $\Gamma$ and that do not permute the components of $\Sg\setminus\Gamma$.) Then we have an exact sequence $$1\to \Tw(\Gamma)\to\MCG(\Sgn,\Gamma)\to\prod_{\text{$v$ vertex of $\tau_\Gamma$}}\MCG(\Sgn\langle v\rangle)\to1.$$
\end{lem}
\begin{lem}[Adapted from {\cite[Prop. 2.6]{HubbardKoch2014}}]\label{lem:StabilizerOfPointInTbarExactSequence}
    Let $\Gamma$ be a multicurve in $\Sgn$, and let $(C,\phi)\in\T_\Gamma(\Sgn).$ Let $\Stab(C,\phi)\subseteq\MCG(\Sgn)$ denote the stabilizer. Then we have an exact sequence $$1\to \Tw(\Gamma)\to\Stab(C,\phi)\to\Aut(C)\to1.$$
\end{lem}

As noted above, the action of $\MCG(\Sgn)$ on $\Tbar(\Sgn)$ is not properly discontinuous. However, the action is otherwise well-behaved in several ways:
\begin{itemize}
    \item The orbits are discrete.
    \item The action is compatible with the boundary stratification, in the sense that an element of $\MCG(\Sgn)$ sends each boundary stratum in $\Tbar(\Sgn)$ biholomorphically to another boundary stratum.
    \item Furthermore, the stabilizer of a stratum $\T_\Gamma(\Sgn)$ acts on $\T_\Gamma(\Sgn)$ with infinite kernel $\Tw(\Gamma)$, but the induced action of the quotient on $\Tbar_\Gamma(\Sgn)$ is properly discontinuous.
\end{itemize}
The above properties imply that the quotient $\Tbar(\Sgn)/\MCG(\Sgn)$ inherits a boundary stratification from $\Tbar(\Sgn)$, whose strata are locally closed connected subsets, indexed by the set of orbits of the $\MCG(\Sgn)$-action on the set of boundary strata of $\Tbar(\Sgn)$. Since each stratum of $\Tbar(\Sgn)$ has the structure of a complex manifold, the above properties imply that each boundary stratum in $\Tbar(\Sgn)/\MCG(\Sgn)$ inherits the structure of a complex orbifold.

This does not prove that $\Tbar(\Sgn)/\MCG(\Sgn)$ is itself a complex orbifold. However, this turns out to be the case --- it is a nontrivial theorem:
\begin{thm}[\cite{Harvey1974,Braungardt2001,HinichVaintrob2010,HubbardKoch2014}]\label{thm:MgnComplexOrbifold}
    $\Tbar(\Sgn)/\MCG(\Sgn)$ is naturally endowed with the structure of a complex orbifold, and this orbifold is canonically isomorphic to $\Mbar_{g,n}$. With respect to this structure,
    \begin{itemize}
        \item The restriction of the quotient map to each complex manifold $T_\Gamma(\Sgn)\subseteq\Tbar(\Sgn)$ is a holomorphic covering map over its image.
        \item The preimage of a locally closed (resp. closed) boundary stratum $\M_\tau\subseteq\Mbar_{g,n}$ (resp. $\Mbar_\tau$) under the quotient map is the union $\bigsqcup_{\Gamma}\T_\Gamma(\Sgn)$ (resp. $\bigsqcup_{\Gamma}\Tbar_\Gamma(\Sgn)$) where $\Gamma$ runs over all multicurves in $\Sgn$ with dual graph isomorphic to $\tau.$
    \end{itemize}
\end{thm}

\begin{ex}[Running example, continued]\label{ex:V111.5}
    Using the identifications in Examples \ref{ex:V110} and \ref{ex:V111}, the quotient $(\mathbb{H}\cup\Q\cup\{\infty\})/\SL_2(\Z)$ can be identified with $\Mbar_{1,1},$ by Theorem \ref{thm:MgnComplexOrbifold}. This special case is precisely the classical construction of the compactified modular curve, see e.g. \cite[Sec. 1.5]{Shimura1971}.
\end{ex}

\begin{rem}\label{rem:InheritedStratification}
    We generalize Theorem \ref{thm:MgnComplexOrbifold} in Section \ref{sec:ComplexStructuresTechnical}. Specifically, for any subgroup $G\subseteq\MCG(\Sgn),$ we show that there is a nonempty maximal $G$-invariant open subset $U(G)\subseteq\Tbar(\Sgn)$ such that $U(G)/G$ admits the natural structure of a complex orbifold. The open set $U(G)$ is a union of locally closed boundary strata.
\end{rem}

\subsection{Simplicial complexes}\hfill\label{sec:Complexes}

\subsubsection{The curve complex} Associated to an $n$-pointed genus-$g$ surface $\Sgn$ is an infinite simplicial complex $\mathcal{C}(\Sgn)$ of dimension $3g-3+n$, called the \newword{curve complex}, whose vertices are isotopy classes of essential simple closed curves in $\Sgn$, and whose simplices are multicurves. In barycentric coordinates, a \emph{point} of $\mathcal{C}(\Sgn)$ can be interpreted as a multicurve $\Gamma$ in $\Sgn$, with positive real numbers that sum to 1 labeling the simple closed curves in $\Gamma$. Note that $\MCG(\Sgn)$ acts naturally on the right on $\mathcal{C}(\Sgn)$ by $\Gamma\cdot h=h^{-1}(\Gamma)$, see Remark \ref{rem:HInverseMulticurve}. This is an action by simplicial automorphisms.

\subsubsection{The disk complex}\label{sec:DiskComplex} Associated to an $n$-pointed genus $g$ handlebody $\Xgn$ is an infinite simplicial complex $\mathcal{D}(\Xgn)$ of dimension $3g-3+n$, called the \newword{disk complex}, whose vertices are $n$-pointed isotopy classes of essential disks, and whose simplices are multidisks. In barycentric coordinates, a \emph{point} of $\mathcal{D}(\Xgn)$ can be interpreted as a multidisk $\Delta$ in $\Xgn$, with positive real numbers that sum to 1 assigned to the disks in $\Delta$. The \newword{pure locus} is $\mathcal{D}^{\pure}(\Xgn)\subseteq\mathcal{D}(\Xgn)$ is the dense open set consisting of all simplices corresponding to pure multidisks $\Delta$ --- that is, multidisks $\Delta$ such that every component of $\Xg\setminus\Delta$ is simply connected. $\D^{\pure}(\Xgn)$ is not a subcomplex of $\D(\Xgn)$ --- rather it is the complement of a subcomplex.
\begin{prop}[\cite{McCullough1991,Giansiracusa2011}]\label{prop:SimpleDisksContractible}
    If $g>0$, then $\mathcal{D}(\Xgn)$ and $\mathcal{D}^{\pure}(\Xgn)$ are contractible.
\end{prop}

\begin{prop}\label{prop:DiskComplexInjectsIntoCurveComplex}
    The natural map $\mathcal{D}(\Xgn)\to\mathcal{C}(\partial\Xgn)$, which sends a multidisk to its boundary multicurve, is an embedding of simplicial complexes.
\end{prop}
\begin{proof}
    Since the map is clearly compatible with face inclusions, the statement is immediate from the fact that a simple closed curve in $\partial\Xgn$ bounds at most one disk in $\Xgn$ up to homotopy \cite[Lem. 2.3]{Hensel2020}.
\end{proof}
\begin{rem}
    The image of the above embedding is a full subcomplex of $\mathcal{C}(\partial\Xgn)$ --- this follows from the fact that two disjoint meridians can be made to bound disjoint disks, which follows from an innermost-circles argument.
\end{rem}

\begin{rem}
    $\MCG(\Xgn)$ acts naturally on the right on $\mathcal{D}(\Xgn)$ by $\Delta\cdot h=h^{-1}(\Delta)$, see Remark \ref{rem:HInverseMulticurve}. This is an action by simplicial automorphisms, and the embedding from Proposition \ref{prop:DiskComplexInjectsIntoCurveComplex} is equivariant with respect to the inclusion $\MCG(\Xgn)\into\MCG(\partial\Xgn)$ from Proposition \ref{prop:MCGInjection}. (In other words, if we view $\mathcal{D}(\Xgn)$ as a subset of $\mathcal{C}(\partial\Xgn)$ via Proposition \ref{prop:DiskComplexInjectsIntoCurveComplex}, then this subset is invariant under the $\MCG(\Xgn)$-action on $\mathcal{C}(\partial\Xgn)$.)
\end{rem}

\subsubsection{Outer Space}\label{sec:OuterSpaceBackground}
Outer space $\CV_g$ for $g\geq 2$ is the moduli space of marked stable metric graphs of genus $g$ without leaves whose edge lengths sum to one \cite{CullerVogtmann1986}. These graphs are not weighted, so  the word ``stable" amounts to saying that every vertex has valence at least three. Culler and Vogtmann proved that $\CV_g$ is contractible, and the group $\Out(F_g),$ realized as the group $\pi_0(HE(R_g))$ of homotopy classes of self-homotopy equivalences of the $g$-petaled rose $R_g$, acts with finite stabilizers.  This picture was extended by Hatcher \cite{Hatcher1995}  to include graphs with leaves for all $g,n$ with $2g+n\geq 3$.  Specifically, we define  $\CV_{g,n}$ to be the space of marked stable metric graphs of genus $g$ with $n$  labeled leaves, where a marking is an $n$-pointed homotopy equivalence from the $n$-thorned rose  $R_{g,n}$. Hatcher proved the space $\CV_{g,n}$ is contractible when $g>0$, and the group $A_{g,n}:=\pi_0(HE(R_{g,n}))$ acts on $\CV_{g,n}$ on the right by the formula $(\tau,r)\cdot\alpha=(\tau,r\circ\alpha)$ \cite{CullerVogtmann1986,Hatcher1995}. In particular, the stabilizer of a point is the group of automorphisms of the graph that preserve edge-lengths and leaves.

$\CV_{g,n}$ is a union of open simplices, and for $g\geq 2$ there are infinitely many of them. The barycentric coordinates on each simplex correspond to all possible assignments of positive real lengths that sum to 1 on the edges of a given marked graph $(\tau,r)$. For $g>0$ it is not a simplicial complex, since some faces are missing, but Culler-Vogtmann and Hatcher \cite{CullerVogtmann1986,Hatcher1995} observed that there is a simplicial complex $\CV_{g,n}^*$ of dimension $3g-3+n$ into which $\CV_{g,n}$ embeds as the complement of a natural subcomplex. Furthermore, Hatcher's argument that $\CV_{g,n}$ is contractible when $g>0$ extends to show that $\CV_{g,n}^*$ is contractible when $g>0$ \cite{Hatcher1995}. The simplicial complex $\CV_{g,n}^*$ is called the \newword{simplicial completion} of $\CV_{g,n}$. The action of $A_{g,n}$ on $\CV_{g,n}$ extends to a simplicial  action on $\CV^*_{g,n},$ but is no longer properly discontinuous since the stabilizer of any point in $\CV_{g,n}^*\setminus \CV_{g,n}$ is infinite.

For the purposes of this paper, we define   $\CV_{g,n}^*$ to be the space of equivalence classes of marked stable metric graphs of genus $g$ with $n$ leaves, with edge lengths allowed to be zero. (Again, these graphs are not weighted.) The equivalence relation is generated by the following elementary equivalences:
\begin{enumerate}
    \item $(\tau_1,r_1)$ is equivalent to $(\tau_2,r_2)$ if there is an $n$-pointed, edge-length-preserving isomorphism $\Phi:\tau_1\to\tau_2$ such that $\Phi\circ r_1$ is $n$-pointed homotopic to $r_2$, and
    \item $(\tau_1,r_1)$ is equivalent to $(\tau_2,r_2)$ if $(\tau_1,r_1)$ is obtained from $(\tau_2,r_2)$ by contracting a length-zero edge in $\tau_2$ which is not a loop. (The marking $r_1$ is the composition of $r_2:R_{g,n}\to\tau_2$ and the edge contraction map $\tau_2\to\tau_1$.)
\end{enumerate}
\begin{rem}\label{rem:SimplicesOfCV}
    The \emph{set} of equivalence classes of marked graphs above is the underlying topological \emph{space} of the simplicial complex $\CV_{g,n}^*$. The simplices of $\CV_{g,n}^*$ are indexed by equivalence classes $(\tau,r)$ of marked stable genus-$g$ (unweighted) $n$-leafed graphs, \emph{without} a metric, but with the specification of a collection of ``length-zero edges'', where the equivalence relation is the same as above. For such an equivalence class $(\tau,r)$, the dimension of the corresponding simplex is $e-1$, where $e$ denote the number of edges \emph{not} specified to have length zero.
\end{rem}
\begin{rem} In Appendix \ref{app:SimplicialCompletion} we relate the above characterization of $\CV_{g,n}^*$ to    Hatcher's original definition in terms of spheres in a doubled handlebody (See Proposition \ref{prop:PointsOfSimplicialCompletion}).
     
\end{rem}

\subsubsection{The moduli space of tropical curves}\label{sec:ModuliOfGraphs}

A \newword{tropical curve} is a stable weighted metric graph, as defined in Section \ref{sec:GraphConventions}. The \newword{moduli space of tropical curves} $\M_{g,n}^{\trop}$ is the space of   tropical curves of genus $g$ with $n$ leaves, up to isomorphism preserving the metric and each of the labeled leaves. The \newword{link} $\Link\M_{g,n}^{\trop}\subseteq\M_{g,n}^{\trop}$ is the locus of metric graphs whose edge lengths sum to 1.

There are canonical isomorphisms, as noted in \cite{ChanGalatiusPayne2021}:
\begin{itemize}
    \item $\mathcal{C}(S_{g,n})/\MCG(S_{g,n})\xrightarrow{\cong}{}\Link\M_{g,n}^{\trop}$, sending a multicurve $\Gamma$ to its dual graph $\tau_\Gamma$. The barycentric coordinates on $\mathcal{C}(S_{g,n})$ determine an assignment of a positive real length to each edge of $\tau_\Gamma$, so that the edge lengths sum to 1.
    \item $\D(\Xgn)/\MCG(\Xgn)\xrightarrow{\cong}{}\Link\M_{g,n}^{\trop}$, sending a multidisk $\Delta$ to its dual graph $\tau_\Delta$.
    \item $\CV_{g,n}^*/A_{g,n}\xrightarrow{\cong}{}\Link\M_{g,n}^{\trop}$, forgetting the marking, contracting all length-zero edges, and adding appropriate vertex weights.
\end{itemize}

Since the action of $A_{g,n}$ on $\CV_{g,n}^*$ is simplicial, $\M_{g,n}^{\trop}$ has the structure of a generalized cone complex, or a cone stack over polyhedral complexes, of dimension $3g-3+n$, see \cite{AbramovichCaporasoPayne2012,CavalieriChanUlirschWise2020}. (A generalized cone complex is to a cone complex what an orbifold is to a manifold.) For each isomorphism class $\tau$ of stable $n$-leafed genus-$g$ weighted graphs, the space $\M_{g,n}^{\trop}$ has a ``folded cone'' isomorphic to $\R_{\ge0}^{\abs{E(\tau)}}/\Aut(\tau)$ (recording the possible vectors of internal edge length assignments), and the cones are glued together by identifying a graph with internal edge $e$ having length zero with the graph obtained by contracting that edge.

The \newword{pure locus} $\M_{g,n}^{\trop,\pure}\subset\M_{g,n}^{\trop}$ is the locus of graphs where all vertex weights are zero. The complement of $\M_{g,n}^{\trop,\pure},$ i.e. the set of graphs with at least one positive weight vertex, is a subcomplex of $\M_{g,n}^{\trop}$. We also define the pure locus $\Link\M_{g,n}^{\trop,\pure}=\M_{g,n}^{\trop,\pure}\cap\Link\M_{g,n}^{\trop}$; this is the quotient of $\CV_{g,n}$ by $A_{g,n}.$

\subsection{Boundary complexes} 
\label{sec:BoundaryComplexes}
The space $\Mbar_{g,n}\setminus\M_{g,n}$ is the topological boundary of $\M_{g,n}$ as a subset of $\Mbar_{g,n}$ It is a \newword{normal crossings hypersurface}, and  the \newword{boundary complex} of the pair $\Mbar_{g,n}\supseteq\M_{g,n}$ is obtained by taking a vertex for each irreducible boundary hypersurface and a $k$-simplex for each irreducible component of each $(k+1)$-fold intersection of irreducible boundary hypersurfaces, glued in the natural way \cite{AbramovichCaporasoPayne2012}. (The resulting object is a symmetric $\Delta$-complex.) The boundary complex of $\Mbar_{g,n}\supseteq\M_{g,n}$ is naturally identified with $\Link(\M_{g,n}^{\trop})$.

For $\Sgn$ an $n$-pointed genus-$g$ surface, recall that, like $\Mbar_{g,n}$, $\Tbar(\Sgn)$ has a natural boundary stratification whose locally closed strata are complex manifolds (Section \ref{sec:AugmentedTeichmullerSpace}). Unlike the compactification $\M_{g,n}\subseteq\Mbar_{g,n}$, the boundary of the partial compactification $\T(\Sgn)\subseteq\Tbar(\Sgn)$ is not a normal crossings hypersurface, since, unlike $\Mbar_{g,n},$ $\Tbar(\Sgn)$ is not even an orbifold. However, there is still a reasonable notion of boundary complex:
\begin{Def}\label{Def:BoundaryComplex}
The \emph{boundary complex} of the pair $\T(\Sgn)\subseteq\Tbar(\Sgn)$ is the simplicial complex whose vertex set is the set of codimension-1 locally closed strata in $\Tbar(\Sgn)$, and whose $k$-simplices are collections of $k$ codimension-1 strata whose closures have nonempty intersection. The boundary complex of $\T(\Sgn)\subseteq\Tbar(\Sgn)$ is naturally identified with the curve complex $\mathcal{C}(\Sgn)$.
\end{Def}

\begin{rem}\label{rem:InheritedBoundaryComplex}
    Suppose $U\subset\Tbar(\Sgn)$ is an open subset that is a union of locally closed strata $\T_\Gamma(\Sgn).$ (If $U\ne\emptyset,$ then in particular $\T
    (\Sgn)\subseteq U.$) Then $U$ defines a closed subcomplex $\mathcal{C}_U(\Sgn)\subseteq\mathcal{C}(\Sgn),$ which is the ``boundary complex of $U$'' in the same sense as Definition \ref{Def:BoundaryComplex}.

    Suppose further that $U$ is invariant under some subgroup $G\subseteq\MCG(\Sgn)$. Then $\mathcal{C}_U(\Sgn)$ is $G$-invariant. Furthermore, $U/G$ inherits a stratification from $U$, whose locally closed are naturally complex orbifolds. The boundary complex of this stratification is precisely $\mathcal{C}_U(\Sgn)/G$, again in the sense of Definition \ref{Def:BoundaryComplex}. This quotient $\mathcal{C}_U(\Sgn)/G$ has the structure of a symmetric $\Delta$-complex, but may not be a simplicial complex.
\end{rem}

\section{The mapping class group exact sequence of an \texorpdfstring{$n$}{n}-pointed handlebody}\label{sec:LuftMarkedPoints}
Let $\Xgn$ be an $n$-pointed genus-$g$ handlebody, and let $S_{g,n}=\partial\Xgn$. Recall the three groups $\MCG(\Xgn)$, $\hMCG(\Xgn)$ and $\Tw(\Xgn)$ associated to $\Xgn$, from Section \ref{sec:MappingClassGroups}. Luft proved the following:
\begin{prop}[\cite{Luft1978}]\label{prop:Luft}
    For any $g\ge1$, and $n\in\{0,1\}$, there is a natural exact sequence \begin{align}
        1\to\Tw(\Xgn)\to\MCG(\Xgn)\to\hMCG(\Xgn)\to1.
    \end{align}
\end{prop}
In this section, we prove that this sequence remains exact in the presence of any number of labeled points:
\begin{thm}\label{thm:LuftMarkedPoints}
    For any $g$ and $n$, there is a natural exact sequence
    \begin{align}\label{eq:HandlebodyGroupExactSequence}
        1\to\Tw(\Xgn)\to\MCG(\Xgn)\to\hMCG(\Xgn)\to1.
    \end{align}
\end{thm}
We will prove Theorem \ref{thm:LuftMarkedPoints} by induction on $n$, using Proposition \ref{prop:Luft} as the base case $n=0.$ We set up the induction in several steps. In what follows, assume $n>0$.

\begin{lemma}\label{lem:conjugate}
The subgroup of $\pi_1(\Sg\setminus\{x_1,\ldots,x_{n-1}\},x_n)$ generated by simple loops that bound an embedded disk is a normal subgroup.
\end{lemma}
(Note that the image of such a loop is a simple closed curve in $S_{g,n-1}$. This simple closed curve need not be essential as it may encircle a single $x_i$.)
\begin{proof}
    Let $\beta\in\pi_1(\Sg\setminus\{x_1,\ldots,x_{n-1}\},x_n)$ be a simple loop based at $x_n$ that bounds an embedded disk, and let $\alpha\in\pi_1(\Sg\setminus\{x_1,\ldots,x_{n-1}\},x_n)$ be any loop. We prove that $\gamma:=\alpha\beta\alpha^{-1}$ is homotopic to a simple closed loop that bounds an embedded disk. The proof is illustrated in Figure~\ref{fig:push}.
    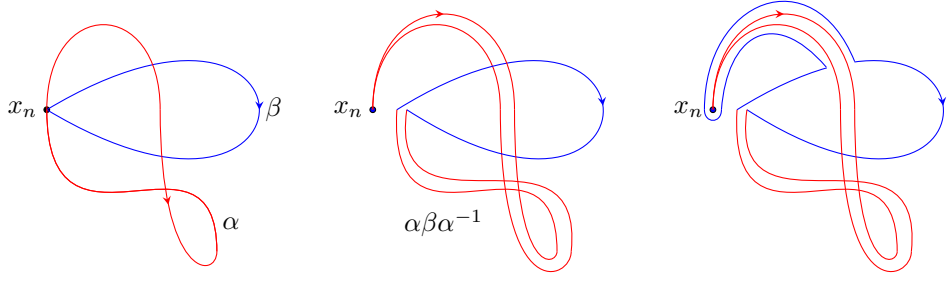
\begin{figure}
\begin{center}
\begin{tikzpicture} 
\begin{scope}[scale=.75]
\begin{scope}[decoration={markings,mark = at position 0.5 with {\arrow{stealth}}}];
\draw [fill=black] (0,0) circle(.05) node[left] {$x_n$};
  \midarrow [blue] ( 0,0) .. controls ( 5,3) and ( 5,-3) .. ( 0,0);
     \draw[red] (3,-2.5)  .. controls (3,0) and (0,-3) .. (0,0);
    \draw [red] ( 0,0) .. controls ( 0,2) and (2,2) .. ( 2,0);
    \midarrow [red] ( 2,0) .. controls ( 2,-3) and (3,-3) .. ( 3,-2.5);
    \draw[red] (3,-2.5)  .. controls (3,0) and (0,-3) .. (0,0);
 \node (a) at (3.25,-2) {$\alpha$}; 
  \node (d) at (4,0) {$\beta$};
  
   \begin{scope}[xshift = 6 cm] 
    \draw [fill=blue] (-.25,0) circle(.05) node[left] {$x_n$};
     \draw [red] ( -.25,0) .. controls ( -.25,2) and (2,2) .. ( 2,0);
    \draw[red] ( 2,0) .. controls ( 2,-3.5) and (3.25,-3) .. ( 3.25,-2.5);
    \draw[red] (3.25,-2.5)  .. controls (3.5,0) and (0,-2.5) .. (.35,0);
      \midarrow [blue] (.175,0) .. controls ( 5,3) and ( 5,-3) .. ( .35,0);
 \midarrow [red] ( -.25,0) .. controls ( -.25,2.25) and (2.25,2.25) .. ( 2.25,0);
    \draw [red] ( 2.25,0) .. controls ( 2.25,-3) and (3,-2.75) .. ( 3,-2.5);
    \draw[red] (3,-2.5)  .. controls (3,0) and (0,-3) .. (.175,0);
     
     \node (ada) at (1,-2) {$\alpha\beta\alpha^{-1}$}; 
      \end{scope}
      
      \begin{scope}[xshift = 12 cm] 
    \draw [fill=blue] (-.25,0) circle(.05) node[left] {$x_n$};
    \midarrow [blue] (.175,0) .. controls ( 5,3) and ( 5,-3) .. ( .35,0);
       \fill[white] ( 2,.8) circle (.25);
     \draw [red] ( -.25,0) .. controls ( -.25,2) and (2,2) .. ( 2,0);
    \draw[red] ( 2,0) .. controls ( 2,-3.5) and (3.25,-3) .. ( 3.25,-2.5);
    \draw[red] (3.25,-2.5)  .. controls (3.5,0) and (0,-2.5) .. (.35,0);
     
 \midarrow [red] ( -.25,0) .. controls ( -.25,2.25) and (2.25,2.25) .. ( 2.25,0);
    \draw [red] ( 2.25,0) .. controls ( 2.25,-3) and (3,-2.75) .. ( 3,-2.5);
     \draw[red] (3,-2.5)  .. controls (3,0) and (0,-3) .. (.175,0);
  
     \draw [blue] (-.1,0)  .. controls (0,1.3 ) and (.8,1.85 ) .. (1.75,.75);
     \draw [blue] (-.4,0)  .. controls (-.4,1.95) and (1.5,2.7 ) .. (2.25,.85);
       \draw [blue] (-.4,0)  .. controls (-.4,-.25) and (-.1,-.25) .. (-.1,0);
      \end{scope}
    \end{scope}
    \end{scope}
 \end{tikzpicture}
 \end{center}
 \caption{Eliminating intersections of $\alpha\beta\alpha^{-1}$ with itself (first push).}\label{fig:push}
 \end{figure}
    
    We may slightly perturb $\gamma$ within its homotopy class so that (1) it touches $x_n$ only at the start/end point, (2) $\alpha$ is replaced by two parallel arcs, and (3) there are finitely many self-intersection points of of $\gamma$, all of them transverse double points, which occur in pairs with these two parallel arcs. (See Figure~\ref{fig:push}.) Now, starting at $x_n$ and walking along $\alpha$, eliminate pairs of self-intersection points one at a time by pushing them in pairs along $\alpha$ and over $x_n$. These operations are all local and therefore can be done without pushing any curves over of the punctures $x_1,\ldots,x_{n-1}$. The result is therefore a simple loop $\gamma'$ in $\Sg\setminus\{x_1,\ldots,x_{n-1}\}$ that is homotopic to $\gamma$, hence maps trivially to $\pi_1(\Xg)$. By Dehn's lemma \cite{Papakyriakopoulos1957}, $\delta'$ bounds an embedded disk in $\Xg$, as desired.  
\end{proof}

\begin{lem}\label{lem:KernelOfAlphaGeneratedByMeridians}
    The kernel $K$ of the natural map $$\pi_1(\Sg\setminus\{x_1,\ldots,x_{n-1}\},x_n)\to\pi_1(\Xg\setminus\{x_1,\ldots,x_{n-1}\},x_n)$$ induced by inclusion is precisely the subgroup generated by simple loops that bound an embedded disk (which is normal by Lemma \ref{lem:conjugate}).
\end{lem}
    \begin{proof}
Let $\Delta=\{\delta_1,\ldots,\delta_g,\delta_1',\ldots,\delta_{n-1}'\}$ be a set of disjoint smoothly embedded disks in $\Xg\setminus\{x_1,\ldots,x_{n-1}\}$, with boundary curves $\{\gamma_1,\ldots,\gamma_g,\gamma_1',\ldots,\gamma_n'\}$, such that $\Xg\setminus\{\delta_1,\ldots,\delta_g\}$ is homeomorphic to a $3$-ball (minus $2g$ disks on the boundary) and $\delta_i'$ separates $x_i$ from all other labeled points $x_j$ and from all of the $\delta_j$s. (See Figure \ref{fig:disks}.)
By the Seifert-Van Kampen theorem, we have $$\pi_1((\Sg\setminus\{x_1,\ldots,x_{n-1}\})\cup\Delta,x_n)\cong\pi_1(\Sg\setminus\{x_1,\ldots,x_{n-1}\},x_n)/\left\langle\!\left\langle\gamma_1,\ldots,\gamma_g,\gamma_1',\ldots,\gamma_n'  \right\rangle\!\right\rangle,$$ where $\left\langle\!\left\langle\gamma_1,\ldots,\gamma_g,\gamma_1',\ldots,\gamma_n'  \right\rangle\!\right\rangle$ denotes the normal closure of the subgroup generated by $\gamma_1,\ldots,\gamma_g,\gamma_1',\ldots,\gamma_n'$. (More precisely, we must conjugate each of these simple closed curves by an arc connecting it to the basepoint, then use the method of Lemma \ref{lem:conjugate} to find a homotopy-equivalent simple loop based at $x_n$ that bounds an embedded disk.)
\begin{figure}
\begin{center}
\begin{tikzpicture} 
\begin{scope}[scale=.7]
   \draw []  plot [smooth cycle, tension=.8] coordinates 
 { (-6,0)   (-4.5 ,2)  (-2,1.25)   (0,2)   (2,1.25)   (4.5,2)  (6,0) (4.5,-2) (2,-1.25)  (0,-2) (-2,-1.25) (-4.5,-2)
    };
\draw  (.70,.1) arc (25:155:.75);
\draw   (.75,.3) arc (-10:-170:.75);
\begin{scope}[xshift=3.75cm]
\draw (.70,.1) arc (25:155:.75);
\draw (.75,.3) arc (-10:-170:.75);
\end{scope}
\begin{scope}[xshift=-3.75cm]
\draw (.70,.1) arc (25:155:.75);
\draw  (.75,.3) arc (-10:-170:.75);
\end{scope}

 \draw[thick, red, dashed, fill=red!10] (5.2,0) ellipse (.8cm and .25cm); 
 \draw[thick, red] (6,0) arc(0:180:.8cm and .25cm); 
   \draw[thick, red, dashed, fill=red!10] (0,1.25) ellipse (.25cm and .725cm);
    \draw[thick, red] (0,.55) arc(-90:90:.25cm and .725cm); 
 \draw[thick, red, dashed, fill=red!10] (-5.2,0) ellipse (.8cm and .25cm);  
  \draw[thick, red] (-4.4,0) arc(0:180:.8cm and .25cm);  
  \node [above] (d1) at  (-5,.15) {$\delta_1$};  \node (dg1) at  (-6.25,0) {$\gamma_1$};
     \node [right] (d2) at (.1,1.125) {$\delta_2$};\node (dg2) at  (0,2.25) {$\gamma_2$};
       \node [above] (d3) at  (5,.15) {$\delta_3$};\node (dg3) at  (6.3,0) {$\gamma_3$};

\fill[red!10] (4.6,1.95) arc(140:335:.5cm) arc(290:180:.6cm and .6cm);
\fill[red!20] (4.6,1.95) arc(180:290:.6cm) arc(312:160:.52cm and .52cm);
\draw[thick,red] (4.6,1.95) arc(140:335:.5cm);
\draw[thick,red,dashed] (4.6,1.95) arc(160:312:.52cm and .52cm);
\draw (5,1.6) node {$\bullet$} node[above right] {$x_1$};
\draw (4.3,1.3) node {$\delta_1'$};
\draw (4.5,2.3) node {$\gamma_1'$};
 \end{scope}
 \end{tikzpicture} 
 \end{center}
 \caption{Disks $\delta_i,\delta_i',$ labeled points $x_i$ and associated simple curves $\gamma_i,\gamma_i'$}\label{fig:disks}
 \end{figure}
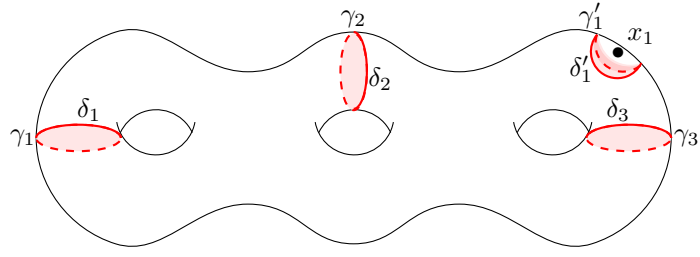

Note that $(\Sg\setminus\{x_1,\ldots,x_{n-1}\})\cup\Delta$ is homotopy equivalent to $\Sg\cup\{\delta_1,\ldots,\delta_g\}$, and $\Xg$ is obtained from $\Sg\cup\{\delta_1,\ldots,\delta_g\}$ by attaching a single $3$-ball along its boundary sphere. Thus another application of the Seifert-Van Kampen theorem shows that the inclusion $$(\Sg\setminus\{x_1,\ldots,x_{n-1}\})\cup\Delta\into\Xg\setminus\{x_1,\ldots,x_{n-1}\}$$ is the identity on $\pi_1$. Thus $$\pi_1(\Xg\setminus\{x_1,\ldots,x_{n-1}\},x_n)\cong\pi_1(\Sg\setminus\{x_1,\ldots,x_{n-1}\})/\left\langle\!\left\langle\gamma_1,\ldots,\gamma_g,\gamma_1',\ldots,\gamma_n'  \right\rangle\!\right\rangle,$$ i.e. $K=\left\langle\!\left\langle\gamma_1,\ldots,\gamma_g,\gamma_1',\ldots,\gamma_n'  \right\rangle\!\right\rangle.$
Since $K$ clearly contains every simple loop in $\pi_1(\Sg\setminus\{x_1,\ldots,x_{n-1}\},x_n)$ that bounds an embedded disk, Lemma~\ref{lem:conjugate} implies that $K\subseteq\pi_1(\Sg\setminus\{x_1,\ldots,x_{n-1}\},x_n)$ is precisely the subgroup generated by such loops.
\end{proof}
\begin{proof}[Proof of Theorem \ref{thm:LuftMarkedPoints}]
    We proceed by induction on $n$, with base case $n=0$ established in \cite[Cor. 2.5]{Luft1978}.

    Let $n>0.$ We have the following diagram:
    $$\begin{tikzcd}
        &&\MCG(\Xgn)\arrow[r,"\zeta"]\arrow[d,"\xi_1"]&\hMCG(\Xgn)\arrow[d,"\xi_2"]\\
        1\arrow[r]&\Tw(\X_{g,n-1})\arrow[r]&\MCG(\X_{g,n-1})\arrow[r]&\hMCG(\X_{g,n-1})\arrow[r]&1
    \end{tikzcd}$$
    where $\xi_1$ sends a homeomorphism to itself, $\xi_2$ sends a homotopy equivalence to itself, and $\zeta$ sends a homeomorphism to its homotopy class.
    The bottom row is exact by the inductive hypothesis, and it is easy to show that $\xi_1$ and $\xi_2$ are surjective. A straightforward diagram chase therefore gives the following diagram, where rows and columns are exact:
    $$\begin{tikzcd}
        &1\arrow[d]&1\arrow[d]&1\arrow[d]&\\
        1\arrow[r]&K\arrow[r]\arrow[d]&\ker(\xi_1)\arrow[r]\arrow[d]&\ker(\xi_2)\arrow[d]&\\
        1\arrow[r]&\ker(\zeta)\arrow[r]\arrow[d]&\MCG(\Xgn)\arrow[r,"\zeta"]\arrow[d,"\xi_1"]&\hMCG(\Xgn)\arrow[d,"\xi_2"]&\\
        1\arrow[r]&\Tw(\X_{g,n-1})\arrow[r]\arrow[d]&\MCG(\X_{g,n-1})\arrow[r]\arrow[d]&\hMCG(\X_{g,n-1})\arrow[r]\arrow[d]&1\\
        &1&1&1&
    \end{tikzcd}$$
    We first check that $\zeta$ is surjective. By Proposition \ref{prop:hMCGAddPoint}, we have canonical isomorphisms \begin{align*}
        \ker(\xi_2)&\cong\pi_1(\Xg\setminus\{x_1,\ldots,x_{n-1}\},x_n)
    \end{align*}
    and
    \begin{align*}    \ker(\xi_1)&\cong\pi_1(\Sg\setminus\{x_1,\ldots,x_{n-1}\},x_n),
    \end{align*}
    and the map $\ker(\xi_1)\to\ker(\xi_2)$ is identified with the natural map $$\pi_1(\Sg\setminus\{x_1,\ldots,x_{n-1}\},x_n)\to\pi_1(\Xg\setminus\{x_1,\ldots,x_{n-1}\},x_n)$$ induced by inclusion. It is well-known that this map is surjective, from which a straightforward diagram chase implies that $\zeta$ is surjective. We now have the diagram:
    $$\begin{tikzcd}
        &1\arrow[d]&1\arrow[d]&1\arrow[d]&\\
        1\arrow[r]&K\arrow[r]\arrow[d]&\pi_1(\Sg\setminus\{x_1,\ldots,x_n\},x_n)\arrow[r]\arrow[d]&\pi_1(\Xg\setminus\{x_1,\ldots,x_n\},x_n)\arrow[r]\arrow[d]&1\\
        1\arrow[r]&\ker(\zeta)\arrow[r]\arrow[d]&\MCG(\Xgn)\arrow[r,"\zeta"]\arrow[d,"\xi_1"]&\hMCG(\Xgn)\arrow[r]\arrow[d,"\xi_2"]&1\\
        1\arrow[r]&\Tw(\X_{g,n-1})\arrow[r]\arrow[d]&\MCG(\X_{g,n-1})\arrow[r]\arrow[d]&\hMCG(\X_{g,n-1})\arrow[r]\arrow[d]&1\\
        &1&1&1&
    \end{tikzcd}$$
    
    We wish to show that $\ker(\zeta)=\Tw(\Xgn).$ Since Dehn twists around disks in $\Xgn$ act trivially up to homotopy, we have $\Tw(\Xgn)\subseteq\ker(\zeta).$
    
    Since any disk in $\X_{g,n-1}$ can be lifted to a disk in $\Xgn,$ we have $$\xi_1(\Tw(\Xgn))=\Tw(\X_{g,n-1}).$$ It is therefore sufficient (by the correspondence theorem for groups) to show that $K\subseteq\Tw(\Xgn).$

    By Lemma \ref{lem:KernelOfAlphaGeneratedByMeridians}, $K\subseteq\pi_1(\Sg\setminus\{x_1,\ldots,x_{n-1}\},x_n)$ is generated by simple loops that bound an embedded disk. The isomorphism $\pi_1(\Sg\setminus\{x_1,\ldots,x_{n-1}\},x_n)\cong\ker(\xi_1)$ identifies such a loop $\gamma$ with the product of Dehn twists $\tw_{1}\cdot\tw_{2},$ where $\delta_1$ and $\delta_1$ are the two disks in $\Xgn$ whose boundaries are obtained by perturbing $\gamma$ away from $x_n$ in the two possible directions, and $\tw_{1},\tw_{2}$ denote the corresponding Dehn twists (see \cite[Fact 4.7]{FarbMargalit2011}, \cite[Lem. 3.3]{Hensel2020}). (If $\gamma$ is just a loop around a labeled point $x_i$, then one of these Dehn twists is trivial --- this does not cause any problems with the argument.) Thus $K\subseteq\Tw(\Xgn)$ as desired.
\end{proof}

\begin{rem}
    Since $\Xgn$ is $n$-pointed-homotopy equivalent to the $n$-thorned rose $R_{g,n}$, we have $\hMCG(\Xgn)\cong A_{g,n}.$
\end{rem}

\section{Complex handlebodies and the Teichm\"uller space of a handlebody}\label{sec:T(Vgn)}
\begin{Def}
    An $n$-pointed genus-$g$ \emph{complex handlebody} is a triple $(Y,C,\iota),$ where $Y$ is an $n$-pointed genus-$g$ handlebody, $C$ is an $n$-pointed genus-$g$ Riemann surface, and $\iota:C\to Y$ is an $n$-pointed homeomorphism of $C$ onto the boundary of $Y$. An isomorphism of complex handlebodies $(Y_1,C_1,\iota_1)$ and $(Y_2,C_2,\iota_2)$ is an $n$-pointed homeomorphism $\mu:Y_1\to Y_2$ such that $\iota_2^{-1}\circ\mu\circ\iota_1:C_1\to C_2$ is an isomorphism of $n$-pointed Riemann surfaces. 
\end{Def}
Let $\Xgn$ be an $n$-pointed genus-$g$ handlebody, and let $ \iota_{g,n}\colon\Sgn=\partial\Xgn\hookrightarrow \Xgn$  be the inclusion.  
\begin{Def}
    A \emph{complex structure} on $\Xgn$ is an $n$-pointed homeomorphism from $\Xgn$ to a complex handlebody. In other words, a complex structure on $\Xgn$ is a tuple $(Y,C,\iota,\sigma)$, where $(Y,C,\iota)$ is an $n$-pointed genus-$g$ complex handlebody, and $\sigma:\Xgn\to Y$ is an $n$-pointed homeomorphism.
\end{Def}

\begin{Def}\label{def:TX}
    The Teichm\"uller space $\T(\Xgn)$ is the set of complex structures on $\Xgn$ up to $n$-pointed homeomorphism $n$-pointed-isotopic to the identity. Explicitly, two complex structures $(Y_1,C_1,\iota_1,\sigma_1)$ and $(Y_2,C_2,\iota_2,\sigma_2)$ are considered equivalent if 
    there exists an $n$-pointed isomorphism $\mu:Y_1\to Y_2$ of complex handlebodies and an $n$-pointed homeomorphism $\tilde\alpha:\Xgn\to \Xgn$ that is $n$-pointed-isotopic to the identity, such that the diagram commutes:
    \begin{align}\label{eq:TeichmullerEquivalenceHandlebody}
        \begin{tikzcd}[ampersand replacement=\&]
            \Xgn\arrow[r,"\sigma_1"]\arrow[d,"\tilde\alpha"]\&Y_1\arrow[d,"\mu"]\\
            \Xgn\arrow[r,"\sigma_2"]\&Y_2
        \end{tikzcd}
    \end{align}
\end{Def}
There is a well-defined map $\T(\Xgn)\to\T(\Sgn)$ sending $(Y,C,\iota,\sigma)\mapsto(C,\iota^{-1}\circ\sigma\circ\iota_{g,n})$. Note that $\sigma\circ\iota_{g,n}$ lands in $\partial Y$, where $\iota^{-1}$ is well-defined. The fact that the map is well-defined is immediate from the above notion of equivalence, and the fact that an isotopy from $\tilde\alpha:\Xgn\to\Xgn$ to the identity restricts to an isotopy from $\tilde\alpha|_{\Sgn}:\Sgn\to\Sgn$ to the identity. 
\begin{prop}\label{prop:TeichmullerBijection}
    The natural map $\T(\Xgn)\to\T(\Sgn)$ is a bijection. 
\end{prop}
\begin{proof}
    Surjectivity is clear --- given a complex structure $\phi:\Sgn\to C$, the tuple $(\X,C,\iota_{g,n}\circ\phi^{-1},\id_\X)\in\T(\Xgn)$ is a complex structure on $\Xgn$ that maps to $(C,\phi)\in\T(S).$

    For injectivity, suppose we have two complex structures $(Y_1,C_1,\iota_1,\sigma_1)$ and $(Y_2,C_2,\iota_2,\sigma_2)$ on $\Xgn$ that restrict to the same complex structure on $\Sgn$. Then there is an isomorphism $\beta:C_1\to C_2$ and homeomorphism $\alpha:\Sgn\to\Sgn$ isotopic to the identity such that $\sigma_2|_{\Sgn}\circ\alpha=\beta\circ\sigma_1|_{\Sgn}$.

    We show that $\alpha$ extends to a homeomorphism $\tilde\alpha:\Xgn\to \Xgn$ that is isotopic to the identity. Let $U\supset\Sgn$ be a tubular neighborhood in $\Xgn$ homeomorphic to $\Sgn\times[0,1]$, with $\Sgn$ identified with $\Sgn\times\{0\}$. Since $\alpha$ is isotopic to the identity, we may extend $\alpha$ to a map $U\to U$ that restricts to the identity on $\Sgn\times\{1\}$. Extending by the identity on $\Xgn\setminus U$ gives a homeomorphism $\tilde\alpha:\Xgn\to \Xgn$ with a clear isotopy to the identity.

    Let $\mu=\sigma_2\circ\tilde\alpha\circ\sigma_1^{-1}.$ The restriction of $\mu$ to $C_1$ is $\beta$, which is an isomorphism of Riemann surfaces, and the diagram \eqref{eq:TeichmullerEquivalenceHandlebody} commutes by construction, showing that \begin{align*}
        (Y_1,C_1,\iota_1,\sigma_1)&=(Y_2,C_2,\iota_2,\sigma_2)\in\T(\Xgn).\qedhere
    \end{align*}
\end{proof}

From now on, we identify $\T(\Xgn)$ with $\T(\Sgn)$ via this map.

\section{Stable complex handlebodies and the augmented Teichm\"uller space of a handlebody}\label{sec:Tbar(Vgn)}
\begin{Def}\label{Def:StableComplexHandlebody}
    An \emph{stable $n$-pointed genus-$g$ complex handlebody} is a triple $(Y,C,\iota),$ where $Y$ is a finite disjoint union of handlebodies, glued together at finitely many points, together with $n$ points on the boundary of $Y$ disjoint from the locus of glued points, $C$ is an $n$-pointed genus-$g$ stable curve, and $\iota:C\to \partial Y\subset Y$ is an $n$-pointed homeomorphism of $C$ to the boundary of $Y$. An isomorphism of stable complex handlebodies is an $n$-pointed homeomorphism $\mu:Y_1\to Y_2$ such that the restriction $\iota_2^{-1}\circ\mu\circ\iota_1:C_1\to C_2$ is an isomorphism of $n$-pointed stable curves. 

    Note that an $n$-pointed genus-$g$ stable complex handlebody $(Y,C,\iota)$ has a (stable $n$-pointed genus-$g$) weighted dual graph $\tau_Y$ constructed in the same way as the dual graph of a stable curve --- in fact, the dual graph of $(Y,C,\iota)$ coincides with the dual graph of $C$. As before, there is a canonical continuous map   $\dgm_Y$, unique up to homotopy, from $Y$ to the topological  realization of $\tau_Y$. 
\end{Def}

As before, let $\Xgn$ be an $n$-pointed genus-$g$ handlebody,     $\Sgn=\partial\Xgn$ and   $\iota_{g,n}\colon\Sgn\hookrightarrow \Xgn$ the inclusion.  
\begin{Def}
    The \emph{augmented Teichm\"uller space} $\Tbar(\Xgn)$ is the set of tuples $\{(Y,C,\iota,\sigma)\}/\sim$, where $(Y,C,\iota)$ is a stable $n$-pointed genus-$g$ complex handlebody, and $\sigma:\Xgn\to Y$ is a surjective $n$-pointed map that is a homeomorphism away from some multidisk $\Delta$ in $\Xgn$, and contracts each connected component of $\Delta$ to a node of $C$. The equivalence relation is the same as for $\T(\Xgn).$
\end{Def}
\begin{prop}\label{prop:AugmentedTeichmullerInjection}
    The natural map $\Tbar(\Xgn)\to\Tbar(\Sgn)$ sending $(Y,C,\iota,\sigma)\mapsto(C,\sigma|_{\Sgn})$ is injective. 
\end{prop}
\begin{proof}
    The proof is identical to the proof of Proposition \ref{prop:TeichmullerBijection}.
\end{proof}
\begin{prop}\label{prop:AugmentedTeichmullerLocus}
    The image of $\Tbar(\Xgn)\to\Tbar(\Sgn)$ is the locus consisting of pairs $(C,\phi)\in\Tbar(\Sgn)$ such that the multicurve $\Delta$ at which $\phi$ fails to be a homeomorphism is a multimeridian.
\end{prop}
\begin{proof}
    Let $(C,\phi)\in\Tbar(\Sgn)$ be such that the multicurve $\Gamma\subset S$ at which $\phi$ fails to be a homeomorphism is a multimeridian, bounding some multidisk $\Delta$. Note that $\phi:\Sgn\to C$ induces a homeomorphism $\underline{\phi}:\Sgn/\Gamma\to C,$ where $\Sgn/\Gamma$ is obtained from $\Sgn$ by contracting each component of $\Gamma$ to a point.
    
    Let $q_\Delta:\Xgn\to \Xgn/\Delta$ denote the quotient map that contracts each component of $\Delta$ to a point. By the universal property of the quotient topology, $q_\Delta\circ\iota_{g,n}$ descends to a map $F:\Sgn/\Gamma\to \Xgn/\Delta$, since $q_\Delta\circ\iota_{g,n}:\Sgn\to \Xgn/\Delta$ is constant on each component of $\Gamma$. It is immediate that  $(\Xgn/\Delta,C,F\circ\underline{\phi}^{-1},q_\Delta)\in\Tbar(\Xgn)$ maps to $(C,\phi).$ 
    
    On the other hand, for any $(Y,C,\iota,\sigma)\in\Tbar(\Xgn),$ we know $\sigma$ fails to be a homeomorphism on a multidisk, by definition. Thus $\sigma|_{\Sgn}$ fails to be a homeomorphism on a multimeridian.
\end{proof}
\begin{rem}
    From now on, we identify $\Tbar(\Xgn)$ with its image in $\Tbar(\Sgn)$, via the injection of Proposition \ref{prop:AugmentedTeichmullerInjection}. Under this identification, $\Tbar(\Xgn)$ is a dense $\MCG(\Xgn)$-invariant open set of $\Tbar(\Sgn)$ that contains $\T(\Xgn)=\T(\Sgn)$. 
\end{rem}

\begin{rem}\label{rem:TbarVgnContractible}
    $\Tbar(\Xgn)$ is contractible, as it is a geodesically convex open subset of the $CAT(0)$ space $\Tbar(\Sgn)$, and therefore inherits a $CAT(0)$ structure \cite{Wolpert2009}.
\end{rem}

\begin{ex}[Running example, continued]\label{ex:V112}
    In Example \ref{ex:V111}, we saw that $\Tbar(S_{1,1})$ is homeomorphic to $\mathbb{H}\cup\Q\cup\{\infty\}$, where $\mathbb{H}$ is the (open) upper half-plane in $\C$. By Proposition \ref{prop:AugmentedTeichmullerLocus}, under this identification, $\Tbar(V_{1,1})$ is identified with $\mathbb{H}\cup\{\infty\}$, since $V_{1,1}$ contains a unique essential disk up to isotopy. As expected, $\mathbb{H}\cup\{\infty\}$ is invariant under the subgroup $$\MCG(\X_{1,1})\cong\left\{\pm\begin{pmatrix}
        1&0\\
        a&1
    \end{pmatrix}:a\in\Z\right\}\subset\SL_2(\Z)$$ from Example \ref{ex:V110}; in fact, as noted in Example \ref{ex:V111}, this subgroup is precisely the stabilizer of $\infty$.
\end{ex}

By definition, $\Tbar(\Xgn)$ is naturally partitioned into subsets indexed by multidisks in $\Xgn$.
\begin{Def}[Boundary strata of $\Tbar(\Xgn)$]
Let $\Delta\subset \Xgn$ be a multidisk. We denote by $\T_\Delta(\Xgn)$ the set of tuples $(Y,C,\iota,\sigma)$ such that the multidisk where $\sigma$ fails to be a homeomorphism is ($n$-pointed-)isotopic to $\Delta$. We denote by $\Tbar_\Delta(\Xgn)$ the set of tuples $(Y,C,\iota,\sigma)$ such that the multidisk where $\sigma$ fails to be a homeomorphism is isotopic to a multidisk containing $\Delta$.
\end{Def}

\begin{rem}[Intersection of Boundary Strata]
    Let $\Delta_1$ and $\Delta_2$ be two multidisks in $\Xgn$. Then $$\Tbar_{\Delta_1}(\Xgn)\cap\Tbar_{\Delta_2}(\Xgn)=\begin{cases}
        \Tbar_{\Delta_1\cup\Delta_2}(\Xgn)&\Delta_1\cup\Delta_2\text{ is a multidisk}\\
        \emptyset&\text{else}.
    \end{cases}$$
\end{rem}

The following is immediate from Proposition \ref{prop:AugmentedTeichmullerLocus} (see also Remark \ref{rem:InheritedBoundaryComplex}):
\begin{prop}
    The boundary complex of $\Tbar(\Xgn)$, in the sense of Definition \ref{Def:BoundaryComplex}, is canonically isomorphic to the disk complex $\mathcal{D}(\Xgn)$.
\end{prop}
The recursive descriptions \eqref{eq:TSBoundaryRecursion} and \eqref{eq:TbarSBoundaryRecursion} of the boundary strata of $\Tbar(\Sgn)$ restrict to recursive descriptions of the boundary strata of $\Tbar(\Xgn)$, as follows
\begin{prop}\label{prop:TBarVRecursion}
    There are natural isomorphisms
    $$\T_\Delta(\Xgn)\cong\prod_{\text{$v$ a vertex of $\tau_\Delta$}}\T(\Xgn\langle v\rangle) $$ and 
    $$\Tbar_\Delta(\Xgn)\cong\prod_{\text{$v$ a vertex of $\tau_\Delta$}}\Tbar(\Xgn\langle v\rangle),$$ where for each connected component $v$ of $\Xgn\setminus\Delta,$ $\Xgn\langle v\rangle$ denotes the pointed handlebody obtained by contracting each boundary disk in (the closure in $\Xgn$ of) $v$ to a point, and labeling that point. (If such a boundary disk is nonseparating for $v$ then this results in two labeled points.)
\end{prop}
In particular, $\T_\Delta(\X)$ is naturally a complex manifold whose codimension is the number of disks in $\Delta$.
\begin{rem}\label{rem:TbarStrataContractible}
    Recall that $\T(\Xgn)$ and $\Tbar(\Xgn)$ are contractible. (See Remark \ref{rem:TbarVgnContractible} for the latter.) Proposition \ref{prop:TBarVRecursion} therefore implies that the strata $\T_\Delta(\Xgn)$ and $\Tbar_\Delta(\Xgn)$ are contractible.
\end{rem}

\section{The moduli space of complex handlebodies}\label{sec:Mghbar}
\begin{Def}
    The \emph{moduli space of stable $n$-pointed genus-$g$ complex handlebodies} $$\Mbargnh=\{(Y,C,\iota)\}/\cong$$ is the set of stable $n$-pointed genus-$g$ complex handlebodies up to isomorphism. We define the subset $$\Mgnh=\{(Y,C,\iota)\in\Mbargnh:\text{$C$ is smooth}\}.$$
\end{Def}

For an element $(Y,C,\iota,\sigma)\in\Tbar(\Xgn)$, we may forget $\sigma$ to get a stable complex handlebody $(Y,C,\iota)\in\Mbargnh,$ giving a natural map $J_1:\Tbar(\Xgn)\to\Mbargnh$.

\begin{prop}\label{prop:UnmarkedBijection}
    The natural map $J_1:\Tbar(\Xgn)\to\Mbargnh$ descends to a bijection $$\bar{J}_1:\Tbar(\Xgn)/\MCG(\Xgn)\xrightarrow{\cong}{}\Mbargnh.$$
\end{prop}
\begin{proof}
    \noindent\textbf{The map $\bar{J}_1$ is well-defined.} Let $(Y,C,\iota,\sigma)\in\Tbar(\Xgn)$, and let $h\in\MCG(\Xgn).$ By definition of the $\MCG(\Xgn)$-action on $\Tbar(\Xgn)$, we have $$J_1((Y,C,\iota,\sigma)\cdot h)=J_1((Y,C,\iota,\sigma\circ h))=(Y,C,\iota)=J_1((Y,C,\iota,\sigma)),$$ so $J_1$ descends to a map $\bar{J}_1:\Tbar(\Xgn)/\MCG(\Xgn)\to\Mbargnh.$

    \medskip

    \noindent\textbf{Surjectivity of $\bar{J}_1$} Let $(Y,C,\iota)\in\Mbargnh,$ and let $y_1,\ldots,y_k$ be the nodes of the stable curve $C$. For each node $y_i$, delete a small open neighborhood of $\iota(y_i)\in Y.$ The result is a disjoint union of handlebodies with $2k$ boundary disks altogether. We glue pairs of boundary disks that came from the same node to get an $n$-pointed genus-$g$ handlebody $\tilde Y$, together with a multidisk $\Delta\subset\tilde Y,$ such that $\tilde Y\setminus\Delta$ is canonically (up to homotopy) homeomorphic to $Y\setminus\{\iota(y_1),\ldots,\iota(y_k)\}.$ Choose an arbitrary ($n$-pointed) homeomorphism $\sigma:\Xgn\xrightarrow{\cong}{}\tilde Y$. Then $(Y,C,\iota,\sigma)\in\bar{J}_1^{-1}(Y,C,\iota),$ i.e. $\bar{J}_1$ is surjective.

    \medskip

    \noindent\textbf{Injectivity of $\bar{J}_1$} Let $(Y_1,C_1,\iota_1,\sigma_1)$ and $(Y_2,C_2,\iota_2,\sigma_2)$ be two elements of $\Tbar(\Xgn)$ such that $$\bar{J}_1(Y_1,C_1,\iota_1,\sigma_1)=\bar{J}_1(Y_2,C_2,\iota_2,\sigma_2).$$ Then by definition of isomorphism of stable complex handlebodies, there exists a homeomorphism $\mu:Y_1\to Y_2$ such that $\iota_2^{-1}\circ\mu\circ\iota_1:C_1\to C_2$ is an isomorphism of stable curves. Consider the homeomorphism $$\alpha=\sigma_2^{-1}\circ\mu\circ\sigma_1:\Xgn\to \Xgn$$ as an element of $\MCG(\Xgn).$ Then $\alpha\cdot(Y_2,C_2,\iota_2,\sigma_2)=(Y_2,C_2,\iota_2,\sigma_2\circ\alpha)=(Y_2,C_2,\iota_2,\mu\circ\sigma_1)$. By definition of the equivalence relation on $\Tbar(\Xgn),$ we have $(Y_2,C_2,\iota_2,\mu\circ\sigma_1)=(Y_1,C_1,\iota_1,\sigma_1)$, taking $\tilde\alpha$ to be the identity in \eqref{eq:TeichmullerEquivalenceHandlebody}. Thus $(Y_1,C_1,\iota_1,\sigma_1)$ and $(Y_2,C_2,\iota_2,\sigma_2)$ are in the same $\MCG(\Xgn)$-orbit in $\Tbar(\Xgn)$. This shows $\bar{J}_1$ is injective.
\end{proof}

\begin{ex}[Running example, continued]\label{ex:V113}
    Recall from Example \ref{ex:V112} that $\Tbar(V_{1,1})\cong\mathbb{H}\cup\{\infty\}$ and $\MCG(V_{1,1})\cong\Stab_{\SL_2(\Z)}(\infty)\cong\Z\times\Z/2\Z,$ where $\Z/2\Z$ acts trivially on $\mathbb{H}\cup\{\infty\}$. Note that the quotient $$\T(V_{1,1})/\MCG(V_{1,1})\cong\mathbb{H}/\Stab_{\SL_2(\Z)}(\infty)$$ is homeomorphic to a cylinder, since $\Stab_{\SL_2(\Z)}(\infty)$ acts on $\mathbb{H}$ by integer horizontal translations. Note that as an \emph{orbifold}, the quotient $\T(V_{1,1})/\MCG(V_{1,1})$ ``remembers'' the fact that every point of $\T(V_{1,1})$ has nontrivial stabilizer isomorphic to $\Z/2\Z,$ coming from the fact that $\MCG(V_{1,1})$ does not act freely. To be precise, we should say that the \emph{underlying topological space} of the orbifold $\T(V_{1,1})/\MCG(V_{1,1})$ is homeomorphic to a cylinder.
    
    The quotient $\Tbar(V_{1,1})/\MCG(V_{1,1})\cong\mathbb{H}\cup\{\infty\}/\Stab_{\SL_2(\Z)}(\infty)$ is homeomorphic to a disk, with the image of $\infty$ at the center. (The above comment about orbifolds does not apply, as $\MCG(V_{1,1})$ acts on $\Tbar(V_{1,1})$ with infinite stabilizers, so a priori there is no orbifold structure on the quotient, only the structure of a topological space. In Section \ref{sec:ComplexStructures}, we will endow the quotient with the structure of a complex orbifold.)

\end{ex}

We note the following, which is immediate from the proof:
\begin{obs}\label{obs:BijectionRestrictsToInterior}
    The bijection in Proposition \ref{prop:UnmarkedBijection} restricts to a bijection $\Mgnh\cong\T(\Xgn)/\MCG(\Xgn).$
\end{obs}

\begin{rem}\label{rem:TopologyOfMgh}
    The bijection in Proposition \ref{prop:UnmarkedBijection} naturally endows $\Mbargnh$ with a topology. Applying Remark \ref{rem:InheritedBoundaryComplex} to the pair $(U,G)=(\Tbar(\Xgn),\MCG(\Xgn))$, we also have the following structures:
    \begin{itemize}
        \item A boundary stratification of $\Mbargnh$ by connected, locally closed strata, each with the natural structure of a complex orbifold, with strata indexed by $\MCG(\Xgn)$-orbits of boundary strata of $\Tbar(\Xgn),$ and
        \item A boundary complex of $\Mgnh\subseteq\Mbargnh$, canonically isomorphic to $\mathcal{D}(\Xgn)/\MCG(\Xgn).$
    \end{itemize}
\end{rem}
Combining this with the canonical isomorphism $\D(\Xgn)/\MCG(\Xgn)\cong\Link\M_{g,n}^{\trop}$ in \cite{Grasse1989,McCullough1991}, we have the following:
\begin{cor}\label{cor:BoundaryMghLinkMgTrop}
    The boundary complex of $\Mgnh\subseteq\Mbargnh$ is canonically isomorphic to $\Link\M_{g,n}^{\trop}$. That is, the boundary strata of $\Mbargnh$ are indexed by stable $n$-pointed genus-$g$ weighted graphs, and the stratum $\overline{\mathcal{MH}}_\tau\subseteq\Mbargnh$ corresponding to a graph $\tau$ is the locus consisting of stable complex handlebodies with dual graph $\tau$. (See Definition \ref{Def:StableComplexHandlebody}.)
\end{cor}
Proposition \ref{prop:TBarVRecursion} implies that the recursive nature of the boundary of $\Tbar(\Xgn)$ descends to $\Mbargnh$:
\begin{cor}\label{cor:MHRecursion}
    Let $\tau$ be a stable genus-$g$ $n$-leafed graph, and let $\mathcal{MH}_\tau,\overline{\mathcal{MH}}_\tau\subseteq\Mbargnh$ be the corresponding (locally closed/closed) boundary strata. Then $$\mathcal{MH}_\tau\cong\prod_{\text{$v$ a vertex of $\tau$}}\mathcal{MH}_{g_v,n_v}$$ and $$\overline{\mathcal{MH}}_\tau\cong\prod_{\text{$v$ a vertex of $\tau$}}\overline{\mathcal{MH}}_{g_v,n_v},$$ where $g_v$ is the vertex weight, and $n_v$ is the valence of $v$.
\end{cor}

\section{Homotopy-marked complex handlebodies and their moduli}\label{sec:Mgmhbar}
\begin{Def}
    A \emph{stable $n$-pointed genus-$g$ homotopy-marked complex handlebody} is a tuple $(Y,C,\iota,\eta)$, where $(Y,C,\iota)$ is a stable $n$-pointed genus-$g$ complex handlebody and $\eta:\Xgn\to Y$ is an $n$-pointed homotopy equivalence. Two stable homotopy-marked handlebodies $(Y_1,C_1,\iota_1,\eta_1)$ and $(Y_2,C_2,\iota_2,\eta_2)$ are isomorphic if there is a commutative diagram
    $$
        \begin{tikzcd}
            C_1\arrow[r,"\beta"]\arrow[d,"\iota_1"]&C_2\arrow[d,"\iota_2"]\\
            Y_1\arrow[r,"\mu"]&Y_2\\
            \Xgn\arrow[u,"\eta_1"]\arrow[ur,"\eta_2"']&
        \end{tikzcd}
    $$ where $\beta$ is an isomorphism of $n$-pointed Riemann surfaces and $\mu$ is an $n$-pointed homeomorphism, such that $\eta_2$ is $n$-pointed-homotopic to $\mu\circ\eta_1$.
\end{Def}
\begin{Def}
    The moduli space of stable $n$-pointed genus-$g$ homotopy-marked complex handlebodies is the set $$\Mbargnmh=\{(Y,C,\iota,\eta)\}/\cong$$ of stable $n$-pointed genus-$g$ homotopy-marked complex handlebodies up to isomorphism. We define the subset $$\Mgnmh=\{(Y,C,\iota,\eta)\in\Mbargnmh:\text{$C$ is smooth}\}.$$
\end{Def}

For $(Y,C,\iota,\sigma)\in\Tbar(\Xgn),$ note that the map $\sigma:\Xgn\to Y$ is a homeomorphism away from some multidisk, which it contracts. In particular, it is a homotopy equivalence. We therefore have a natural map $J_2:\Tbar(\Xgn)\to\Mbargnmh$ that sends $(Y,C,\iota,\sigma)\mapsto(Y,C,\iota,\sigma).$
\begin{prop}\label{prop:MarkedBijection}
The map $J_2:\Tbar(\Xgn)\to\Mbargnmh$ descends to a bijection $$\bar{J}_2:\Tbar(\Xgn)/\Tw(\Xgn)\xrightarrow{\cong}{}\Mbargnmh.$$    
\end{prop}
\begin{proof}
    \noindent\textbf{The map $\bar{J}_2$ is well-defined.} Let $(Y,C,\iota,\sigma)\in\Tbar(\Xgn)$, and let $h\in\Tw(\Xgn)$. By Proposition \ref{prop:Luft}, $h$ is in the kernel of the map $\MCG(\Xgn)\to\hMCG(\Xgn),$ so $h$ is $n$-pointed-homotopic (though not necessarily isotopic) to the identity. Thus $\sigma\circ h$ is $n$-pointed-homotopic to $\sigma.$ Hence $$J_2((Y,C,\iota,\sigma)\cdot h)=J_2(Y,C,\iota,\sigma\circ h)=J_2(Y,C,\iota,\sigma),$$ so $\bar{J}_2$ is well-defined.

    \medskip

    \noindent\textbf{Surjectivity of $\bar{J}_2$.} Let $(Y,C,\iota,\eta)\in\Mbargnmh$. As in the proof of surjectivity of $\bar{J}_1,$ there exists a map $\sigma:\Xgn\to Y$ that is an $n$-pointed homeomorphism away from some multidisk $\Delta$, each component of which is contracted to a node of $Y$. This map $\sigma$ is invertible up to homotopy; let $\underline{\sigma}:Y\to \X$ be an $n$-pointed homotopy inverse of $\sigma$. The composition $\underline{\sigma}\circ\eta:\Xgn\to\Xgn$ is an $n$-pointed homotopy equivalence. By Theorem \ref{thm:LuftMarkedPoints}, there exists an $n$-pointed homeomorphism $h:\Xgn\to \Xgn$ that is $n$-pointed homotopic to $\underline{\sigma}\circ\eta$. Thus $\sigma\circ h$ is homotopic to $\eta$. Because $h$ is a homeomorphism and $\sigma$ is a homeomorphism away from a multidisk, $\sigma\circ h$ is a homeomorphism away from a multidisk. Thus $J_2(Y,C,\iota,\sigma\circ h)=(Y,C,\iota,\eta)$. This shows that $J_2$, and hence $\bar{J}_2$, is surjective.

    \medskip

    \noindent\textbf{Injectivity of $\bar{J}_2$.} Let $(Y_1,C_1,\iota_1,\sigma_1)$ and $(Y_2,C_2,\iota_2,\sigma_2)$ be two elements of $\Tbar(\Xgn)$ such that $$J_2(Y_1,C_1,\iota_1,\sigma_1)=J_2(Y_2,C_2,\iota_2,\sigma_2).$$ By definition of isomorphism of stable homotopy-marked complex handlebodies, there exists an $n$-pointed homeomorphism $\mu:Y_1\to Y_2$ such that $\iota_2^{-1}\circ\mu\circ\iota_1:C_1\to C_2$ is an isomorphism of $n$-pointed stable curves and the two maps $\mu\circ\sigma_1:\Xgn\to Y_2$ and $\sigma_2:\Xgn\to Y_2$ are homotopic. Since $\sigma_2:\Xgn\to Y_2$ is a homeomorphism away from a contracted multidisk, it is in particular a homotopy equivalence. Therefore $\sigma_2$ has a homotopy inverse $\underline{\sigma_2}:Y_2\to\Xgn$. Now by construction, $\underline{\sigma_2}\circ\mu\circ\sigma_1:\Xgn\to\Xgn$ is a homotopy equivalence homotopic to the identity. By Theorem \ref{thm:LuftMarkedPoints}, it lifts to an element $h\in\MCG(\Xgn)$ homotopic to the identity, which implies (again by Theorem \ref{thm:LuftMarkedPoints}) that actually $h\in\Tw(\Xgn).$ Then $$(Y_1,C_1,\iota_1,\sigma_1)=(Y_2,C_2,\iota_2,\sigma_2\circ h)=(Y_2,C_2,\iota_2,\sigma_2)\cdot h,$$ i.e. $(Y_1,C_1,\iota_1,\sigma_1)$ and $(Y_2,C_2,\iota_2,\sigma_2)$ are in the same $\Tw(\Xgn)$-orbit. We conclude that $\bar{J}_2$ is injective.  
\end{proof}
In analogy to Observation \ref{obs:BijectionRestrictsToInterior}, we note:
\begin{obs}
    The bijection in Proposition \ref{prop:MarkedBijection} restricts to a bijection $\Mgnmh\cong\T(\Xgn)/\Tw(\Xgn).$
\end{obs}
\begin{ex}[Running example, continued]\label{ex:V114} 
    Note that $\Tw(V_{1,1})\subset\MCG(V_{1,1})$ is identified with the subgroup $$\left\{\begin{pmatrix}
        1&0\\
        a&1
    \end{pmatrix}:a\in\Z\right\}\subset\left\{\pm\begin{pmatrix}
        1&0\\
        a&1
    \end{pmatrix}:a\in\Z\right\}\cong\MCG(\X_{1,1}).$$
    By the same reasoning as in Example \ref{ex:V113}, $$h\mathcal{T}(\X_{1,1})\cong\T(\X_{1,1})/\Tw(\X_{1,1})$$ is homeomorphic to a cylinder, and $$\bar{h\mathcal{T}}(\X_{1,1})\cong\Tbar(\X_{1,1})/\Tw(\X_{1,1})$$ is homeomorphic to a disk. Unlike in Example \ref{ex:V113}, orbifold subtleties do not apply, as $\Tw(\X_{1,1})$ acts freely on $\T(\X_{1,1})$.
\end{ex}

In analogy to Remark \ref{rem:TopologyOfMgh}, we have
\begin{rem}\label{rem:TopologyOfMgmh}
    The bijection in Proposition \ref{prop:MarkedBijection} naturally endows $\Mbargnmh$ with a topology. Applying Remark \ref{rem:InheritedBoundaryComplex} to the pair $(U,G)=(\Tbar(\Xgn),\Tw(\Xgn))$, we also have the following structures:
    \begin{itemize}
        \item A boundary stratification of $\Mbargnmh$ by connected, locally closed strata, each with the natural structure of a complex orbifold\footnote{We will show in Section \ref{sec:ComplexStructures} that the strata, as well as $\Mbargnmh$ itself, are complex \emph{manifolds}.}, with strata indexed by $\Tw(\Xgn)$-orbits of boundary strata of $\Tbar(\Xgn),$ and
        \item A boundary complex of $\Mgnmh\subseteq\Mbargnmh$, canonically isomorphic to $\mathcal{D}(\Xgn)/\Tw(\Xgn).$
    \end{itemize}
\end{rem}

\subsection{The \texorpdfstring{$\hMCG(\Xgn)$}{hMCG(Vgn)}-action on \texorpdfstring{$\Mbargnmh$}{hT(Vgn)}}. 

\begin{prop}
The $\MCG(\Xgn)$-action on $\Tbar(\Xgn)$ descends to a properly discontinuous action of $$\MCG(\Xgn)/\Tw(\Xgn)\cong\hMCG(\Xgn)\cong A_{g,n}$$ on $\Mbargnmh$. The quotient $\Mbargnmh/\hMCG(\Xgn)$ is canonically identified with $\Mbargnh$, and under this identification, the quotient map $\Mbargnmh\to\Mbargnh$ coincides with the map that forgets the marking.
\end{prop}

\begin{ex}[Running example, continued]\label{ex:V115}
    The map $h\mathcal{T}(\X_{1,1})\to\mathcal{MH}_{1,1}$ is the quotient map by the \emph{trivial} action of $$\MCG(\X_{1,1})/\Tw(\X_{1,1})=\hMCG(\X_{1,1})\cong\Z/2\Z.$$ It is therefore a degree-2 covering map of orbifolds, though it is a homeomorphism of underlying topological spaces.
\end{ex}
\begin{rem}
    The cases $A_{0,n}\cong\{1\}$ (Remark \ref{rem:GenusZero}) and $A_{1,1}\cong\Z/2\Z$ (Example \ref{ex:V115}) are the only choices of $(g,n)$ for which $A_{g,n}\cong\hMCG(\Xgn)$ is finite. That is, in all other cases, the map $\Mbargnmh\to\Mbargnh$ is infinite-to-one.
\end{rem}

\section{The boundary complex of \texorpdfstring{$\Mgnmh\subset\Mbargnmh$}{hT(Vgn)}}\label{sec:MgmhbarBoundary}
In Proposition \ref{prop:MarkedBijection}, and specifically Remark \ref{rem:TopologyOfMgmh}, we noted that the boundary complex of $\Mbargnmh$ is canonically isomorphic to $\D(\Xgn)/\Tw(\Xgn)$. In this section, we will show (Theorem \ref{thm:QuotientOfDiskComplex}) that $\D(\Xgn)/\Tw(\Xgn)$ is isomorphic to $\CV_{g,n}^*,$ the simplicial completion of Culler-Vogtmann  space. (See Section \ref{sec:OuterSpaceBackground}.) The isomorphism will be canonical after choosing a $n$-pointed homotopy equivalence $R_{g,n}\to\Xgn$.

For the remainder of this section, we fix an $n$-pointed homotopy equivalence $$\rho:R_{g,n}\to\Xgn.$$ Let $\underline{\rho}:\Xgn\to R_{g,n}$ be an $n$-pointed homotopy inverse. 
 
We define a map $J_3:\D(\Xgn)\to\CV_{g,n}^*$ as follows. A point in $\D(V_{g,n})$ is a pair $\Delta=(\{\delta_i\},\{\ell_i\})$, where $\{\delta_i\}$ is a multidisk, with barycentric coordinates $\{\ell_i\}$ that sum to 1. Choose arbitrarily a pure multidisk $\Delta'\supseteq\Delta.$ We define \begin{align}\label{eq:J3Def}
    J_3(\Delta)=(\tau_{\Delta'},\dgm_{\Xgn,\Delta'}\circ\rho)\in\CV_{g,n}^*,
\end{align}
where the edge of $\tau_{\Delta'}$ corresponding to a disk $\delta_i\in\Delta$ is assigned length $\ell_i$, and all other edges (corresponding to disks in $\Delta'\setminus\Delta$) are assigned length $0$. Note that because $\Delta'$ is a pure multidisk, the dual graph is $\tau_{\Delta'}$ is not weighted (i.e. all vertex weights are zero), as required for an element of $\CV_{g,n}^*$.

\begin{prop}\label{prop:J3WellDefined}
    $J_3$ is well-defined, i.e. the point $(\tau_{\Delta'},\dgm_{\Xgn,\Delta'}\circ\rho)\in\CV_{g,n}^*$ is independent of choice of $\Delta'\supseteq\Delta$.
\end{prop}
\begin{proof}
    Let $\Delta',\Delta''$ be two pure multidisks containing $\Delta,$ with barycentric coordinates extending those of $\Delta$. We cut $\X$ along $\Delta$, and contract each ``boundary patch'' obtained this way to a single point, which we treat as a labeled point -- the result is a disjoint union of pointed handlebodies $\{\X^v\},$ indexed by vertices $v$ of $\tau_\Delta$. The multidisk $\Delta'$ defines a pure multidisk $\Delta'_v$ in each $\X^v,$ and similarly $\Delta''$ defines a pure multidisk $\Delta''_v$ in each $\X^v$.

    For each vertex $v$ of $\tau_\Delta$, $\Delta'_v$ and $\Delta''_v$ correspond to (open) simplices of $\D^{\pure}(\X^v)$. By Proposition \ref{prop:SimpleDisksContractible}, $\D^{\pure}(\X^v)$ is contractible, hence in particular it is connected. Thus there is a finite sequence of pure multidisks $\Delta'_v=\Delta_{v,1},\ldots,\Delta_{v,m}=\Delta''_v$, where $\Delta_{v,j}$ is obtained from $\Delta_{v,j-1}$ by either adding or deleting disks --- we may assume without loss of generality that $\Delta_{v,j}$ is obtained from $\Delta_{v,j-1}$ by adding or deleting a \emph{single} disk. Combining these sequences over all vertices $v$ of $\tau_\Delta$ gives a finite sequence $\Delta'=\Delta_1,\Delta_2,\ldots,\Delta_M=\Delta''$ of pure multidisks containing $\Delta$, where $\Delta_j$ is obtained from $\Delta_{j-1}$ by either adding or deleting a disk (not in $\Delta$).

    The marked dual graphs of $\Delta_j$ and of $\Delta_{j-1}$ are therefore related by the contraction, in the sense of Section \ref{sec:OuterSpaceBackground}, of a length-zero non-loop edge. These two marked graphs are therefore in the same equivalence class in the definition of $\CV_{g,n}^*$ (Section \ref{sec:OuterSpaceBackground}), i.e. they define the same point in $\CV_{g,n}^*.$ Applying this argument to each step in the sequence $\Delta'=\Delta_1,\Delta_2,\ldots,\Delta_M=\Delta''$, we conclude that the marked dual graphs of $\Delta'$ and of $\Delta''$ define the same point in $\CV_{g,n}^*.$ Thus $J_3$ is well-defined.
\end{proof}

The following properties are immediate:
\begin{itemize}
    \item $J_3$ is a \textbf{simplicial} map. 
    \item $J_3$ sends each simplex of $\D(\Xgn)$ isomorphically to a simplex of $\CV_{g,n}^*.$
    \item $J_3$ is compatible with the maps $\kappa_1:\D(\Xgn)\to\Link\M_{g,n}^{\trop}$ and $\kappa_2:\CV_{g,n}^{*}\to\Link\M_{g,n}^{\trop}$ described in Section \ref{sec:ModuliOfGraphs}, i.e. the following diagram commutes:
    \begin{align*}
        \begin{tikzcd}[ampersand replacement = \&]
            \D(\Xgn)\arrow[r,"J_3"]\arrow[dr,"\kappa_1"]\&\CV_{g,n}^{*}\arrow[d,"\kappa_2"]\\
            \&\Link\M_{g,n}^{\trop}
        \end{tikzcd}
    \end{align*}
    \item $J_3$ restricts to a map $\D^{\pure}(\Xgn)\to\CV_{g,n}$.
\end{itemize}

 \begin{thm}\label{thm:QuotientOfDiskComplex}
    $J_3$ is $\Tw(\Xgn)$-invariant, and descends to an isomorphism $$\bar{J}_3:\D(\Xgn)/\Tw(\Xgn)\to\CV_{g,n}^*.$$ 
\end{thm}
 
For the purposes of the proof, we first introduce a convenient alternate interpretation of $\CV_{g,n}^*.$ Observe that having fixed the above homotopy equivalence $\rho:R_{g,n}\to\Xgn,$ any marked graph $(\tau,r)$ (i.e. a graph together with a homotopy equivalence from $R_{g,n}$) can be converted into a ``$Xgn$-marked graph'' $(\tau,r\circ\underline{\rho})$ (i.e. a graph together with a homotopy equivalence from $\Xgn$), and similarly any $\Xgn$-marked graph may be converted into a marked graph. Therefore, we may identify $\CV_{g,n}^*$ with the space of equivalence classes of $\Xgn$-marked stable genus-$g$ $n$-leafed metric graphs, where the equivalence relation is identical to that in Section \ref{sec:OuterSpaceBackground}. Under this identification, the definition \eqref{eq:J3Def} of $J_3$ becomes simply: $$J_3(\Delta)=(\tau_{\Delta'},\dgm_{\Xgn,\Delta'})\in\CV_{g,n}^*,$$ where $\Delta'$ is chosen as above.

Under the original interpretation, $\CV_{g,n}^*$ had a natural right action by $A_{g,n}$ by precomposing the marking with a homotopy equivalence $R_{g,n}\to R_{g,n}$, as described in Section \ref{sec:OuterSpaceBackground}. Under the new interpretation, we acquire a natural right $\hMCG(\Xgn)$-action on $\CV_{g,n}^*$ by precomposing the $\Xgn$-marking by a homotopy equivalence $\Xgn\to\Xgn.$ The two interpretations, together with their respective natural actions, are compatible via the isomorphism $A_{g,n}\xrightarrow{\cong}\hMCG(\Xgn)$ that is conjugation (up to homotopy) by $\rho$.
\begin{proof}[Proof of Theorem \ref{thm:QuotientOfDiskComplex}]

    \medskip\noindent\textbf{$J_3$ is equivariant.} We have two natural right actions: of $\MCG(\Xgn)$ on $\D(\Xgn)$ and of $\hMCG(\Xgn)$ on $\CV_{g,n}^*$. As in Section \ref{sec:MappingClassGroups}, there is a natural map $\MCG(\Xgn)\to\hMCG(\Xgn)$. In what follows, for $h\in\MCG(\Xgn),$ we denote by $\bar{h}$ the image in $\hMCG$. We claim that the following diagram is equivariant:
   \begin{align*}
       \begin{tikzcd}[ampersand replacement = \&,column sep=0pt,row sep=25pt]
           \D(\Xgn)\arrow[d,"J_3"]\&\circlearrowleft\&\MCG(\Xgn)\arrow[d]\\
            \CV_{g,n}^*\&\circlearrowleft\&\hMCG(\Xgn)
       \end{tikzcd}
   \end{align*}
   In other words, we claim that given $\Delta\in\D(\Xgn)$ and $h\in\MCG(\Xgn)$, we have $$J_3(\Delta\cdot h)=J_3(\Delta)\cdot \overline{h}.$$ By continuity, and by the density of $\D^{\pure}(\Xgn)$ it suffices to prove the claim for $\Delta\in\D^{\pure}(\Xgn)$, so in what follows, we assume that $\Delta\in\D^{\pure}(\Xgn)$ and that $h\in\MCG(\Xgn)$. Then
    \begin{align*}
        J_3(\Delta\cdot h)&=(\tau_{\Delta\cdot h}, \dgm_{\Xgn,\Delta\cdot h})\\
        &=(\tau_{h^{-1}(\Delta)}, \dgm_{\Xgn,h^{-1}(\Delta)}),
    \end{align*}
    where the second equality is by definition of the $\MCG(\Xgn)$-action on $\D(\Xgn)$, see Section \ref{sec:DiskComplex}.

On the other hand 
    \begin{align*}
        J_3(\Delta)\cdot \overline{h}&=(\tau_{\Delta}, \dgm_{\Xgn,\Delta})\cdot\overline{h}\\
        &=(\tau_{\Delta}, \dgm_{\Xgn,\Delta}\circ\overline{h})\\
        &=(\tau_{\Delta}, \dgm_{\Xgn,\Delta}\circ h).
    \end{align*}

Now, observe that $h$ induces an $n$-pointed, edge-length-preserving graph isomorphism $$h_\tau:\tau_{\Xgn, h^{-1}(\Delta)}\to \tau_{\Xgn, \Delta}$$ such that the following square commutes (up to $n$-pointed homotopy):
\begin{align*}
    \begin{tikzcd}[ampersand replacement = \&]
        \Xgn\arrow[r,"h"]\arrow[d,"\dgm_{\Xgn,h^{-1}(\Delta)}",swap]\&\Xgn\arrow[d,"\dgm_{\Xgn,\Delta}"]\\
        \tau_{h^{-1}(\Delta)}\arrow[r,"h_\tau"]\&\tau_\Delta
    \end{tikzcd}
\end{align*}
(In fact, it is easy to choose representatives of the homotopy classes $\dgm_{\Xgn,h^{-1}(\Delta)}$, $h_\tau$, and $\dgm_{\Xgn,\Delta}$ so that the diagram actually commutes, not just up to homotopy.) By definition, $(\tau_{h^{-1}(\Delta)}, \dgm_{\Xgn,h^{-1}(\Delta)})$ and $(\tau_{\Delta}, \dgm_{\Xgn,\Delta}\circ h)$ represent the same point in $\CV_{g,n}$, proving the claim.

    \medskip\noindent\textbf{$J_3$ descends to $\D(\Xgn)/\Tw(\Xgn)$.} Since $J_3$ is equivariant as above, and $\Tw(\Xgn)$ is the kernel of the map $\MCG(\Xgn)\to\hMCG(\Xgn)$ by Theorem \ref{thm:LuftMarkedPoints}, we have that $J_3$ descends to a map $$\bar{J}_3:\D(\Xgn)/\Tw(\Xgn)\to\CV_{g,n}^*.$$
    Note also that $\D^{\pure}(\Xgn)\subseteq\D(\Xgn)$ is $\Tw(\Xgn)$-invariant, and $J_3$ descends to a map $\D^{\pure}(\Xgn)/\Tw(\Xgn)\to\CV_{g,n}.$

     \medskip\noindent\textbf{$\bar{J}_3$ is surjective.} The surjectivity of $J_3$ follows by combining the equivariance established above with the following:
 \begin{itemize}
     \item the quotient maps $\kappa_1:\D(\Xgn)\to\Link(\M_{g,n}^{\trop})$ and $\kappa_2\to\Link(\M_{g,n}^{\trop})$ are surjective, 
     \item $\kappa_1=\kappa_2\circ J_3$, and
     \item the homomorphism $\MCG(\Xgn)\to\hMCG(\Xgn)$ is surjective.
 \end{itemize}

In more detail: suppose we have $y\in\CV_{g,n}$. Since $\kappa_1$ is surjective, there exists $x\in\D(\Xgn)$ such that $\kappa_1(x)=\kappa_2(y)\in\Link\M_{g,n}^{\trop}$. Then $\kappa_2(J_3(x))=\kappa_1(x)=\kappa_2(y)$. Since $\kappa_2$ is the quotient map by $\hMCG(\Xgn)$, there exists $f\in\hMCG(\Xgn)$ such that $y=J_3(x)\cdot f$. Since the homomorphism $\MCG(\X_{g,n})\to\hMCG(\Xgn)$ is surjective, we can lift $f$ to $\hat{f}\in \MCG(\Xgn)$. By equivariance of $J_3$, $J_3(x\cdot \hat{f})=J_3(x)\cdot f=y$. We conclude that the image of $J_3$ contains $\CV_{g,n}\subseteq\CV_{g,n}^*$, so by density $J_3$ is surjective. Therefore $\bar{J}_3$ is surjective.

    \medskip
    
    To prove injectivity of $\bar{J}_3$, we first prove the following lemma:
    \begin{lem}\label{lem:StabilizerSurjectivity}
        Let $\Delta\in\D^{\pure}(\X)$. Then the natural homomorphism $$\Stab_{\MCG(\Xgn)}(\Delta)\to\Stab_{\hMCG(\Xgn)}(J_3(\Delta)),$$ which exists by equivariance of $J_3$, is surjective.
    \end{lem}
    \begin{proof}[Proof of Lemma \ref{lem:StabilizerSurjectivity}]
        Let $f\in\Stab_{\hMCG(\Xgn)}(J_3(\Delta))$. By definition, this means there exists an edge-length-preserving automorphism $\beta:\tau_\Delta\to\tau_\Delta$ such that the diagram commutes up to homotopy:
        \begin{equation}
            \begin{tikzcd}[column sep=40pt]
                \Xgn\arrow[d,"f"]\arrow[r,"\dgm_{\Xgn,\Delta}"]&\tau_{\Delta}\arrow[d,"\beta"]\\
                \Xgn\arrow[r,"\dgm_{\Xgn,\Delta}"]&\tau_\Delta
            \end{tikzcd}
        \end{equation}
        By \cite{CullerVogtmann1986,Hatcher1995}, there is an isomorphism $\Stab_{\hMCG(\Xgn)}(J_3(\Delta))\to\Aut(\tau_\Delta),$ where $\Aut(\tau_\Delta)$ is the group of edge-length-preserving automorphisms of $\tau_\Delta$, that sends $f\mapsto\beta.$ (In particular, $f$ is uniquely determined by $\beta$.)

        We use $\beta$ to construct a homeomorphism $\widehat{f}:\Xgn\to\Xgn$ as follows. Cutting $\Xgn$ along $\Delta$ yields a disjoint union of handlebodies $\X^v$, where $v$ runs over vertices of $\tau_\Delta$, with pairs of boundary ``patches'' corresponding to the edges of $\tau_\Delta.$ The automorphism $\beta:\tau_\Delta\to\tau_\Delta$ sends vertices to vertices and half-edges to half-edges. For each vertex $v$ of $\tau$, choose a homeomorphism $\X^v\to\X^{\beta(v)}$ that sends each boundary patch on $\X^v$ (corresponding to a half-edge $\nu$ incident to $v$) to the boundary patch on $\X^{\beta(v)}$ corresponding to $\beta(\nu).$ After possibly adjusting these homeomorphisms near each boundary patch to ensure they agree on paired patches, we may glue them to obtain a homeomorphism $\widehat{f}:\Xgn\to\Xgn$ that fixes $\Delta$ as a multidisk (with its barycentric coordinates), though it may permute the disks in $\Delta$. That is, $\widehat{f}\in\Stab_{\MCG(\Xgn)}(\Delta).$
        
        By construction, $\widehat{f}$ induces the graph automorphism $\beta:\tau_\Delta\to\tau_\Delta,$ hence maps to $f$ under the map $\MCG(\Xgn)\to\hMCG(\Xgn)$ since $f\in\Stab_{\hMCG(\Xgn)}(J_3(\Delta))$ is the unique element that induces $\beta$, as above. This proves the lemma.
    \end{proof}

\noindent\textbf{$\bar{J}_3$ is injective on the pure locus.}
    Let $\Delta_1,\Delta_2\in\D^{\pure}(\Xgn)$, such that $J_3(\Delta_1)=J_3(\Delta_2).$ Then in particular, $\kappa_1(\Delta_1)=\kappa_1(\Delta_2)$, and so there exists $h\in\MCG(\Xgn)$ such that $\Delta_2=\Delta_1\cdot h$. Let $\bar h\in\hMCG(\Xgn)$ denote the image of $h$ under the natural map $\MCG(\Xgn)\to\hMCG(\Xgn)$. Since $J_3$ is equivariant and $J_3(\Delta_1)=J_3(\Delta_2),$ we have $J_3(\Delta_1)= J_3(\Delta_2)=J_3(\Delta_1)\cdot \bar{h}$, in other words, $\bar{h}\in\Stab_{\hMCG(\Xgn)}(J_3(\Delta_1))$. By Lemma \ref{lem:StabilizerSurjectivity}, we may lift $\bar h$ to an element $h' \in\Stab_{\MCG(\Xgn)}(\Delta_1)$. Since $h$ and $h'$ have the same image in $\hMCG(\Xgn)$, $(h')^{-1}\cdot h$ is in the kernel of the homomorphism $\MCG(\Xgn)\to\hMCG(\Xgn)$, i.e. $(h')^{-1}\cdot h\in\Tw(\Xgn)$. Now, observe that
    $$\Delta_1\cdot ((h')^{-1}\cdot h)=(\Delta_1\cdot (h')^{-1})\cdot h=\Delta_1\cdot h=\Delta_2.$$

    Given two elements of $\D^{\pure}(\Xgn)$ with the same image under $J_3$, we have produced an element of $\Tw(\Xgn)$ that takes one to the other. This shows that the restriction $$\left.\bar{J_3}\right|_{\D^{\pure}(\Xgn)/\Tw(\Xgn)}:\D^{\pure}(\Xgn)/\Tw(\Xgn)\to\CV_{g,n}$$ is injective.

In order to show that $\bar{J}_3$ is injective on all of $\D(\Xgn)/\Tw(\Xgn)$, we require the following lemma:

\begin{lem}\label{lem:StrongSurjectivity}
    Let $\Sigma$ be an open simplex in $\D(\Xgn),$ and let $\Lambda$ be an open simplex of $\CV_{g,n}$ whose closure in $\CV_{g,n}^*$ contains $J_3(\Sigma).$ Then there exists an open simplex $\widehat\Lambda$ in $\D^{\pure}(\Xgn)$ such that $J_3(\widehat\Lambda)=\Lambda$ and such that $\Sigma$ is contained in the closure of $\widehat\Lambda.$
\end{lem}
\begin{proof}
Let $\Sigma$ be an open simplex in $\D(\Xgn),$ with corresponding multidisk $\Delta$ in $\Xgn.$ Let $\Lambda$ be an open simplex of $\CV_{g,n}$ whose closure in $\CV_{g,n}^*$ contains $J_3(\Sigma),$ and let $(\tau,r)$ be the $\Xgn$-marked graph corresponding to $\Lambda,$ where $r:\Xgn\to\tau$ is a homotopy equivalence. Since $J_3(\Sigma)$ is a face of $\Lambda,$ and the barycentric coordinates on $\Lambda$ are in canonical bijection with the edges of $\tau$, we may decompose the edges of $\tau$ as $$E(\tau)=E_\Delta\sqcup E',$$ where $E_\Delta$ is the set of edges whose barycentric coordinates along $J_3(\Sigma)$ are nonzero, and $E'$ is the set of edges whose barycentric coordinates along $J_3(\Sigma)$ are zero.

Cutting $\Xgn$ along $\Delta$ and contracting boundary patches, as in the proof of Proposition \ref{prop:J3WellDefined}, yields pointed handlebodies $\{\X^v\}$ indexed by the vertices of $\tau_\Delta$. For each disk $\delta\in\Delta$, we have two corresponding paired labeled points $x_\delta,x_\delta^*$, each in one of the $\X^v$s. Let $g_v$ be the genus of $\X^v$, and let $n_v$ be the number of labeled points on $\X^v$.

Similarly, let $\tau^v$ be the graphs obtained by cutting $\tau$ along $E_\Delta$. For each edge in $e\in E_\Delta,$ we have two corresponding paired labeled leaves $p_e,p_e^*,$ each in one of the $\tau^v$s.

Let $c:\tau\to\underline{\tau}$ be the map that contracts all of the edges in $E'$. (This need \emph{not} be a homotopy equivalence.) The edges of $\underline{\tau}$ are in canonical bijection with $E_\Delta.$ It follows from the definition of $J_3$ and the equivalence relation on $\CV_{g,n}^*$ that there is an isomorphism $\Phi:\underline{\tau}\to\tau_\Delta$ such that the maps $\Phi\circ c\circ r:\Xgn\to\tau_\Delta$ and $\dgm_{\Xgn,\Delta}:\Xgn\to\tau_\Delta$ are homotopic. Thus after a homotopy, we may assume that the marking $r:\Xgn\to\tau$ restricts to maps $r^v:\X^v\to\tau^v$ sending each ``cut point'' $x_\delta\in\X^v$ to the corresponding ``cut leaf'' of $p_e\in\tau^v$, where $e\in E_\Delta$ is the edge that $\Phi$ identifies with edge of $\tau_\Delta$ corresponding to the disk $\delta$.

By surjectivity of $J_3$ (as applied to the map $\D(\X^v)\to\CV_{g_v,n_v}^*$), there exists a (necessarily pure) multidisk $\Delta_v$ in $\X^v$ such that $(\tau_{\Delta_v},\dgm_{\X^v,\Delta_v})=(\tau^v,r^v)$ as elements of $\CV_{g_v,n_v}^*$. Let $$\Delta'=\Delta\cup\bigcup_{\substack{v\text{ a vertex}\\\text{of $\tau_\Delta$}}}\Delta_v,$$ where we identify $\Delta_v$ with a multidisk in $\Xgn$ in the obvious way. Then $\Delta'$ is a pure multidisk containing $\Delta.$ Thus $\Sigma$ is in the closure of the open simplex $\widehat\Sigma$ of $\D(\Xgn)$ corresponding to $\Delta'$. The marked dual graph $(\tau_{\Delta'},\dgm_{\Xgn,\Delta'})$ is obtained by gluing together all of the marked graphs $(\tau_{\Delta_v},\dgm_{\X^v,\Delta_v})$ via the pairs of labeled leaves corresponding to disks $\delta\in\Delta.$ By construction, the resulting marked graph is the same as the marked graph obtained by gluing together the marked graphs $(\tau^v,r^v)$ along the pairs of labeled leaves $p_e,p_e^*$ corresponding to edges in $E_\Delta,$ which we showed above is equivalent to $(\tau,r)$ in $\CV_{g,n}.$

It follows that if $\widehat\Sigma$ is the open simplex in $\D(\Xgn)$ corresponding to $\Delta'$, then $\Sigma$ is contained in the closure of $\widehat\Sigma$, and $J_3(\widehat\Sigma)=\Lambda$, as desired.
\end{proof}

\noindent\textbf{$\bar{J}_3$ is injective.}
    Finally, suppose $\Delta_1,\Delta_2\in\D(\Xgn)$ such that $J_3(\Delta_1)=J_3(\Delta_2).$ Pick any open simplex $\Lambda$ of $\CV_{g,n}$ whose closure contains $J_3(\Delta_1)=J_3(\Delta_2)$. By Lemma \ref{lem:StrongSurjectivity}, there exist open simplices $\widehat\Lambda_1$ and $\widehat\Lambda_2$ of $\D^{\pure}(\Xgn)$ such that $J_3(\widehat\Lambda_1)=J_3(\widehat\Lambda_2)=\Lambda$, $\Delta_1$ is in the closure of $\widehat\Lambda_1$, and $\Delta_2$ is in the closure of $\widehat\Lambda_2$. Since the restriction $$\left.\bar{J_3}\right|_{\D^{\pure}(\Xgn)/\Tw(\Xgn)}:\D^{\pure}(\Xgn)/\Tw(\Xgn)\to\CV_{g,n}$$ is injective, there exists an element $h\in\Tw(\Xgn)$ such that $\widehat\Lambda_1\cdot h=\widehat\Lambda_2.$ We then have the commutative triangle of isomorphisms of closed simplices:
    \begin{equation}
        \begin{tikzcd}
            \bar{\widehat\Lambda_1}\arrow[r,"h"]\arrow[dr,"J_3"]&\bar{\widehat\Lambda_2}\arrow[d,"J_3"]\\
            &\bar{\Lambda}
        \end{tikzcd}
    \end{equation}
    where bar denotes closure. Since $J_3(\Delta_1)=J_3(\Delta_2),$ we must have $\Delta_1\cdot h=\Delta_2$. This proves that $\bar{J_3}$ is injective.
\end{proof}
Theorem \ref{thm:QuotientOfDiskComplex} and Remark \ref{rem:TopologyOfMgmh} imply:
\begin{cor}\label{cor:MgmhBoundaryComplexCV}
    After choosing a homotopy equivalence $\rho:R_{g,n}\to\Xgn$, the boundary complex of $\Mgnmh\subseteq\Mbargnmh$ is canonically isomorphic to $\CV_{g,n}^*.$ 
\end{cor}

\begin{rem}
    This isomorphism in Corollary \ref{cor:MgmhBoundaryComplexCV}, which is given as a composition $$(\text{boundary complex of $\Mbargnmh$})\cong\D(\Xgn)/\Tw(\Xgn)\cong\CV_{g,n}^*,$$ has a simple description. Let $(Y,C,\iota,\eta)$ be a stable genus-$g$ homotopy-marked complex handlebody. The \emph{marked dual graph} of $(Y,C,\iota,\eta)$ is the marked graph $(\tau_Y,\dgm_Y\circ\eta\circ\rho)$, where $\tau_Y$ is the weighted dual graph of $(Y,C,\iota),$ and $\dgm_Y:Y\to\tau$ is as in Definition \ref{Def:StableComplexHandlebody}. Then the isomorphism in Corollary \ref{cor:MgmhBoundaryComplexCV} sends a stable $n$-pointed genus-$g$ homotopy-marked complex handlebody to (the equivalence class of) its marked dual graph.
\end{rem}

The simplices of $\CV_{g,n}^*$ are indexed by marked stable genus-$g$ $n$-leafed graphs, with some edges specified as having length zero, up to the equivalence relation in Section \ref{sec:OuterSpaceBackground}. Similarly to Corollary \ref{cor:MHRecursion}, we immediately have:
\begin{cor}\label{cor:hTRecursion}
    Let $(\tau,r)$ be an equivalence class of marked stable genus-$g$ $n$-leafed graphs, and let $\underline{\tau}$ be the stable genus-$g$ $n$-leafed graph obtained by contracting all length-zero edges and adding appropriate vertex weights, as in Section \ref{sec:ModuliOfGraphs}. Let $h\mathcal{T}_{(\tau,r)}(\Xgn),\overline{h\mathcal{T}}_{(\tau,r)}(\Xgn)\subseteq\Mbargnmh$ be the corresponding (locally closed/closed) boundary strata. Then $$h\mathcal{T}_{(\tau,r)}(\Xgn)\cong\prod_{\text{$v$ a vertex of $\underline{\tau}$}}h\mathcal{T}(\X_{g_v,n_v})$$ and $$\overline{h\mathcal{T}}_{(\tau,r)}(\Xgn)\cong\prod_{\text{$v$ a vertex of $\underline{\tau}$}}\overline{h\mathcal{T}}(\X_{g_v,n_v})$$ where $g_v$ is the vertex weight, and $n_v$ is the valence of $v$.
\end{cor}

\section{\texorpdfstring{Proof that $\Mbargnmh$ is simply connected}
{Proof that hT(Vgn) is simply connected}}

In this section we prove:
\begin{thm}\label{thm:hTBarSimplyConnected}
    $\Mbargnmh$ is simply connected for all $g$ and $n$.
\end{thm}

To simplify notation, we temporarily set $V=\Xgn$.  We will prove Theorem~\ref{thm:hTBarSimplyConnected} by writing $\overline{h\T}(V)$ as a union of two open sets, namely the ``interior" $h\TV$ and  a neighborhood $U$ of the boundary $\overline{h\T}(V)\setminus h\TV$.  We will  compute the fundamental groups  of these open sets  and of their intersection $U^o$,  and then apply  the Seifert-Van Kampen theorem together with induction on $g$ and $n$. 

Since $h\TV$ is the quotient of the contractible space $\TV$ by the free and proper action of  $\Tw(V)$, the fundamental group of $h\TV$ is isomorphic to $\Tw(V)$. In order to compute the fundamental groups of $U$ and $U^o$ we will analyze their preimages $\widetilde U\subset \TbarV$ and $\widetilde U^o\subset \TV$ and the action of the twist group $\Tw(V)$ on these lifts.

\subsection{\texorpdfstring{Definition of the sets $U$, $U^o$, $\widetilde U$ and $\widetilde U^o$.}{Definition of the sets U, Ucirc, Utilde and Ucirctilde.}}\label{sec:Notation}\label{sec:defU} Let $\Delta$ be a multidisk in $V$, and let $\Tbar_\Delta(V)\subseteq\Tbar(V)$ denote the corresponding closed stratum consisting of stable complex handlebodies $X$ where all disks in $\Delta$ (and possibly more) have been contracted. Let $V\setminus\Delta$ denote the possibly disconnected handlebody with (paired-up) boundary patches obtained by cutting $V$ along $\Delta$, and let $V\sslash\Delta$ denote the space obtained from $V$ by contracting each disk in $\Delta$ to a point. Then $V\sslash\Delta$ is a collection of pointed handlebodies, glued together at their labeled points (corresponding to contracted disks).

Given $X\in\T(V)$, the complex structure on    $\partial X$   corresponds by the Uniformization Theorem  to a unique hyperbolic metric. Given an isotopy class $\gamma$ of essential simple closed curves in $\partial V$, we define the \emph{length of $\gamma$} to be the length of the unique geodesic representative of $\gamma$ with respect to the hyperbolic metric corresponding to $X$. This defines a continuous function $$\abs{\gamma}:\T(V)\to\R_{>0}.$$ This extends to a continuous function $\abs{\gamma}:\Tbar(V)\to[0,+\infty],$ where $\abs{\gamma}^{-1}(0)$ consists of points $X$ for which $\gamma$ is contracted, and $\abs{\gamma}^{-1}(\infty)$ consists of points $X$ where $\gamma$ intersects the contracted multicurve (even after isotopy).

Fix $\epsilon$ less than the Margulis constant for hyperbolic surfaces, so isotopy classes of essential simple closed curves $\gamma$ with $\abs{\gamma}<\epsilon$ are disjoint. For a disk $\delta$ in $V$, we denote by $$\widetilde U_\delta\subseteq\Tbar(V)$$ the open subset consisting of elements $X\in\Tbar(V)$ such that $\abs{\partial\delta}<\epsilon$. We also define
\begin{align*}
    \widetilde U_\delta^{\circ}&:=\widetilde U_\delta\cap\T(V).
\end{align*}
Since $\epsilon$ is less than the Margulis constant, the Margulis Lemma implies that for nonisotopic essential disks $\delta,\delta'$ in $V$, we have $\widetilde U_\delta\cap\widetilde U_{\delta'}\ne\emptyset$ (if and) only if $\delta$ and $\delta'$ are disjoint (up to isotopy). Therefore, for any multidisk $\Delta$ in $V$, we define nonempty open subsets: 
\begin{align*}
    \widetilde U_\Delta&:=\bigcap_{\delta\in\Delta}\widetilde U_\delta\subseteq\Tbar(V)
    &&\text{and}&
    \widetilde U_\Delta^{\circ}&:=\bigcap_{\delta\in\Delta}\widetilde U_\delta^{\circ}=\widetilde U_\Delta\cap\T(V)\subseteq\T(V).
\end{align*}
That is, $\widetilde U_\Delta$ consists of points of $\Tbar(V)$ where \emph{all} disks in $\Delta$ have boundary length $<\epsilon$. We also define
\begin{align*}
    \widetilde U&:=\bigcup_{\delta}\widetilde U_\delta\subseteq\Tbar(V) &&\text{and}&
    \widetilde U^{\circ}&:=\bigcup_{\delta}\widetilde U_\delta^{\circ}=\widetilde U\cap\T(V)\subseteq\T(V)
\end{align*}
where $\delta$ runs over isotopy classes of essential disks in $V$. Note that $\widetilde U^{\circ}$ consists of the ``thin part" of $\T(V),$ i.e. points where some essential disk has boundary length $<\epsilon$, and $\widetilde U$ is an open neighborhood of the entire boundary of $\Tbar(V)$.

Observe that $\widetilde U$ and $\widetilde U^{\circ}$ are $\MCG(V)$-invariant, so in particular they are $\Tw(V)$-invariant. We define: 
\begin{align*}
    U:&=\widetilde U/\Tw(V)\subseteq\overline{h\T}(V)&&\text{and}&U^{\circ}:&=\widetilde U^{\circ}/\Tw(V)=U\cap h\T(V)\subseteq h\T(V).
\end{align*}
These are both open sets, since the quotient map by $\Tw(V)$ is an open map. (See Remark \ref{rem:OpenMaps}.)

Finally, let $\sigma$ be a simplex of $\CV_{g,n}^*$.  By Section \ref{sec:MgmhbarBoundary}, $\sigma$ corresponds to a $\Tw(V)$-orbit of multidisks. Let 
\begin{align*}
    \widetilde U_{\sigma}:&=\bigcup_{\Delta}\widetilde U_\Delta\subseteq\Tbar(V),
\end{align*}
where $\Delta$ runs over multidisks in the orbit. Then $\widetilde U_{\sigma}$ is a $\Tw(V)$-invariant open subset of $\Tbar(V)$, and we define 
\begin{align*}
    U_{\sigma}:&=\widetilde U_{\sigma}/\Tw(V)\subseteq h\TbarV.
\end{align*}
This is again open, as the quotient map is an open map. 

By Corollary \ref{cor:NotMultidisk}, if $\Delta$ and $\Delta'$ are distinct representatives of the $\Tw(V)$-orbit corresponding to $\sigma$, then $\Delta\cup\Delta'$ is not a multidisk. 
It follows $\widetilde U_\Delta\cap\widetilde U_{\Delta'}=\emptyset.$ Therefore for any multidisk $\Delta$ in the orbit $\sigma$ we have  
$$U_{\sigma}=\widetilde U_\Delta/\Stab(\Delta).$$ 
There is a natural decomposition 
\begin{align}\label{eq:StabDecomposition}
    \Stab(\Delta)\cong\Tw(\Delta)\times\Tw(V\sslash\Delta),
\end{align} where $\Tw(\Delta)\cong\Z^{\abs{\Delta}}$ is the group of Dehn twists around the disks of $\Delta$ and $\Tw(V\sslash\Delta)$ is the group generated by Dehn twists around disks disjoint from $\Delta.$

\subsection{Proof of Theorem~\ref{thm:hTBarSimplyConnected}}
We will apply the Seifert-van Kampen theorem to the open cover $h\TbarV=h\TV\cup U.$ To do so, we need to understand the fundamental groups of the three spaces $h\TV$, $U$, and $U^{\circ}=U\cap h\TV$, as well as the pushforwards along the inclusions $U^{\circ}\into h\TV$ and $U^{\circ}\into U$. We compute all of these objects via the $\Tw(V)$-equivariant topology of $\T(V),$ $\widetilde U$, and $\widetilde U^{\circ}$.
\begin{lem}\label{lem:NeighborhoodRetracts}
\begin{enumerate}
    \item For any multidisk $\Delta$ in $V,$ $\widetilde U_{\Delta}$ and $\widetilde U_{\Delta}^{\circ}$ are contractible.\label{lem:UTildeDeltaContractible}
    \item For any simplex  $\sigma$ in $\CV_{g,n}^*$, the inclusion $\overline{h\T}_{\sigma}(\X)\into U_{\sigma}$ is a homotopy equivalence.\label{lem:hTSigmaHomotopyEquivalence}
\end{enumerate}
    
\end{lem}

As we will see, the proof of \eqref{lem:hTSigmaHomotopyEquivalence} is essentially a $\Tw(\X)$-equivariant version of the proof of \eqref{lem:UTildeDeltaContractible}.
\begin{proof}
    Let $\Delta$ be a multidisk in $V.$ Let $V\setminus\Delta$ and $V\sslash\Delta$ be as in Section \ref{sec:Notation}. 

    Let $\Tbar(V\setminus\Delta)$ be the augmented Teichm\"uller space of the possibly disconnected handlebody with boundary patches obtained by cutting along $\Delta,$ where isotopies are required to fix boundary circles setwise, but not pointwise. Note the boundary patches naturally come in pairs. 
    We have a natural restriction map 
    \begin{align}\label{eq:UTildeToTbar}
        F:\widetilde U_\Delta\to\Tbar(V\setminus\Delta).
    \end{align}
    The image $\Im(F)$ is the locus in $\Tbar(V\setminus\Delta)$ 
     satisfying:
    \begin{enumerate}
        \item All boundary circles have length in the interval $[0,\epsilon)$, and 
        \item If two boundary circles correspond to the two sides of the same disk $\delta\in\Delta,$ then they have the same length.
    \end{enumerate}
    The fiber of $F$ over a point in $\Im(F)$ is isomorphic to $\R^k$, where $k$ is the number of pairs of boundary circles with \emph{nonzero length} at that point. Each factor $\R$ is a twist parameter that parametrizes all ways of gluing the corresponding paired disks. The various fibers fit together locally as follows. All fibers $\R^k$ should be thought of as quotients of $\R^{\Delta}$, where we quotient by coordinates corresponding to length-zero boundary circles.

    By \cite[Thm. 5.4]{Mondello2011}, we have a homeomorphism
    \begin{align}\label{eq:MapToProduct}
        \Im(F)
        \to[0,\epsilon)^{\Delta}\times \Tbar_\Delta(V).
    \end{align}
    Here the first factor records the lengths of the boundary circles, and the second factor intuitively collapses each boundary circle to a point --- this can be described as follows, see \cite{Mondello2011}. We identify $\Tbar(\X\setminus\Delta)$ with the space of conformal structures on the boundary of $\X\setminus\Delta$, a surface with boundary circles. There is a ``grafting map'' that glues an infinite flat cylinder to each boundary circle. The result is conformally equivalent to a (possibly disconnected) union of \emph{punctured} surfaces, with paired-up punctures corresponding to the paired-up boundary patches of $\X\setminus\Delta$. Gluing these punctures in pairs yields a point of $\Tbar_\Delta(\X).$

    (Note: The map \eqref{eq:MapToProduct} is the restriction of the map of \cite{Mondello2011} to $\widetilde U_\Delta$ --- the more complicated space $\widetilde\T$ in \cite{Mondello2011} is needed to allow arbitrarily short arcs between circles of $\partial\Delta$, which cannot happen in $\widetilde U_\Delta.$)
    
    Composing 
    the homeomorphism  \eqref{eq:MapToProduct} with the projection onto the second factor then gives a map 
    \begin{align}\label{TbarBalToTbar}
        G:\Im(F)\to \Tbar_\Delta(V).
    \end{align}
     The map $G\circ F:\widetilde U_\Delta\to\Tbar_\Delta(V)$ is a retraction of the open neighborhood $\widetilde U_\Delta\supset\Tbar_\Delta(V)$ onto the stratum. By the above descriptions of the fibers of $F$ and $G$, the fibers of $G\circ F$ are copies of $(([0,\epsilon)\times\R)/(0,t)\sim(0,t'))^{\abs{\Delta}}\cong D_\star^{\abs{\Delta}},$ where $D_\star$ is an open disk with one point on the boundary. Our above description of the way the fibers of $F$ fit together implies that the map $G\circ F$ is actually a $D_\star^{\abs{\Delta}}$-bundle. Restricting to $\widetilde U_\Delta^\circ$, we see that $\widetilde U_\Delta^\circ$ is a $D^{\abs{\Delta}}$-bundle over $\Tbar_\Delta(V),$ where $D$ is the disk obtained from $D_\star$ by removing the boundary point.

    Since $\Tbar_\Delta(V)$ is contractible (Remark \ref{rem:TbarStrataContractible}), and $D_\star$ and $D$ are contractible, the bundles $\widetilde U_\Delta$ and $\widetilde U_\Delta^\circ$ are contractible. This proves \eqref{lem:UTildeDeltaContractible}.

    To prove \eqref{lem:hTSigmaHomotopyEquivalence}, recall from Section \ref{sec:Notation} that $U_\sigma=\widetilde U_\Delta/\Stab(\Delta),$ and that $\Stab(\Delta)\cong\Tw(\Delta)\times\Tw(\X/\Delta).$ We first see how the above maps $F,G$ descend under to quotient by $\Tw(\Delta).$ Consider the intermediate space $$W_\Delta:=\widetilde U_\Delta/\Tw(\Delta),$$ which satisfies $U_\sigma=W_\Delta/\Tw(\X/\Delta).$ The map $F$ from Equation \eqref{eq:UTildeToTbar} is $\Tw(\Delta)$-invariant and descends to a map $$\underline{F}:W_\Delta\to\Im(F).$$ The fibers of this map are described just as in the proof of \eqref{lem:UTildeDeltaContractible}, except that the ways of gluing noncontracted boundary curves are now parametrized by $(S^1)^k$ instead of $\R^k$. (This is because a full twist in the gluing between a pair of boundary circles is the same as precomposition by a Dehn twist around a disk of $\Delta$, which acts trivially on $W_\Delta$.) By the same argument as before, the composition $$G\circ\underline{F}:W_\Delta\to\Tbar_\Delta(\X)$$ is a bundle with fibers isomorphic to $D_\epsilon^{\abs{\Delta}},$ where $D_\epsilon$ is a disk of radius $\epsilon$; these are precisely the quotients of the fibers $(([0,\epsilon)\times\R)/(0,t)\sim(0,t'))^{\abs{\Delta}}$ of $G\circ F$ by $\Tw(\Delta)\cong\Z^{\abs{\Delta}}$ acting on the copies of $\R$ by translation.

    We now quotient further by $\Tw(\X\sslash\Delta).$ The map $\underline{F}$ is clearly $\Tw(\X\sslash\Delta)$-equivariant. By \cite[Thm. 5.4]{Mondello2011}, the map \eqref{eq:MapToProduct}is $\Tw(\X\sslash\Delta)$-equivariant, hence so is $G$. (In fact, more is true --- they are $\MCG(\X\sslash\Delta)$-equivariant.) 
    It follows that the bundle $G\circ\underline{F}:W_\Delta\to\Tbar_\Delta(\X)$ descends to a 
    map $U_\sigma\to h\mathcal{T}_\sigma(\X)$. The $\Tw(\X\sslash\Delta)$-actions on $W_\Delta$ and $\Tbar_\Delta(\X)$ are not free --- there are (infinite) stabilizers along the locus where additional disks (other than $\Delta$) are contracted. However, if an element of $\Tw(\X\sslash\Delta)$ stabilizes a point of $\Tbar(\X\sslash\Delta)$, then it stabilizes the entire fiber of $W_\Delta\to\Tbar(\X\sslash\Delta)$, since this element acts by the identity on a neighborhood of $\Delta$. Hence the quotient $U_{\sigma}\to \overline{h\T}_{\sigma}(\X)$ inherits the structure of a $D_\epsilon^{\abs{\Delta}}$-bundle.

    As with the original fiber bundle $G\circ F,$ this last map is a retraction of the neighborhood $U_\sigma\supset h\mathcal{T}_\sigma(\X)$ onto the stratum. Since this retraction has contractible fibers, the inclusion $h\mathcal{T}_\sigma(\X)\into U_\sigma$ is a homotopy equivalence. This proves \eqref{lem:hTSigmaHomotopyEquivalence}.
\end{proof}

\begin{lem} If $g>0$, then
        $\widetilde U$ and $\widetilde U^\circ$ are contractible.\label{lem:UTildeDeformationRetract}
\end{lem}
\begin{proof}
    Consider the open cover $\widetilde U=\bigcup_{\delta}\widetilde U_\delta$, where $\delta$ runs over isotopy classes of essential disks in $\Xgn.$ By Lemma \ref{lem:NeighborhoodRetracts}\eqref{lem:UTildeDeltaContractible}, all intersections of these open sets are contractible. By the Nerve Theorem, $\widetilde U$ is homotopy equivalent to the nerve of the open cover, which is the disk complex $\D(V)$. Thus $U$ is contractible by Proposition \ref{prop:SimpleDisksContractible}.

    The exact same argument, replacing $\widetilde U_\delta$ with $\widetilde U_\delta^\circ,$ shows that $U^\circ$ is contractible.
\end{proof}

\begin{lem}\label{lem:InclusionIsomOnPi1}
    The inclusion $U^{\circ}\into h\TV$ induces an isomorphism on fundamental groups $$\pi_1(U^{\circ})\xrightarrow{\cong}{}\pi_1(h\TV).$$
\end{lem}
\begin{proof}
    Since $\Tw(V)$ acts freely on $\T(V)$, and $\widetilde U^{\circ}$ is $\Tw(V)$-invariant, we have a restriction of covering maps:
    \begin{align*}
        \begin{tikzcd}[ampersand replacement=\&]
            \widetilde U^{\circ}\arrow[r]\arrow[d]\&\T(V)\arrow[d]\\
            U^{\circ}\arrow[r]\&h\TV
        \end{tikzcd}
        \end{align*}
        Since $\widetilde U^{\circ}$ and $\T(V)$ are contractible, the result follows from elementary covering space theory.\qedhere
\end{proof}

\begin{proof}[Proof of Theorem \ref{thm:hTBarSimplyConnected}]
    We proceed by induction on $g$ and $n$, with base cases $g=0$ (for all $n\ge3$). 
    Indeed, by Remark \ref{rem:GenusZero}, $\bar{h\mathcal{T}}(V_{0,n})\cong\Mbar_{0,n},$ which is simply connected \cite{BoggiPikaart2000}.

    Now fix $(g,n)$ with $g\ge1$, and suppose that $\overline{h\T}(V_{g',n'}
    )$ is simply connected for all $(g',n')$ with $g'<g$ or $n'<n$.  
     Consider the  open cover $\Mbargnmh=\Mgnmh\cup U,$ where $U$ is the neighborhood of $\partial\Tbar(V_{g,n})$ defined in Section~\ref{sec:defU}.   The Seifert-van Kampen theorem implies that $\pi_1(\Mbargnmh)$ is the amalgamated free product:
    \begin{align}\label{eq:AmalgamatedFreeProduct}
        \pi_1(\Mbargnmh)\cong\pi_1(\Mgnmh)*_{\pi_1(U^{\circ})}\pi_1(U).
    \end{align}

We claim that   $U$   is simply connected.  To see this, 
        consider the open cover of $U$ by the sets $U_{\sigma_0}$, where $\sigma_0$ runs over \emph{vertices} of $\CV_{g,n}^*$. 
        Any intersection of such sets is equal to $U_{\sigma}$ for some simplex $\sigma$   of $CV_{g,n}^*$. By Corollary \ref{cor:hTRecursion}, $\bar{h\T}_\sigma(V_{g,n})$ is isomorphic to a product of spaces $\bar{h\T}_{g',n'}$ with $3g'-3+n'<3g-3+n$, so in particular $g'<g$ or $n'<n.$ By assumption, $\bar{h\T}_{g',n'}$ is simply connected. Thus $\bar{h\T}_\sigma(V)$ is simply connected. By Lemma \ref{lem:NeighborhoodRetracts}\eqref{lem:hTSigmaHomotopyEquivalence}, we have an isomorphism $\pi_1(\bar{h\T}_\sigma(V))\to\pi_1(U_\sigma)$, so $U_\sigma$ is simply connected.
       Thus $\pi_1(U)$ is isomorphic to the fundamental group of the nerve of the open cover, which is precisely $\CV_{g,n}^*$. Since $g\ge1,$ $\CV_{g,n}^*$ is contractible (Section \ref{sec:OuterSpaceBackground}), so in particular $\pi_1(\CV_{g,n}^*)$ is trivial.

      Since $U$ is simply-connected, \eqref{eq:AmalgamatedFreeProduct} simplifies to:
    \begin{align}\label{eq:AmalgamatedFreeProduct2}
        \pi_1(\Mbargnmh)\cong\pi_1(\Mgnmh)/\NCl(\pi_1(U^{\circ})),
    \end{align}
    where $\NCl$ denotes normal closure of the image. By Lemma \ref{lem:InclusionIsomOnPi1}, $\NCl(\pi_1(U^{\circ}))=\pi_1(\Mgnmh)$, so \eqref{eq:AmalgamatedFreeProduct2} implies $\pi_1(\Mbargnmh)$ is trivial.
\end{proof}

\begin{rem}\label{rem:HTbarNotContractible}
    Though $\Mbargnmh$ is simply connected for all $g,n$, it is not contractible unless $(g,n)\in\{(0,3),(1,1)\}.$ To see this, first note that $\Mbar_{0,n}$ is not contractible if $n\ge4.$ If $(g,n)\not\in\{(0,3),(1,1)\},$ then there exists a stratum $S\subseteq\Mbargnmh$ isomorphic to a product of copies of $\Mbar_{0,m}$ for various $m$, with at least one factor with $m\ge4.$ We claim that the fundamental class $[S]\in H_*(\Mbargnmh)$ is nonzero. The pushforward of $[S]$ in $H_*(\Mbar_{g,n})$ is, up to a nonzero multiple, the class of a boundary stratum in $\Mbar_{g,n}$. Since $\Mbar_{g,n}$ is projective, the fundamental class of any subvariety is nontrivial. Thus $[S]\ne0.$
\end{rem}

\section{Tropicalization maps \`a la Chan-Galatius-Payne}\label{sec:CGP}
In \cite{ChanGalatiusPayne2019}, the authors define a canonical (up to homotopy) continuous surjection $\lambda_{g,n}:\M_{g,n}\to\M_{g,n}^{\trop}$, as follows. By the Uniformization Theorem, there is a bijection between complex structures and hyperbolic metrics on a given surface, so given $(C,p_1,\ldots,p_n)\in\M_{g,n},$ we may endow it with a unique hyperbolic metric. If $\epsilon>0$ is the Margulis constant for the hyperbolic plane, then simple closed curves of length less than epsilon are all disjoint, so form a multicurve $\Gamma$ (see e.g.  \cite[Cor. 13.7]{FarbMargalit2011}).  If $\len(\gamma_e)$ is the length of the simple closed curve $\gamma\in\Gamma,$ the graph $\lambda_{g,n}(C,p_1,\ldots,p_n)$ is defined to be the dual graph of $\Gamma$, where the length of the edge $e$ corresponding to $\gamma$ is $-\log(\len(\gamma_e)/\epsilon)$.

\begin{rem}
    As noted in \cite[Sec. 3.4]{ramadas2024thurston}, this construction naturally defines an $\MCG(\Sgn)$-equivariant surjection $\tilde\lambda_{g,n}:\T(\Sgn)\to\Cone(\mathcal{C}(\Sgn)),$ where $\Cone(\mathcal{C}(\Sgn))$ denotes the cone over $\mathcal{C}(\Sgn)$, from which the map of \cite{ChanGalatiusPayne2019} is descended.
\end{rem}

We now construct analogous maps $\lambda_{g,n}^{\mathcal{MH}}:\Mgnh\to\M_{g,n}^{\trop}$ and $\lambda_{g,n}^{h\mathcal{T}}:\Mgnmh\to\Cone(\CV_{g,n}^*)$.

\begin{Def}
    Let $\epsilon>0$ be smaller than the Margulis constant. Given $(Y,C,\iota,\sigma)\in\T(\X)$, consider the hyperbolic metric on $C\setminus\{p_1,\ldots,p_n\}$ corresponding to the complex structure. The set of disks in $Y$ whose boundaries have hyperbolic length at most $\epsilon$ is a multidisk $\Delta$, by the Margulis Lemma. We define $\tilde\lambda_{\Xgn}(Y,C,\iota,\sigma)\in\Cone(\D(\Xgn))$ to be the weighted multidisk $\Delta$, with the weight of a disk $\delta\in\Delta$ equal $-\log(\len(\partial\delta)/\epsilon)$, where $\len(\partial\delta)$ is the hyperbolic length of the boundary of $\delta$.
\end{Def}
The proof of the following proposition is similar to the proof of \cite[Prop. 4.1]{HainautPetersen2025}:
\begin{prop}\label{prop:CGPAnalogy}
    The map $\tilde\lambda_{\Xgn}:\T(\Xgn)\to\Cone(\D(\Xgn))$ is continuous, surjective, and $\MCG(\Xgn)$-equivariant. In particular, by Propositions \ref{prop:UnmarkedBijection} and \ref{prop:MarkedBijection}, $\tilde\lambda_{\Xgn}$ descends to continuous, surjective maps $$\lambda_{g,n}^{h\mathcal{T}}:\Mgnmh\to\Cone(\CV_{g,n}^*)$$ and $$\lambda_{g,n}^{\mathcal{MH}}:\Mgnh\to\M_{g,n}^{\trop}.$$
\end{prop}

\section{Complex structures on certain quotients of Teichm\"uller space}\label{sec:ComplexStructuresTechnical}
In this section, we construct complex structures on quotients of certain open subsets of $\Tbar(\Sgn)$ by subgroups of $\MCG(\Sgn)$. There is an existing well-known example of such a construction, namely the complex structure on $\Mbar_{g,n}$ via ``plumbing coordinates'' \cite{Marden1987,Masur1976,Wolpert1990}. In order to generalize this construction, we recall the basic setup.

Let $\Gamma$ be a multicurve in $\Sgn$, and let $$U_\Gamma=\bigsqcup_{\Gamma'\subseteq\Gamma}\T_{\Gamma'}(\Sgn)\subseteq\Tbar(\Sgn),$$ where $\T_{\Gamma'}(\Sgn)$ is as in Section \ref{sec:AugmentedTeichmullerSpace}.
Then $U_\Gamma$ is an open subset of $\Tbar(\Sgn)$ containing $\T_\Gamma(\Sgn)$ (but not $\Tbar_\Gamma(\Sgn)$). $\MCG(\Sgn)$ induces maps between open sets $U_\Gamma$ --- specifically, if $h\in\MCG(\Sgn),$ then $h$ induces a homeomorphism $h:U_\Gamma\to U_{\Gamma\cdot h}$ for any $\Gamma.$ 

Let $$\Tw(\Gamma)\subseteq\MCG(\Sgn)$$ denote the subgroup generated by Dehn twists around the simple closed curves in $\Gamma.$ Then $U_\Gamma$ is invariant under the action of $\Tw(\Gamma)$. Let $$Q_\Gamma=U_\Gamma/\Tw(\Gamma)$$ denote the topological quotient space, with quotient map $t_\Gamma:U_\Gamma\to Q_\Gamma$. The action of $\MCG(\Sgn)$ descends to $Q_\Gamma$, since for any $h\in\MCG(\Sgn),$ the above homeomorphism $h:U_\Gamma\to U_{\Gamma\cdot h}$ is equivariant with respect to the $\Tw(\Gamma)$-action on $U_\Gamma$ and the $\Tw(\Gamma\cdot h)$-action on $U_{\Gamma\cdot h}$, hence induces a homeomorphism $$h:Q_\Gamma\to Q_{\Gamma\cdot h}.$$ Observe that if $\Gamma'\subseteq\Gamma$ are two multicurves, then we have an induced map $t_{\Gamma',\Gamma}:Q_{\Gamma'}\to Q_{\Gamma}$ since the composition $$U_{\Gamma'}\into U_{\Gamma}\xrightarrow{t_\Gamma}{}Q_\Gamma$$ is $\Tw(\Gamma')$-invariant. (Note $t_\Gamma=t_{\emptyset,\Gamma}.$) Furthermore, for any $\Gamma,$ the inclusion $\Tw(\Gamma)\into\MCG(\Sgn)$ induces a map $$q_{\Gamma}:Q_\Gamma\to\Mbar_{g,n},$$ which commutes with system of maps $Q_{\Gamma'}\to Q_{\Gamma}$. We have the following, due to Hubbard and Koch:
\begin{thm}[{\cite[Thms. 10.1, 13.1]{HubbardKoch2014}}]\label{thm:HubbardKoch}
For every multicurve $\Gamma$ on $\Sgn$, $Q_\Gamma$ has the natural structure of a complex manifold. With respect to this structure,
\begin{enumerate}
    \item The maps $t_{\Gamma',\Gamma}:Q_{\Gamma'}\to Q_{\Gamma}$ and $q_\Gamma:Q_\Gamma\to\Mbar_{g,n}$ are \'etale maps of complex orbifolds,\label{it:tGammaGammaEtale}
    \item For any multicurve $\Gamma'\subseteq\Gamma,$ the stratum $\T_{\Gamma'}(\Sgn)/\Tw(\Gamma)\subseteq Q_\Gamma$ is a locally closed complex submanifold, and the map $\T_{\Gamma'}(\Sgn)\to\T_{\Gamma'}(\Sgn)/\Tw(\Gamma)$ is a holomorphic covering map, \label{it:QGammaStratumSubmanifold}
    \item As $\Gamma'$ runs over multicurves contained in $\Gamma,$ the submanifolds $\T_{\Gamma'}(\Sgn)/\Tw(\Gamma)\subseteq Q_\Gamma$ meet at simple normal crossings, and\label{it:QGammaStrataSNC}
    \item For any $h\in\MCG(\Sgn),$ the induced maps $h:Q_\Gamma\to Q_{\Gamma\cdot h}$ are biholomorphisms.\label{it:QGammaInducedActionBiholomorphic}
\end{enumerate}

\end{thm}
\begin{rem}
    Hubbard and Koch used Theorem \ref{thm:HubbardKoch} to give a natural functorial identification of $\Tbar(\Sgn)/\MCG(\Sgn)$ with the moduli space $\Mbar_{g,n}$ of stable curves.
\end{rem}
We now adapt the construction of the complex structure on $\Mbar_{g,n}$ to a more general setting, as follows. For a point $x\in\Tbar(\Sgn),$ let $\Gamma(x)$ denote the associated contracted multicurve in $\Sgn$.
\begin{lem}\label{lem:ComplexOrbifold}
    Let $G\subseteq\MCG(\Sgn)$ be a subgroup, and define $$U(G)=\{x\in\Tbar(\Sgn):\Tw(\Gamma(x))\subseteq G\}\subseteq\Tbar(\Sgn).$$ Then
    \begin{enumerate}
        \item $U(G)$ is a $G$-invariant open subset of $\Tbar(\Sgn)$.\label{it:InvariantOpens}
        \item The quotient $U(G)/G$ admits the natural structure of a complex orbifold, with respect to which the natural map $U(G)/G\to\Mbar_{g,n}$ induced by the inclusion $G\subseteq\MCG(\Sgn)$ is an \'etale morphism of complex orbifolds.\label{it:U1G}
        \item \label{it:IsotropyU1} For $x\in U(G)$ with image $\bar{x}\in U(G)/G,$ the isotropy group at $\bar{x}$ is canonically isomorphic to $\Stab_G(x)/\Tw(\Gamma(x)).$
        \item \label{it:StratumSuborbifold} If $\Gamma$ is a multicurve with $\Tw(\Gamma)\subseteq G,$ then the image $\T_\Gamma(\Sgn)/G\subseteq U(G)/G$ is a locally closed complex suborbifold, and the map $\T_\Gamma(\Sgn)\to\T_\Gamma(\Sgn)/G$ is a holomorphic covering map.
    \end{enumerate}  
\end{lem}
\begin{proof}
    If $x\in\Tbar(\X)$ and $h\in G,$ then $\Tw(\Gamma(x\cdot h))=\Tw(\Gamma(x)\cdot h)=h^{-1}\Tw(\Gamma(x))h,$ which shows that $U(G)$ is $G$-invariant. Furthermore we have $$U(G)=\bigsqcup_{\Gamma:\Tw(\Gamma)\subseteq G}\T_\Gamma(\Sgn)=\bigcup_{\Gamma:\Tw(\Gamma)\subseteq G}U_\Gamma,$$ so $U(G)$ is open. This proves \eqref{it:InvariantOpens}

    We now construct the structure of a complex orbifold on $U(G)/G.$ Let $$\pi_G:U(G)\to U(G)/G$$ denote the quotient map, and for $\Gamma$ such that $\Tw(\Gamma)\subseteq G$, let $\pi_{\Gamma,G}:=\pi_G|_{U_\Gamma}$ denote its restriction to $U_\Gamma$. Let $\beta_G:U(G)/G\to\Mbar_{g,n}$ denote the map induced by the inclusion $G\into\MCG(\Sgn)$. Let $P\in U(G)/G$ and let $\underline{P}=\beta_G(P).$ Fix a lift $\widetilde{P}\in U(G)$ with $\pi_G(\widetilde{P})=P.$ Then $$\widetilde{P}\in U_{\Gamma(\widetilde{P})}\subseteq U(G).$$ Let $\widehat P=t_{\Gamma(\widetilde{P})}(\widetilde{P}).$ That is, we have the diagram:
    $$\begin{tikzcd}[column sep=50pt]
        \widetilde{P}\arrow[r,|->]\arrow[d,phantom,sloped,"\in"]&\widehat{P}\arrow[r,|->]\arrow[d,phantom,sloped,"\in"]&P\arrow[r,|->]\arrow[d,phantom,sloped,"\in"]&\underline{P}\arrow[d,phantom,sloped,"\in"]\\
        U_{\Gamma(\widetilde{P})}\arrow[r,"t_{\Gamma(\widetilde{P})}"]\arrow[rr,bend right,"\pi_{\Gamma(\widetilde{P}),G}",swap]&Q_{\Gamma(\widetilde{P})}\arrow[r,"q_{\Gamma(\widetilde{P}),G}"]\arrow[rr,bend right,swap,"q_{\Gamma(\widetilde{P})}",crossing over]&U(G)/G\arrow[r,"\beta_G"]&\Mbar_{g,n}
    \end{tikzcd}$$

    Let $\Isot(\underline{P})$ denote the (finite) automorphism group of the stable curve $\underline{P}.$ (This is precisely the isotropy group of the point $\underline{P}\in\Mbar_{g,n}$.) By Lemma \ref{lem:StabilizerOfPointInTbarExactSequence}, we have $$\Isot(\underline{P})\cong\frac{\Stab_{\MCG(\Sgn)}(\widetilde{P})}{\Tw(\Gamma(\widetilde{P}))}.$$

    By Theorem \ref{thm:HubbardKoch}(1), the map $q_{\Gamma(\widetilde{P})}:Q_{\Gamma(\widetilde{P})}\to\Mbar_{g,n}$ is an \'etale morphism from a complex manifold to a complex orbifold. Thus there exists a neighborhood $W(P,\widetilde{P})\subseteq Q_{\Gamma(\widetilde{P})}$ of $\widehat P$, together with an $\Isot(\underline{P})$-action on $W(P,\widetilde{P})$ by biholomorphisms fixing $\widehat P$, such that the restriction $$q_{\Gamma(\widetilde{P})}|_{W(P,\widetilde{P})}:W(P,\widetilde{P})\to\Mbar_{g,n}$$ is the quotient map by $\Isot(\underline{P})$, followed by the inclusion of an open set. 
    
    Since the factorization $$q_{\Gamma(\widetilde{P})}=\beta_G\circ q_{\Gamma(\widetilde{P}),G}$$ is induced by the inclusions of groups $$\Tw(\Gamma(\widetilde{P}))\subseteq G\subseteq\MCG(\Sgn),$$ the restriction $$q_{\Gamma(\widetilde{P}),G}|_{W(P,\widetilde{P})}:Q_{\Gamma(\widetilde{P})}\to U(G)/G$$ is simply the quotient map by the subgroup
    \begin{align}\label{eq:IsotropyGroupAtP}
        \Isot(P):=\frac{\Stab_G(\widetilde{P})}{\Tw(\Gamma(\widetilde{P}))}\subseteq\frac{\Stab_{\MCG(\Sgn)}(\widetilde{P})}{\Tw(\Gamma(\widetilde{P}))},
    \end{align} again followed by the inclusion of an open set. ($\Isot(P)$ is finite by Lemma \ref{lem:StabilizerOfPointInTbarExactSequence} and the fact that $\Stab_G(x)\subseteq\Stab_{\MCG(\Sgn)}(x)$.) This defines the structure of a complex orbifold on a neighborhood $W(P,\widetilde{P})/\Isot(P)$ of $P\in U(G)/G$. This construction works for any $P\in U(G)/G$, giving an open cover by complex orbifolds. 

    It still remains to check that the complex structures on the above open cover are locally compatible, but we observe now that once this is completed, part \eqref{it:IsotropyU1} of the Theorem will be immediate from Equation \eqref{eq:IsotropyGroupAtP}, and part \eqref{it:StratumSuborbifold} will be immediate from Theorem \ref{thm:HubbardKoch}\eqref{it:QGammaStratumSubmanifold}.

    We now check that the complex structures on the open cover are locally compatible. Let $$W(P_1,\widetilde{P}_1)/\Isot(P_1)\hspace{1in}\text{and}\hspace{1in}W(P_2,\widetilde{P}_2)/\Isot(P_2)$$ be two open sets of the above form, with $W(P_1,\widetilde{P}_1)\subseteq Q_{\Gamma(\widetilde{P}_1)}$ and $W(P_2,\widetilde{P}_2)\subseteq Q_{\Gamma(\widetilde{P}_2)}$.

    Let $$P^0\in W(P_1,\widetilde{P}_1)/\Isot(P_1)\cap W(P_2,\widetilde{P}_2)/\Isot(P_2).$$ Choose lifts 
    $$\begin{tikzcd}[row sep=10pt]
        \widetilde{P}^0_1\arrow[r,|->]\arrow[d,phantom,sloped,"\in"]&\widehat{P}^0_1\arrow[r,|->]\arrow[d,phantom,sloped,"\in"]&P^0\arrow[d,phantom,sloped,"\in"]\\
        U_{\Gamma(\widetilde{P}_1)}\arrow[r]&Q_{\Gamma(\widetilde{P}_1)}\arrow[r]&U(G)/G
    \end{tikzcd}$$
    and 
    $$\begin{tikzcd}[row sep=10pt]
        \widetilde{P}^0_2\arrow[r,|->]\arrow[d,phantom,sloped,"\in"]&\widehat{P}^0_2\arrow[r,|->]\arrow[d,phantom,sloped,"\in"]&P^0\arrow[d,phantom,sloped,"\in"]\\
        U_{\Gamma(\widetilde{P}_2)}\arrow[r]&Q_{\Gamma(\widetilde{P}_2)}\arrow[r]&U(G)/G
    \end{tikzcd}$$
    with $\widehat{P}^0_1\in W(P_1,\widetilde{P}_1)$ and $\widehat{P}^0_2\in W(P_2,\widetilde{P}_2)$. In order to show that the complex structures $$W(P_1,\widetilde{P}_1)/\Isot(P_1)\into U(G)/G\hspace{.5in}\text{and}\hspace{.5in}W(P_2,\widetilde{P}_2)/\Isot(P_2)\into U(G)/G$$ are locally compatible at $P^0,$ we must construct a biholomorphism from a neighborhood of $\widehat{P}^0_1\in W(P_1,\widetilde{P}_1)$ to a neighborhood of $\widehat{P}^0_2\in W(P_2,\widetilde{P}_2)$ that is equivariant with respect to the action of $\Stab_{\Isot(P_1)}(\widehat{P}^0_1)\cong\Isot(P^0)$ on $W(P_1,\widetilde{P}_1)$ and the action of $\Stab_{\Isot(P_2)}(\widehat{P}^0_2)\cong\Isot(P^0)$ on $W(P_2,\widetilde{P}_2)$.

    Since $\widetilde{P}^0_1$ and $\widetilde{P}^0_2$ both map to $P^0\in U(G)/G$, there exists an element $h\in G$ such that $\widetilde{P}^0_1\cdot h=\widetilde{P}^0_2$. By Theorem \ref{thm:HubbardKoch}(4), $h$ induces a biholomorphism $h_Q:Q_{\Gamma(\widetilde{P}^0_1)}\to Q_{\Gamma(\widetilde{P}^0_2)}.$ We therefore have a diagram $$\begin{tikzcd}
        U_{\Gamma(\widetilde{P}_1)}\arrow[ddd]&&&&U_{\Gamma(\widetilde{P}_2)}\arrow[ddd]\\
        &U_{\Gamma(\widetilde{P}^0_1)}\arrow[ul,hook']\arrow[rr,"\text{action of }h","\cong"']\arrow[d,"t_{\Gamma(\widetilde{P}^0_1)}"']&&U_{\Gamma(\widetilde{P}^0_2)}\arrow[ur,hook]\arrow[d,"t_{\Gamma(\widetilde{P}^0_2)}"]&\\
        &Q_{\Gamma(\widetilde{P}^0_1)}\arrow[rr,"h_Q","\cong"']\arrow[dl,"\text{\'et}"',"t_{\Gamma(\widetilde{P}^0_1),\Gamma(\widetilde{P}_1)}"]&&Q_{\Gamma(\widetilde{P}^0_2)}\arrow[dr,"\text{\'et}","t_{\Gamma(\widetilde{P}^0_2),\Gamma(\widetilde{P}_2)}"']&\\
        Q_{\Gamma(\widetilde{P}_1)}&&&&Q_{\Gamma(\widetilde{P}_2)}\\
        W(P_1,\widetilde{P}_1)\arrow[u,hook]\arrow[d]&&&&W(P_2,\widetilde{P}_2)\arrow[u,hook]\arrow[d]\\
        W(P_1,\widetilde{P}_1)/\Isot(P_1)\arrow[rr,hook]&&U(G)/G&&W(P_2,\widetilde{P}_2)/\Isot(P_2)\arrow[ll,hook']
    \end{tikzcd}$$
    By Theorem \ref{thm:HubbardKoch}\eqref{it:tGammaGammaEtale} (applied with $\Gamma'=\emptyset$ and $\Gamma=\Gamma(\widetilde{P}^0_1)$), $t_{\Gamma(\widetilde{P}^0_1)}$ is \'etale, so there exists an open set $\widetilde W^0_1\subseteq Q_{\Gamma(\widetilde{P}^0_1)}$ containing $t_{\Gamma(\widetilde{P}^0_1)}(\widetilde{P}^0_1)$ such that the restriction $t_{\Gamma(\widetilde{P}^0_1),\Gamma(\widetilde{P}_1)}|_{\widetilde W^0_1}$ is a biholomorphism onto a $\Stab_{\Isot(P_1)}(\widehat{P}^0_1)$-invariant neighborhood $W^0_1$ of $\widehat{P}^0_1$ in $W(P_1,\widetilde{P}_1)$. Similarly, there is an open set $\widetilde W^0_2\subseteq Q_{\Gamma(\widetilde{P}^0_2)}$ containing $t_{\Gamma(\widetilde{P}^0_2)}$ such that the restriction $t_{\Gamma(\widetilde{P}^0_2),\Gamma(\widetilde{P}_2)}|_{\widetilde W^0_2}$ is a biholomorphism onto a $\Stab_{\Isot(P_2)}(\widehat{P}^0_2)$-invariant neighborhood $W^0_2$ of $\widehat{P}^0_2$ in $W(P_2,\widetilde{P}_2)$. The orbifold transition function is the composition (after possibly shrinking $W^0_1$ and $W^0_2$ further so that $h$ lands in $\widetilde W^0_2$ and is surjective): $$W^0_1\xrightarrow{(t_{\Gamma(\widetilde{P}^0_1),\Gamma(\widetilde{P}_1)}|_{\widetilde W^0_1})^{-1}}{} \widetilde W^0_1\xrightarrow[\cong]{\quad h_Q\quad} \widetilde W^0_2\xrightarrow{t_{\Gamma(\widetilde{P}^0_1),\Gamma(\widetilde{P}_1)}|_{\widetilde W^0_2}}{} W^0_2,$$ which is a biholomorphism by construction.
    
    It remains to show that this biholomorphism is equivariant with respect to the action of $\Stab_{\Isot(P_1)}(\widehat{P}^0_1)\cong\Isot(P^0)$ on $W^0_1$ and the action of $\Stab_{\Isot(P_2)}(\widehat{P}^0_2)\cong\Isot(P^0)$ on $W^0_2.$ The identifications $$\Isot(P_1)\cong\frac{\Stab_G(\widetilde{P}_1)}{\Tw(\Gamma(\widetilde{P}_1))}\hspace{.5in}\text{and}\hspace{.5in}\Isot(P^0)\cong\frac{\Stab_G(\widetilde{P}^0_1)}{\Tw(\Gamma(\widetilde{P}^0_1))}$$ show that the action of $\Stab_{\Isot(P_1)}(\widehat{P}^0_1)$ on $W^0_1$ descends from the natural action of $\frac{\Stab_G(\widetilde{P}^0_1)}{\Tw(\Gamma(\widetilde{P}^0_1))}$ on $Q_{\Gamma(\widetilde{P}^0_1)}$, and similarly the action of $\Stab_{\Isot(P_2)}(\widehat{P}^0_2)$ on $W^0_2$ descends from the natural action of $\frac{\Stab_G(\widetilde{P}^0_2)}{\Tw(\Gamma(\widetilde{P}^0_2))}$ on $Q_{\Gamma(\widetilde{P}^0_2)}$. Equivariance of the orbifold transition function therefore follows from the equalities $$\Tw(\Gamma(\widetilde{P}^0_2))=\Tw(\Gamma(\widetilde{P}^0_1)\cdot h)=h^{-1}\Tw(\Gamma(\widetilde{P}^0_1))h$$
    and $$\Stab_G(\widetilde{P}^0_2)=\Stab(\widetilde{P}^0_1\cdot h)=h^{-1}\Stab(\widetilde{P}^0_1)h.$$
    
    We have shown that the local complex orbifold structures $W(P,\widetilde{P})/\Isot(P)\into U(G)/G$ are locally compatible, so they patch together, giving a global complex orbifold structure on $U(G)/G.$ The map $U(G)/G\to\Mbar_{g,n}$ is an \'etale morphism of complex orbifolds by construction. This proves \eqref{it:U1G}.
\end{proof}

\section{Complex structures on \texorpdfstring{$\Mbargnh$}{MHgn} and \texorpdfstring{$\Mbargnmh$}{hT(Vgn)}}\label{sec:ComplexStructures}
In Sections \ref{sec:Mghbar} and \ref{sec:Mgmhbar}, we endowed $\Mbargnh$ and $\Mbargnmh$ with topologies, inherited from the topology of $\Tbar(\Sgn)$. In this section, we apply Lemma \ref{lem:ComplexOrbifold} to endow them with complex structures.

\begin{thm}\label{thm:MarkedHandlebodiesComplexManifold}
    $\Mbargnmh$ has the natural structure of a complex manifold. With respect to this structure,
    \begin{enumerate}
    \item Each boundary stratum $\Mbargnmh_{(\tau,r)}\subseteq\Mbargnmh$ is a complex submanifold, and the restriction of the map $\Tbar(\Xgn)\to\Mbargnmh$ to any locally closed stratum is a holomorphic covering map over a stratum of $\Mbargnmh$,\label{it:MgmhStrataSuborbifolds}
    \item The submanifolds $\Mbargnmh_{(\tau,r)}$ meet at simple normal crossings,\label{it:SNCStrata}
    \item The natural map $\Mbargnmh\to\Mbar_{g,n}$ is an \'etale morphism of complex orbifolds.\label{it:MgmhMgEtale}
    \item The $\MCG(\Xgn)$-action on $\Tbar(\Xgn)$ induces an $\hMCG(\Xgn)$-action on $\Mbargnmh$ by biholomorphisms, with finite stabilizer groups.\label{it:OutFgActionBiholomorphisms}
    \end{enumerate}
\end{thm}

We have an analogous theorem for $\Mbargnh$, with a few key differences.

\begin{thm}\label{thm:UnmarkedHandlebodiesComplexOrbifold}
    $\Mbargnh$ has the natural structure of a complex orbifold. With respect to this structure,
    \begin{enumerate}
    \item Each stratum $\overline{\mathcal{MH}}_\tau$ is a complex suborbifold, and the strata meet at normal crossings,
    \item The natural map $\Mbargnmh\to\Mbargnh$ is a holomorphic orbifold covering map, and \label{it:MgmhMghCovering}
    \item The natural map $\Mbargnh\to\Mbar_{g,n}$ is an \'etale morphism of complex orbifolds. \label{it:MghMgEtale}
    \end{enumerate}
\end{thm}

We first prove the following, in order to show that Lemma \ref{lem:ComplexOrbifold} applies to our situation.
\begin{lem}\label{lem:StabilizerIsTwistSubgroup}
    Let $x\in\Tbar(\Xgn)\subseteq\Tbar(\Sgn),$ with contracted multidisk $\Delta$. Then $\Tw(\Delta)=
    \Tw(\Xgn)\cap\Stab_{\MCG(\Sgn)}(x)$.
\end{lem}
\begin{proof}
    Clearly $\Tw(\Delta)\subseteq\Tw(\Xgn)\cap\Stab_{\MCG(\Sgn)}(x).$ Let $h\in\Tw(\Xgn)\cap\Stab_{\MCG(\Sgn)}(x)$.
    Then $h$ induces an automorphism of the $n$-leafed dual graph $\tau_\Gamma\cong\tau_\Delta,$ where $\Gamma=\partial\Delta$.

    First, we show this automorphism is trivial. Note that $h\in\Stab_{\MCG(\Xgn)}(x)$, and that the map $\Stab_{\MCG(\Xgn)}(x)\to\Aut(\tau_\Delta)$ factors through $\hMCG(\Xgn)$. Since $h\in\Tw(\Xgn),$ it maps to the identity in $\Aut(\tau_\Gamma).$

    Since $h$ acts trivially on $\tau_\Gamma$, $h$ fixes each simple closed curve of $\Gamma$, and each connected component of $\Sgn\setminus\Gamma,$ up to isotopy. Therefore $h$ determines mapping classes $h^1,\ldots,h^k$ on the components $S^1,\ldots,S^k$ of $\Sgn\setminus\Gamma,$ with boundary components each contracted to a labeled point. By definition, we have $h^i\in\Stab_{\MCG(S^i)}(x^i),$ where $x^i$ is the complex structure on $S^i$ induced by $x$. Since the $\MCG(S^i)$-action on $\T(S^i)$ has finite stabilizers, the mapping classes $h^i$ all have finite order.

    Each $h^i$ has a lift $\tilde h^i$ in $\MCG(\Sgn)$ that acts by the identity on $\Sgn\setminus S^i$ and fixes $\Gamma$ pointwise. By construction, $\tilde h^i$ has finite order in $\MCG(\Sgn).$ 

    Note that the mapping classes $\tilde h^i$ all commute with each other so their product $\tilde h^1\cdots\tilde h^k$ has finite order. Also, they satisfy $h^{-1}\cdot(\tilde h^1\cdots\tilde h^k)\in\Tw(\Gamma)$, since this product is, up to isotopy, the identity away from a small neighborhood of $\Gamma$. Since $h^{-1}\cdot(\tilde h^1\cdots \tilde h^k)\in\Tw(\Gamma)\subseteq\Tw(\Xgn)$ and $h\in\Tw(\Xgn),$ we have $\tilde h^1\cdots\tilde h^k\in\Tw(\Xgn)$. But by Proposition \ref{prop:TwistSubgroupTorsionFree}, $\Tw(\Xgn)$ is torsion-free, so we must have $\tilde h^1\cdots \tilde h^k=1.$ It follows that $h\in\Tw(\Gamma)=\Tw(\Delta),$ as desired.
\end{proof}
\begin{proof}[Proof of Theorem \ref{thm:MarkedHandlebodiesComplexManifold}]
We apply Lemma \ref{lem:ComplexOrbifold} with $G=\Tw(\Xgn)$. Then $U(G)=\Tbar(\Xgn)$ by \cite[Thm. 5.6]{Hensel2020}\footnote{We actually only need $\Tbar(\Xgn)\subseteq\U(G)$, which is obvious from the definition of $U(G)$.}, so Lemma \ref{lem:ComplexOrbifold} endows $\Mbargnmh=\Tbar(\Xgn)/\Tw(\Xgn)$ with the structure of a complex orbifold, satisfying properties \eqref{it:MgmhStrataSuborbifolds} and \eqref{it:MgmhMgEtale}. By Lemma \ref{lem:StabilizerIsTwistSubgroup}, the isotropy groups of this orbifold are all trivial, i.e. $\Mbargnmh$ is actually a complex manifold. Property \eqref{it:SNCStrata} follows from Theorem \ref{thm:HubbardKoch}\eqref{it:QGammaStrataSNC} and the fact that $\CV_{g,n}^*$ is a simplicial complex. Property \eqref{it:OutFgActionBiholomorphisms} follows from Theorem \ref{thm:HubbardKoch}\eqref{it:QGammaInducedActionBiholomorphic}.
\end{proof}

\begin{proof}[Proof of Theorem \ref{thm:UnmarkedHandlebodiesComplexOrbifold}]
This follows from Theorem \ref{thm:MarkedHandlebodiesComplexManifold}, and the fact that the map $\Mbargnmh\to\Mbargnh$ is canonically identified with the quotient map by the properly discontinuous action of $\hMCG(\Xgn)$ on $\Mbargnmh$.
\end{proof}

\appendix

\section{Characterizations of the simplicial completion of Outer space}\label{app:SimplicialCompletion}
The space $\CV_{g,n}^*$ has at least three different characterizations:
\begin{itemize}
    \item A simplicial complex whose simplices correspond to marked graphs,
    \item A moduli space of marked \emph{metric} graphs, and
    \item The sphere complex of a doubled handlebody.
\end{itemize}
For convenience, in this appendix, we compare show the three characterizations and show their equivalence.

 \subsection{Closing the open  simplices of \texorpdfstring{$\CV_{g,n}$}{CVgn}} 
In \cite{CullerVogtmann1986} it was observed that Outer space $\CV_g$ decomposes as the disjoint union of open simplices $\sigma(\tau,r)$: 
$$ \CV_g= \bigsqcup_{(\tau,r)}\sigma(\tau,r),$$
where $(\tau,r)$ runs over equivalence classes of stable (unweighted) marked graphs of genus $g$ and the vertices of $\sigma(\tau,r)$ correspond to the edges of $\tau$. The simplicial completion $\CV_g^*$ of $\CV_g$ was then defined to be the union of corresponding closed simplices $\overline\sigma(\tau,r)$ modulo face relations, i.e.
let $$\widetilde{\CV}_g^*= \bigsqcup_{(\tau,r)}\overline\sigma(\tau,r);   \hbox{\,\,\, then \,\,\,}  \CV_g^*= \widetilde{\CV}_g^*/\sim$$ where  $\overline\sigma(\tau',r')$ is identified with a face of $\overline\sigma(\tau,g)$ if $\tau'$ is obtained from $\tau$ by collapsing each edge of a forest $\Phi\subset\tau$ and $r'$ is the composition of $r$ with the collapse map $c_\Phi\colon\tau\to\tau'$. The   vertices of this face correspond to the edges in $\tau\setminus\Phi$, so its  codimension is the number of edges in $\Phi$. 
The simplicial closure $\CV_g^*$ naturally has the structure of a $\Delta$-complex. It follows from \cite[Lem. 1.3]{SmillieVogtmann1987} that in this $\Delta$-complex, no two faces of a simplex are identified, i.e. $\CV_g^*$ is actually a simplicial complex.

The above construction applies without modification to $\CV_{g,n}$, i.e.
  $$\CV^*_{g,n}=\bigsqcup_{(\tau,r)}\overline\sigma{(\tau,r)}/\sim,$$ where $(\tau,r)$ runs over all marked stable genus-$g$  graphs with $k$ edges and $n$ labeled leaves,  $\overline\sigma{(\tau,r)}$ is a simplex with one vertex for each edge of $\tau$, and, if $ (\tau',r')$ is obtained from $(\tau,r)$ by collapsing a forest $\Phi$, then $\overline\sigma(\tau',r')$ is identified with the face of $\sigma(\tau,r)$ corresponding to the edges in $\tau\setminus\Phi$.
 
\subsection{Edge lengths allowed to be zero} The following is implicit in the definition of $\CV_{g,n}^*$ in Section \ref{sec:OuterSpaceBackground}:
\begin{prop}\label{prop:PointsOfSimplicialCompletion}
    The points of the simplicial complex $\CV_{g,n}^*$ are in natural bijection with the set of isomorphism classes of marked stable genus-$g$ $n$-leafed   metric graphs whose edge lengths sum to 1, with lengths allowed to be zero, modulo the equivalence relation generated by setting $(\tau_1,r_1)\sim(\tau_2,r_2)$ if $(\tau_1,r_1)$ is obtained from $(\tau_2,r_2)$ by contracting a length-zero non-loop edge in $\tau_2$. 
\end{prop}
\begin{proof}
Let $\widetilde{\mathscr{R}}$ denote the set of isomorphism classes of marked stable genus-$g$ $n$-leafed   metric graphs whose edge lengths sum to 1, with lengths allowed to be zero. There is a clear bijection $\widetilde{\mathscr{R}}\to\widetilde{\CV}_{g,n}^*$. It remains to show that under this identification, the equivalence relation on $\widetilde{\mathscr{R}}$ is identified with the face-gluing relation on $\widetilde{\CV}_{g,n}^*$. Indeed, if $(\tau_1,r_1)$ is obtained from $(\tau_2,r_2)$ by contracting a length-zero non-loop edge $e$, then the face in $\Pi_{(\tau_2,r_2)}$ corresponding to contracting $e$ consists of metrics on $\tau_2$ where $e$ has length zero. The equivalence relation generated by identifications of this form is clearly the equivalence relation on $\widetilde{\CV}_{g,n}^*$ corresponding to forest contractions.
\end{proof}

Using this definition of $\CV_{g,n}^*$, the quotient map $\CV_{g,n}^*\to\CV_{g,n}^*/A_{g,n}\cong\Link(\M_{g,n}^{\trop})$ can be described as follows. Given an element $(\tau,r)\in\CV_{g,n}^*$, we forget the marking $r$. We also contract all length-zero edges of $\tau$, and assign to each vertex $v$ a weight equal to the Betti number of the subgraph of $\tau$ that was contracted to $v$. The resulting stable weighted graph is the image of $(\tau,r)$ in $\Link(\M_{g,n}^{trop})$.

\subsection{Spheres in a doubled handlebody}
The first definition of $\CV_{g,n}^*$ for $n>0$ was given by Hatcher, who used a correspondence between sphere systems in a doubled handlebody and their dual graphs, which we now explain.  

Let $\Xgn$ be a genus-$g$ $n$-pointed handlebody, and let $W_{g,n}$ denote the ``double'' of $\Xgn$, i.e. the compact oriented 3-manifold obtained by gluing together two copies of $\X$, with opposite orientations, along their common boundary via the identity map. Since the $n$ labeled points of $\Xgn$ are along the boundary, $W_{g,n}$ is an $n$-pointed 3-manifold. In analogy with Sections \ref{sec:SurfacesBackground} and \ref{sec:HandlebodiesBackground}, we have the following definitions.

A \emph{sphere} in $W_{g,n}$ is an embedded copy of $S^2$ that does not contain any labeled points. A sphere $\omega$ necessarily divides $W_{g,n}$ into two connected components. We call $\omega$ \emph{essential} if each of these components either contains $\ge2$ labeled points, or is not simply connected. A \emph{multisphere} in $W_{g,n}$ is a finite union of disjoint, pairwise non-homotopic, essential spheres in $W_{g,n}$. A multisphere $\Omega$ is called \emph{pure} if every connected component of its complement is simply connected (this is called {\emph{simple} or \emph{complete} in the geometric group theory literature}). The \emph{dual graph} $\tau_\Omega$ of a multisphere $\Omega$ is defined similarly to the dual graph of a multidisk, as is the dual graph map $\dgm_{W_{g,n},\Omega}:W_{g,n}\setminus\{x_1,\ldots,x_n\}\to\tau_\Omega.$

The \emph{sphere complex} $\Sph(W_{g,n})$ is the infinite simplicial complex whose vertices are spheres in $W_{g,n}$, and whose $k$-simplices are multispheres with $k+1$ elements.

Fix an $n$-pointed homotopy equivalence $\rho\colon R_{g,n}\to V_{g,n}$.  This induces a map $f_\rho\colon\Sph(W_{g,n})\to\CV_{g,n}^*$ as follows. Given a multisphere $\Omega\subseteq W_{g,n}$, let $\Omega'\supseteq\Omega$ be a pure multisphere. The map takes the simplex corresponding to $\Omega$ to the simplex corresponding to the dual graph $\tau_{\Omega'},$ marked by the composition 
  $$R_{g,n}\xrightarrow{\rho} V_{g,n}\subset W_{g,n}\xrightarrow{\dgm_{W_{g,n},\Omega}} \tau_\Omega.$$  All of these maps are isomorphisms on $\pi_1$, so the composition is a homotopy equivalence (any map on graphs which induces an isomorphism on $\pi_1$ is a homotopy equivalence).

There is also a ``doubling" map from the disc complex $\mathcal D(V_{g,n})$ to the sphere complex $\Sph(W_{g,n})$, that sends a multidisk to its double, which is a multisphere.  We have the following results from 3-manifold theory.

\begin{lem}\label{lem:DoublingSurjective}
    The doubling map $\D(\Xgn)\to\Sph(W_{g,n})$ is surjective.
\end{lem}
\begin{proof} 
By \cite[Thm. 1.3]{Scharlemann2024}, any embedded $2$-sphere in a punctured doubled handlebody is isotopic to one that intersects the handlebody surface in a single simple closed curve. Since there is only one isotopy class of disks in a handlebody with a given boundary curve \cite[Lem. 2.3]{Hensel2020}, the sphere is the double of this disk.
\end{proof}
\begin{rem}
    Scharlemann's theorem is much more general, applying to any compact 3-manifold with a fixed Heegard splitting.
\end{rem}

\begin{lem}\label{lem:DoublingTwistInvariant}
    Let $\omega$ and $\omega'$ be essential spheres in $W_{g,n}$. Then applying the Dehn twist around $\omega$ to $\omega'$ results in a sphere isotopic to $\omega'$.
\end{lem}
  \begin{proof}     This follows from Laudenbach's theorem that homotopic spheres in a doubled handlebody are isotopic.  See \cite{Hatcher1995}  for a detailed explanation of exactly how it follows.  The relevant papers of Laudenbach are \cite{Laudenbach1973} for a single sphere and \cite{Laudenbach1974} for a multisphere.
\end{proof}

\begin{prop}\label{prop:DiskComplexSphereComplex}
The doubling map induces an isomorphism $\mathcal D(V_{g,n})/\Tw(V_{g,n}) \to \Sph(W_{g,n}).$
\end{prop}

\begin{proof} 
    In Proposition \ref{prop:MarkedBijection}, we defined a map $\D(\Xgn)\to\CV_{g,n}^*$, which is an isomorphism after quotienting by $\Tw(\Xgn)$. This map clearly factors through the doubling map $\D(\Xgn)\to\Sph(W_{g,n})$, identifying the map $\Sph(W_{g,n})\to \CV_{g,n}^*$ with the map $f_\rho$ defined above. The doubling map is twist-invariant by Lemma \ref{lem:DoublingTwistInvariant}, so we have a diagram of simplicial complexes:
    $$\begin{tikzcd}
        \D(\Xgn)/\Tw(\Xgn)\arrow[r]\arrow[rr,bend right, "\cong"]&\Sph(W_{g,n})\arrow[r, "f_\rho"]&\CV_{g,n}^*.
    \end{tikzcd}$$
    The doubling map is surjective by Lemma \ref{lem:DoublingSurjective}, which implies that all maps above are isomorphisms. Thus $\mathcal D(V)/\Tw(V)\cong\Sph(W_{g,n})\cong\CV_{g,n}^*$.
\end{proof}
\begin{rem}
    The isomorphism $\Sph(W_{g,n})\xrightarrow{\cong}{}\CV_{g,n}^*$ is implicit in \cite{Hatcher1995}, where Hatcher takes $\Sph(W_{g,n})$ as the definition of the simplicial completion of $\CV_{g,n},$ and identifies $\CV_{g,n}$ with $\Sph^{\pure}(W_{g,n})$.
\end{rem}

\begin{cor}\label{cor:NotMultidisk}
    If $\Delta$ and $\Delta'$ are multidisks in $\Xgn$ that are in the same $\Tw(\Xgn)$-orbit, then $\Delta\cup\Delta'$ is not a multidisk.
\end{cor}
\begin{proof}
    If $\Delta\cup\Delta'$ were a multidisk, then the quotient map $\D(\Xgn)\to\D(\Xgn)/\Tw(\Xgn)$ would identify two faces of the corresponding cell in $\D(\Xgn),$ contradicting the fact that the quotient is a simplicial complex.
\end{proof}

\bibliography{Handlebodies}
\bibliographystyle{amsalpha}
\end{document}